\documentclass[12pt]{amsart}

\usepackage{latexsym, amssymb, amscd, amsfonts}

\usepackage{microtype}
\usepackage{pstricks}
\usepackage{pst-3d}
\usepackage{multido}
\usepackage{pst-node}
\usepackage{pst-arrow}
\usepackage{pst-plot} 
\usepackage[colorlinks=true,urlcolor=black,citecolor=black,linkcolor=black]{hyperref}
\usepackage{graphicx}
\usepackage[all]{xypic}
\usepackage{array}
\usepackage{ifthen}

\newboolean{bourbaki}
\setboolean{bourbaki}{true}

\newboolean{draft}
\setboolean{draft}{false}

\newboolean{shortnote}
\setboolean{shortnote}{false}

\ifthenelse{\boolean{shortnote}}{%
\newcommand{\point}{\vspace{3mm}\par\refstepcounter{subsection}\noindent{\bf \thesubsection.} }

\newcommand{\bpoint}[1]{\vspace{3mm}\par\refstepcounter{subsection}\noindent{\bf \thesubsection.} 
  {\bf #1.} }

\renewenvironment{equation}{\medskip\noindent\refstepcounter{subsection}\makebox[0pt][l]{({\bf\thesubsection})}\begin{minipage}{\textwidth}$$}{$$\end{minipage}\medskip\noindent}
} 
{
\newcommand{\point}{\vspace{3mm}\par\refstepcounter{subsection}\noindent{\bf \thesubsection.} }

\newcommand{\bpoint}[1]{\vspace{3mm}\par\refstepcounter{subsection}\noindent{\bf \thesubsection.} 
  {\bf #1.} }

\renewenvironment{equation}{\medskip\noindent\refstepcounter{subsubsection}\makebox[0pt][l]{({\bf\thesubsubsection})}\begin{minipage}{\textwidth}$$}{$$\end{minipage}\medskip\noindent}
\newcommand{\neweqnum}{\refstepcounter{subsubsection}({\bf\thesubsubsection})}
\newenvironment{equation-w-subscript}[1]{\medskip\noindent\refstepcounter{subsubsection}\makebox[0pt][l]{$\mbox{({\bf\thesubsubsection})}_{#1}$}\begin{minipage}{\textwidth}$$}{$$\end{minipage}\medskip\noindent}
}

\newcommand{\Fig}{\refstepcounter{subsubsection}{\sc Figure \thesubsubsection.}} 

\newcommand{\bpf}{\noindent {\em Proof.  }}
\newcommand{\epf}{\qed \vspace{+10pt}}

\catcode`\@=11 
\def\addtopunct#1{\expandafter\let\csname punct@\meaning#1\endcsname\let}
\addtopunct{.}    \addtopunct{,}    \addtopunct{?}
\addtopunct{!}    \addtopunct{;}    \addtopunct{:}

\let\seveendformula\]
\def\PunctAndEndFormula #1{#1\seveendformula}
\def\]{\futurelet\punctlet\checkpunct@i}
\def\checkpunct@i{\expandafter\ifx\csname punct@\meaning\punctlet\endcsname\let  
       \expandafter\PunctAndEndFormula 
       \else \expandafter\seveendformula\fi}
\catcode`\@=12 

\makeatletter
\newcommand*{\da@rightarrow}{\mathchar"0\hexnumber@\symAMSa 4B }
\newcommand*{\da@leftarrow}{\mathchar"0\hexnumber@\symAMSa 4C }
\newcommand*{\xdashrightarrow}[2][]{%
  \mathrel{%
    \mathpalette{\da@xarrow{#1}{#2}{}\da@rightarrow{\,}{}}{}%
   }%
}
\newcommand{\xdashleftarrow}[2][]{%
  \mathrel{%
    \mathpalette{\da@xarrow{#1}{#2}\da@leftarrow{}{}{\,}}{}%
  }%
}
\newcommand*{\da@xarrow}[7]{%
  \sbox0{$\ifx#7\scriptstyle\scriptscriptstyle\else\scriptstyle\fi#5#1#6\m@th$}%
  \sbox2{$\ifx#7\scriptstyle\scriptscriptstyle\else\scriptstyle\fi#5#2#6\m@th$}%
  \sbox4{$#7\dabar@\m@th$}%
  \dimen@=\wd0 %
  \ifdim\wd2 >\dimen@
    \dimen@=\wd2 %
  \fi
  \count@=2 %
  \def\da@bars{\dabar@\dabar@}%
  \@whiledim\count@\wd4<\dimen@\do{%
    \advance\count@\@ne
    \expandafter\def\expandafter\da@bars\expandafter{%
      \da@bars
      \dabar@ 
    }%
  }%
  \mathrel{#3}%
  \mathrel{%
    \mathop{\da@bars}\limits
    \ifx\\#1\\%
    \else
      _{\copy0}%
    \fi
    \ifx\\#2\\%
    \else
      ^{\copy2}%
    \fi
  }%
  \mathrel{#4}%
}
\makeatother

\newcommand{\xmapsto}{\mapstochar\relbar\joinrel\xrightarrow}

\newcommand{\xdownarrow}[1]{%
  {\left\downarrow\vbox to #1{}\right.\kern-\nulldelimiterspace}
}

\makeatletter
\newcommand*{\Relbarfill@}{\arrowfill@\Relbar\Relbar\Relbar}
\newcommand*{\xeq}[2][]{\ext@arrow 0055\Relbarfill@{#1}{#2}}
\makeatother

\newgray{lgray}{0.9}
\newgray{vlgray}{0.95}
\newgray{midgray}{0.65}
\newgray{darkgray}{0.3}

\newcommand{\Rot}{%
[ exch dup cos dup 3 -1 roll sin dup -1 mul 4 -1 roll 0 0 ] transform
}

\newcommand{\SetTT}[6]{%
\def \TT {[ #1 #2 #3 #4 6 -1 roll dup #5 mul exch #6 mul ] transform }
}
\SetTT{-0.5}{-0.5}{1}{0}{0}{1}

\newcommand{\Hom}{\operatorname{Hom}}

\renewcommand{\part}{\operatorname{part}}

\newcommand{\st}{\,\,|\,\,}

\newcommand{\Pic}{\operatorname{Pic}}

\newcommand{\Fsh}{\mathcal{F}}
\newcommand{\Gsh}{\mathcal{G}}
\newcommand{\Hsh}{\mathcal{H}}
\newcommand{\Ish}{\mathcal{I}}
\newcommand{\Msh}{\mathcal{M}}
\newcommand{\Osh}{\mathcal{O}}

\renewcommand{\geq}{\geqslant}
\renewcommand{\leq}{\leqslant}

\newcommand{\CComp}{\check{C}}           
\newcommand{\CComplex}{\CComp^{\smallbullet}} 
\newcommand{\HP}{\check{H}_{P}}
\newcommand{\HQ}{\check{H}_{Q}}
\newcommand{\ep}{\epsilon} 
\renewcommand{\emptyset}{\varnothing} 
\renewcommand{\phi}{\varphi}
\newcommand{\Ab}{\operatorname{Ab}}  
\newcommand{\ZZsheaf}{\underline{\ZZ}} 
\newcommand{\und}{\underline}
\renewcommand{\div}{\operatorname{div}} 
\newcommand{\Spec}{\operatorname{Spec}} 
\newcommand{\cdim}{\operatorname{cdim}} 
\newcommand{\codim}{\operatorname{codim}} 
\newcommand{\Eff}{\operatorname{Eff}} 
\newcommand{\sg}{\mbox{\tiny sg}}   
\newcommand{\charx}{\mathfrak{X}}   

\makeatletter
\newcommand{\smallbullet}{} 
\DeclareRobustCommand\smallbullet{%
  \mathord{\mathpalette\smallbullet@{0.65}}
}
\newcommand{\smallbullet@}[2]{%
  \vcenter{\hbox{\scalebox{#2}{$\m@th#1\bullet$}}}%
}
\makeatother

\newcommand{\opp}{\operatornamewithlimits{\oplus}}

\ifthenelse{\boolean{bourbaki}}
{
\renewcommand{\AA}{\mathbf{A}} 
\newcommand{\CC}{\mathbf{C}} 
\newcommand{\FF}{\mathbf{F}} 
\newcommand{\NN}{\mathbf{N}} 
\newcommand{\PP}{\mathbf{P}} 
\newcommand{\RR}{\mathbf{R}} 
\newcommand{\ZZ}{\mathbf{Z}} 

\DeclareSymbolFont{frenchmath}{OT1}{cmr}{m}{n} 

\DeclareMathSymbol{A}{\mathalpha}{frenchmath}{`A}  
\DeclareMathSymbol{B}{\mathalpha}{frenchmath}{`B}  
\DeclareMathSymbol{C}{\mathalpha}{frenchmath}{`C}
\DeclareMathSymbol{D}{\mathalpha}{frenchmath}{`D}
\DeclareMathSymbol{E}{\mathalpha}{frenchmath}{`E}
\DeclareMathSymbol{F}{\mathalpha}{frenchmath}{`F}
\DeclareMathSymbol{G}{\mathalpha}{frenchmath}{`G}
\DeclareMathSymbol{H}{\mathalpha}{frenchmath}{`H}
\DeclareMathSymbol{I}{\mathalpha}{frenchmath}{`I}
\DeclareMathSymbol{J}{\mathalpha}{frenchmath}{`J}
\DeclareMathSymbol{K}{\mathalpha}{frenchmath}{`K}
\DeclareMathSymbol{L}{\mathalpha}{frenchmath}{`L}
\DeclareMathSymbol{M}{\mathalpha}{frenchmath}{`M}
\DeclareMathSymbol{N}{\mathalpha}{frenchmath}{`N}
\DeclareMathSymbol{O}{\mathalpha}{frenchmath}{`O}
\DeclareMathSymbol{P}{\mathalpha}{frenchmath}{`P}
\DeclareMathSymbol{Q}{\mathalpha}{frenchmath}{`Q}
\DeclareMathSymbol{R}{\mathalpha}{frenchmath}{`R}
\DeclareMathSymbol{S}{\mathalpha}{frenchmath}{`S}
\DeclareMathSymbol{T}{\mathalpha}{frenchmath}{`T}
\DeclareMathSymbol{U}{\mathalpha}{frenchmath}{`U}
\DeclareMathSymbol{V}{\mathalpha}{frenchmath}{`V}
\DeclareMathSymbol{W}{\mathalpha}{frenchmath}{`W}
\DeclareMathSymbol{X}{\mathalpha}{frenchmath}{`X}
\DeclareMathSymbol{Y}{\mathalpha}{frenchmath}{`Y}
\DeclareMathSymbol{Z}{\mathalpha}{frenchmath}{`Z}
} 
{ 
\renewcommand{\AA}{\mathbb{A}} 
\newcommand{\CC}{\mathbb{C}} 
\newcommand{\FF}{\mathbb{F}} 
\newcommand{\NN}{\mathbb{N}} 
\newcommand{\PP}{\mathbb{P}} 
\newcommand{\RR}{\mathbb{R}} 
\newcommand{\ZZ}{\mathbb{Z}} 
}

\ifthenelse{\boolean{draft}}
{
\newcommand{\remind}[1]{{\bf[#1]}}
\newcommand{\lremind}[1]{{\bf[label:  #1]}}
\newcommand{\bremind}[1]{{\bf[label:  #1]}}

}
{
\newcommand{\remind}[1]{{}}
\newcommand{\lremind}[1]{{}}
\newcommand{\bremind}[1]{{}}

}

\DeclareMathAlphabet{\matheuler}{U}{zeur}{m}{n}
\DeclareMathAlphabet{\matholdstyle}{OT1}{pplj}{m}{n}
\newcommand{\euler}[1]{\ensuremath{\matheuler{#1}}}  
\newcommand{\Euler}{\euler} 

\newcounter{enumcount}
\newcommand{\PauseEnumerate}{\end{enumerate}\setcounter{enumcount}{\arabic{enumi}}}
\newcommand{\ResumeEnumerate}{\begin{enumerate}\setcounter{enumi}{\theenumcount}}

\newcommand{\AlphaList}{\renewcommand{\labelenumi}{{({\em\alph{enumi}})}}}

\newcommand{\RomanList}{\renewcommand{\labelenumi}{{({\em\roman{enumi}})}}}
\newcommand{\EulerList}{\renewcommand{\labelenumi}{{($\matheuler{\arabic{enumi}}$)}}}
\AlphaList  

\makeatletter
\let\@wraptoccontribs\wraptoccontribs
\makeatother  

\usepackage{tikz-cd}

\newcommand{\Neg}{\operatorname{Neg}}
\renewcommand{\neg}{\operatorname{neg}}
\usepackage[bsmallmatrix]{mathtools}
\usetikzlibrary{calc}
\usepackage[OT1]{fontenc}

\begin{document}
\pagestyle{plain} \title{\large Reduced \v{C}ech complexes and computing higher direct images under toric maps}
\author{Mike Roth}
\address{Dept.\ of Mathematics and Statistics, Queens University, Kingston,
Ontario, Canada} 
\email{mike.roth@queensu.ca} 
\thanks{Research partially supported by an NSERC grant}
\contrib[Appendices by]{Sasha Zotine}
\address{Dept.\ of Mathematics and Statistics, McMaster University, Hamilton, 
Ontario, Canada} 
\email{zotinea@mcmaster.ca}


\maketitle

\section{Introduction}

\point  Let $X$ be a topological space, $\{U_{\alpha}\}_{\alpha\in I}$ an open cover of $X$, and $\Fsh$ a sheaf
of abelian groups on $X$.  Associated to this data is a total \v{C}ech complex, the complex 

\begin{equation}
\label{eqn:total-cech} 
\raisebox{-0.15cm}{\rule{1.3cm}{0cm}
\mbox{$\displaystyle
0 \longrightarrow \prod_{\alpha\in I} \Gamma(U_{\alpha},\Fsh) \longrightarrow \prod_{\alpha,\beta\in I} 
\Gamma(U_\alpha\cap U_{\beta},\Fsh) \longrightarrow \prod_{\alpha,\beta,\gamma\in I} 
\Gamma(U_{\alpha}\cap U_{\beta}\cap U_{\gamma},\Fsh) \longrightarrow \cdots .$}}
\end{equation}

If the higher cohomology groups of $\Fsh$ are zero on all finite intersections of elements of the cover 
then by the well-known theorem of Leray (see e.g., \cite[p.\ 213, Corollaire]{G}) 
the complex \eqref{eqn:total-cech} computes the cohomology groups $H^{i}(X,\Fsh)$. 

For instance, if $X$ is an $n$-dimensional manifold then one may use a `good cover' (one where all finite
intersections are homeomorphic to $\RR^n$) to compute the cohomology of any locally constant sheaf on $X$,
or if $X$ is a scheme, one may use a cover where all finite intersections are affines to compute the cohomology
of any quasi-coherent sheaf on $X$. 

More generally, even without the hypothesis that the higher cohomology on the intersections are zero,  
there is a spectral sequence linking the \v{C}ech complex above and the cohomology of $\Fsh$ on $X$. 

\bpoint{Smaller \v{C}ech complexes}  
\label{sec:Oda-construction} 
In the literature, for certain kinds of toric varieties, there is a construction of a type of \v{C}ech complex
which uses fewer intersections\footnote{ The author believes that this is originally due to Oda, but does not
know a reference.}.

For instance, if $X=\PP^1\times\PP^1$, starting with the four torus-stable affines corresponding to the
four torus fixed points, one may think of these as being associated to the vertices of a square, and then form
a complex where, in degree zero one uses the sections of the sheaf on these open sets, in degree one the 
sections over intersections of two of the open sets only if the corresponding vertices are connected by an edge, and 
in degree two the sections over the intersection of all four of the original open sets 
(this example is given in more detail in \S\ref{sec:first-example-of-reduced}).  
It turns out that this complex again computes the cohomology of coherent sheaves on $X$.

In general a similar construction applies to any projective simplicial toric variety (that is, a projective toric
variety where the fan is simplicial), with the combinatorics of the complex governed by the polytope $P$ associated
to $X$, e.g., under any embedding by an ample line bundle.

One argument that the resulting complex computes the cohomology of
line bundles on $X$ is the following\footnote{ For simplicity we will just discuss the case of line bundles, instead
of arbitrary coherent sheaves.}. 
Given such a variety $X$,  by \cite[Theorem 2.1]{C} $X$ is a GIT quotient 
$U/\!\!/(\CC^{*})^{\rho}$ where $\rho$ is the Picard rank of $X$, and $U\subseteq \AA^{s}$ is an open subset, 
with $s$ being the number of rays in the fan of $X$.

Letting $\pi\colon U\longrightarrow X$ be the quotient map, the $(\CC^{*})^{\rho}$ action on the fibres gives
a decomposition $\pi_{*}\Osh_{U} = \oplus_{L\in\Pic(X)} L$.  Since $\pi$ is affine, to compute the cohomology
$H^{i}(X,L)$ of a given line bundle $L$, it then suffices to find the appropriate graded piece of $H^{i}(U,\Osh_{U})$.

Let $Z$ be the complement of $U$ in $\AA^{s}$.  
Then $Z$ is the zeros of a monomial ideal, whose combinatorics are directly connected with $P$.  
One may then use the combinatorics of the monomial ideal to show
that the smaller \v{C}ech complex computes the cohomology groups $H^{i}(U,\Osh_{U})$, or equivalently, 
local cohomology with support in $Z$. For instance, this is guaranteed by \cite[Theorem 6.2]{Mi2}.

In \cite{Mi2} this complex is called the canonical \v{C}ech complex associated to $X$.   In this paper we will
use the term {\em reduced \v{C}ech complexes}, since these complexes require fewer intersections than the standard
\v{C}ech complex, and since we will give a general method of constructing such complexes which is rather arbitrary,
and so cannot be considered to be canonical.

\bpoint{Objectives} This paper has three main goals : 

\EulerList
\begin{enumerate}
\item To give an axiomatic formulation of the construction of a ``reduced \v{C}ech complex''. 
\PauseEnumerate

This method is different than the one above and applies to arbitrary topological spaces. In particular, it does not
depend on describing $X$ as a GIT quotient, nor on using results on monomial ideals.  

In this approach the combinatorics of the reduced complex are again controlled by a polytope, or more generally a cell 
complex, $P$.  For each $k$-face (or $k$-cell) $\sigma$ of $P$ we must give an open set $U_{\sigma}\subseteq X$.  
The first two axioms require that this association be order reversing, and that the collection of open sets cover $X$. 
Given a sheaf of abelian groups $\Fsh$ on $X$, we may then define a reduced \v{C}ech complex associated 
to $\Fsh$ and the data of the cover. 

There are two additional axioms relating the topology of $P$ and the topology of $X$ as seen by the cover 
(see \S\ref{sec:compatibility-conditions}).  If those conditions hold then, as in the case of a usual \v{C}ech complex,
the reduced \v{C}ech complex computes the cohomology of all sheaves whose higher cohomology vanishes on the open
sets of the cover (Theorem \ref{thm:reduced-cech-computes-cohomology}).
As before, even without the hypothesis of the vanishing of the higher cohomology, 
there is a spectral sequence linking the reduced \v{C}ech complex and
the cohomology of a sheaf on $X$ (Theorem \ref{thm:reduced-cech-spectral-sequence}).

\ResumeEnumerate
\item To extend the construction of \S\ref{sec:Oda-construction} to arbitrary proper toric varieties, and more 
generally, to semi-proper toric varieties.
\PauseEnumerate

When $X$ is a proper toric variety (including non-simplicial toric varieties),  we give a construction of 
a reduced \v{C}ech complex where each of the open sets are torus-stable affines, and which computes the cohomology
of quasi-coherent sheaves on $X$.  For the cell complex $P$ we use the non-negative part of $X$
(see \S\ref{sec:cell-complex-toric-variety}), denoted $X_{\geq 0}$. 
When $X$ is projective, $P$ is the moment polytope of $X$, or equivalently, the polytope associated to any ample
line bundle on $X$. 

In the case of a projective simplicial toric variety the reduced complex produced by this method is the same as the one
given by the method of \S\ref{sec:Oda-construction}.

More generally, for any semi-proper toric variety $X$ (that is, a toric variety such that the canonical map
$X\longrightarrow \Spec(\Gamma(X,\Osh_{X}))$ is proper) we construct a reduced \v{C}ech complex on $X$ with the property
that each open is a torus-stable affine, and such that $\dim(P)$ is the cohomological dimension of $X$ 
(Theorem \ref{thm:construction-of-toric-reduced-cover}).  Here $P$ is not in general a polytope, for instance it
may be made of several polytopes attached in a polyhedral complex. 

\ResumeEnumerate
\item To give an algorithm to compute the higher direct images of line bundles relative to a toric fibration  
$\phi\colon X\longrightarrow Y$ between smooth proper toric varieties. 
\PauseEnumerate

In the next sections we discuss the meaning of this last goal, and the solution given in this paper.

\point 
\label{sec:issue-finding-module}
Let $\phi\colon X\longrightarrow Y$ be a toric morphism between smooth proper toric varieties such that
$\phi_{*}\Osh_{X}=\Osh_{Y}$, and let $L$ be a line bundle on $X$.  The goal above
is to develop an algorithm to calculate the higher direct images sheaves $R^{i}\phi_{*}L$, for any $i\geq 0$. 
To the authors' knowledge, the only previous results of this type are \cite{ES} and \cite{EES}, which compute
the higher direct images sheaves for the case of projective space and the product of projective spaces respectively, as
modules over a ground ring $A$.  Those papers consider the case of an arbitrary coherent sheaf and not just 
line bundles.  However, we hope in the future to be able to use the method of this paper and a resolution of the
sheaf by line bundles to handle that case.

One issue in getting a computer to calculate a sheaf is to be able to describe it using only a finite 
amount of data.  On projective space there is a familiar solution.  If $\Fsh$ is a coherent sheaf on $\PP^n$, then 
$$\Gamma_{*}(\Fsh):= \opp_{m\geq 0} H^{0}(\PP^n, \Fsh\otimes \Osh_{\PP^n}(m))$$
is a finitely generated $\ZZ$-graded module over $S=\CC[y_0, y_1,\ldots, y_n]$.  
Conversely, there is an inverse sheafification
process which, given a finitely generated graded $S$-module $M$, produces a coherent sheaf $\widetilde{M}$ on $\PP^n$. 

For a smooth projective toric variety $Y$ one similarly has the Cox ring $S=\CC[y_0,\ldots , y_{s-1}]$, 
where $s$ is the number of torus fixed irreducible divisors in $Y$. 
This ring is graded by the group $\Pic(Y)$ 
with each $y_i$ having the degree of the associated divisor (i.e., the line bundle associated to that divisor).
Given a finitely generated $\Pic(Y)$-graded module $M$ over $S$, there is a sheafification process which again gives a 
coherent sheaf $\widetilde{M}$ on $Y$, \cite[\S3]{C}.

But, given a coherent sheaf $\Fsh$ on $Y$ it is not always obvious how to find a corresponding $M$.   
For instance, when $\Fsh = R^{i} \phi_{*} L$ is a higher direct image (even when $i = 0$), often the direct sum
$$ \Gamma_{*}(\Fsh):=\opp_{\Msh\in \Eff(Y)} H^0(Y, \Fsh \otimes \Msh) $$
while an $S$-module, is not a finitely generated one.  
Here is a simple example of this behaviour. 

\bpoint{Example} 
\label{ex:not-fg}
Let $Y=\PP\bigl(\Osh_{\PP^1}\oplus \Osh_{\PP^1}(1)\bigr)$ be the first Hirzebruch surface, which we take over
the field $\CC$, although the choice of field is irrelevant to the example.
Let $B\subset Y$ be the unique curve of self-intersection $(-1)$, and $p\in B$ one of the two torus fixed points.  
We take $\Fsh=\Osh_{p}$ (i.e, the structure sheaf of the point, pushed forward to $X$). 

Here the torus-fixed irreducible divisors are $B$,  two fibre-class curves $F_1$ and $F_2$, and a curve $H$,
linearly equivalent to $B+F_i$, of self intersection $1$. Let the variables $y_0$,\ldots, $y_3$ correspond to
these divisors, in the order $F_1$, $B$, $F_2$, and $H$.  
Taking $B$ and the fibre classes $F_i$ as a basis for $\Pic(Y)\cong \ZZ^2$, these
variables have bidegrees $(0,1)$, $(1,0)$, $(0,1)$, and $(1,1)$ respectively. 
In this basis a class $(a,b)\in \Pic(Y)$ is effective if and only if $a$, $b\geq 0$.

Since $\Fsh$ is the structure sheaf of a point, for any line bundle $\Msh$ on $Y$, $H^{0}(Y,\Fsh\otimes \Msh)$ is a 
$1$-dimensional vector space over $\CC$.  Thus 
$$ \Gamma_{*}(\Fsh)=\opp_{\Msh\in \Eff(Y)} H^0(Y, \Fsh \otimes \Msh) = \opp_{(a,b)\in \NN^2} \CC_{a,b},$$
where $\CC_{a,b}$ denotes the field $\CC$, placed in bidegree $(a,b)$.  

The only monomial in $S=\CC[y_0,\ldots, y_3]$ of bidegree $(k,0)$ is $y_1^{k}$.  
Since $y_1$ is zero on $B$, and since $p\in B$, 
multiplication by $y_1$ is zero on $H^0(Y,\Fsh\otimes \Msh)$ for any $\Msh$.  In particular, for each $a\geq 0$, 
the graded piece $\CC_{a,0}$ is not the image of any other graded piece under the action of $S$.   Thus, 
the $\CC_{a,0}$, for all $a\geq 0$ must be included in any $S$-module generating set for $\Gamma_{*}(\Fsh)$, 
and $\Gamma_{*}(\Fsh)$ is therefore not finitely generated. 

This example is an example of a higher direct image.  Let be $X$ be the blowup of $Y$ at $p$, 
$\phi\colon X\longrightarrow Y$ the blowdown map, $E$ the exceptional divisor, and $L=\Osh_{X}(-2E)$.  Then (since
$p$ is a torus-fixed point on $Y$), $X$ is a toric variety, $\phi$ a toric morphism, and 
one can check that $R^{1}\phi_{*}L = \Osh_{p} = \Fsh$.

\point 
\label{sec:Temp-def-of-Hi}
Returning to the general situation, one can similarly look at the $S$-module 
$$\Hsh^{i}(L):=\opp_{\Msh\in \Eff(Y)} H^{i}(X, L \otimes \phi^{*}\Msh),$$
which has a natural edge morphism to $\Gamma_{*}(R^{i}\phi_{*}L)$. 
In general $\Hsh^{i}(L)$ is not a finitely generated $S$-module either, with 
Example \ref{ex:not-fg} again providing an illustration. 

Our solution, as explained below, may be viewed as finding a finitely generated $S$-submodule of $\Hsh^{i}(L)$ 
which still sheafifies to the correct higher direct image.

\point
\label{sec:algorithm-description}
Let $T_K$ be the kernel of the surjective map $T_X\longrightarrow T_Y$ induced by $\phi$, where $T_{X}$ and $T_{Y}$
denote the tori of $X$ and $Y$ respectively.  Given a line bundle
$L$ on $X$, one can lift the $T_X$-action on $X$ to an action on $L$, and thus $T_K$ also acts on $L$. 
Fixing an $i\geq 0$, by the equivariance of $\phi$, the higher direct image sheaf $R^{i}\phi_{*}L$ has a $T_K$ action,
with $T_K$ acting trivially on $Y$, and thus $R^{i}\phi_{*}L$ splits as a direct sum of $T_K$ eigensheaves~:

\begin{equation}
\label{eqn:TK-splitting}
 R^{i}\phi_{*}L =\opp_{\mu\in C(L,i)} \left(R^{i}\phi_{*}L\right)_{\mu},
\end{equation}

with $\mu$ running over a finite subset $C(L,i)$ of the $T_K$-characters. 

Given $L$ on $X$ and an $i\geq 0$, our algorithm does the following~:

\RomanList
\begin{enumerate}
\item Finds the set $C(L,i)$ of $T_K$-characters $\mu$ such that $\left(R^{i}\phi_{*}L\right)_{\mu}\neq 0$;
\item For each $\mu\in C(L,i)$, produces a divisor $D_{\mu}$ on $X$, and a divisor $E_{\mu}$ on $Y$. 
\item For each $D_{\mu}$, constructs a finitely generated $S$-module $M^{i}(D_{\mu})$.
\end{enumerate}

The construction is such that each $M^{i}(D_{\mu})$, with grading shifted by $E_{\mu}$, is naturally an $S$-submodule 
of $\Hsh^{i}(L)_{\mu}$.  Our theorem is then that for each $\mu$ the resulting natural map 
$$\widetilde{M}^{i}(D_{\mu})\otimes \Osh_{Y}(E_{\mu})\longrightarrow \left(R^{i}\phi_{*}L\right)_{\mu}$$
is an isomorphism.   Summing over the finitely many $\mu\in C(L,i)$, we thus obtain an isomorphism of sheaves 
$$\opp_{\mu\in C(L,i)} \left(\widetilde{M}^{i}(D_{\mu})\otimes \Osh_{Y}(E_{\mu})\right)
\stackrel{\sim}{\longrightarrow} R^{i}\phi_{*}L.$$

The module in ({\em iii}) is the $i$-th cohomology of a complex of monomial ideals of $S$.  
This complex is constructed using the machinery of reduced \v{C}ech complexes developed in the earlier part of the 
paper.  In particular, to make the computation more efficient, we take advantage of the fact that if
$\phi\colon X\longrightarrow Y$ is a proper toric morphism, and $V\subseteq Y$ a torus-stable open affine,
then $\phi^{-1}(V)$ is a semi-proper toric variety, and thus we may apply the construction of
Theorem \ref{thm:construction-of-toric-reduced-cover}.

After glueing these local constructions on the base $Y$, and passing to $S$-modules, 
the resulting complex of monomial ideals of $S$ is governed by $P=X_{\geq 0}$ the polytope (or, in general, 
cell complex) associated to the toric variety $X$.  
By way of illustration of the output of the algorithm, although without supplying the intermediate calculations, 
we present the following example.  
(For examples of computations with the intermediate steps included see \S\ref{sec:appendix-examples}.)

\bpoint{Example} 
\label{sec:X-to-F1-example}
Let $Y$ be the first Hirzebruch surface with Cox ring $S=\CC[y_0,\ldots, y_3]$. 
We choose the variables as in Example \ref{ex:not-fg}, in particular so that $y_1$ corresponds to the curve 
$B\subset Y$ of self-intersection $(-1)$.

\vspace{-0.7cm}
\hfill
\begin{tabular}{c}
\psset{unit=0.5cm} 
\begin{pspicture}(-5,-7)(4,4)
\psset{fillstyle=solid,fillcolor=vlgray,opacity=0.90}
\rput(!2 -3 2 \TT){\rnode{a}{}}
\rput(!2 -3 -2 \TT){\rnode{b}{}}
\rput(!2 2 -2 \TT){\rnode{c}{}}
\rput(!2 2 1 \TT){\rnode{d}{}}
\rput(!2 0 2 \TT){\rnode{e}{}}
\rput(!-2 -3 2 \TT){\rnode{f}{}}
\rput(!-2 -3 -2 \TT){\rnode{g}{}}
\rput(!-2 2 -2 \TT){\rnode{h}{}}
\rput(!-2 2 1 \TT){\rnode{i}{}}
\rput(!-2 0 2 \TT){\rnode{j}{}}
\pspolygon(f)(g)(h)(i)(j)
\psline(b)(g)
\pspolygon(a)(b)(c)(d)(e)
\pspolygon(a)(e)(j)(f)
\pspolygon(d)(c)(h)(i)
\rput(!2 -3 -6 \TT){\rnode{p}{}}
\rput(!2 2 -6 \TT){\rnode{q}{}}
\rput(!-2 2 -6 \TT){\rnode{r}{}}
\rput(!-2 -3 -6 \TT){\rnode{s}{}}
\pspolygon(p)(q)(r)(s)
\psline[arrows=->](!0 -0.5 -3.5 \TT)(!0 -0.5 -5.5 \TT)
\end{pspicture} \\
\Fig\label{fig:X-to-Y}
\end{tabular}

\vspace{-5.5cm}
\parshape 1 0cm 11.5cm
Set $Y_1=Y\times \PP^1$, and let $Z$ be the
curve $B\times \{[0:1]\}\subseteq Y\times \PP^1$.  Finally, let $X$ be the blowup of $Y_1$ along $Z$ with 
exceptional divisor $E$, and $\phi\colon X\longrightarrow Y$ the composition of the blowdown map and the projection 
$Y\times \PP^1\longrightarrow Y$.  

\parshape 1 0cm 11.5cm
Figure \ref{fig:X-to-Y} at right shows $P$ ($=X_{\geq 0}$), the polytope\footnote{ Here we are concerned only with 
$P$ as an abstract polyhedral complex, and not the precise realization of $P$ as a lattice polytope, 
and similarly for $Y_{\geq 0}$.}
of $X$, mapping to $Y_{\geq 0}$ the polytope of $Y$. 
(For toric maps passing to the ``non-negative part'' is functorial, see \S\ref{sec:functoriality-initial-discussion}.)

\parshape 1 0cm 11.5cm
By composing the blowdown map with the other projection we obtain a map
$X\longrightarrow \PP^1$. Let $H$ be the class of $\Osh_{\PP^1}(1)$ pulled back via this map, and set
$L=\Osh_{X}(-2H-2E)$.   Let us consider the computation of $R^{1}\phi_{*}L$. 

In this case it turns out that $C(L,1)=\{-1,\,0,\,1\}$ (via an isomorphism $\chi(T_{K})\cong \ZZ$).
This part of the calculation is explained later in the paper, see \S\ref{ex:X-to-F1-characters}.

We will use the following labelling system for the vertices and faces of $P$~:

\newcommand{\lblshape}[1]{\psframebox[fillcolor=white,opacity=1,framearc=0.2]{\small #1}}
\newcommand{\facelbl}[1]{\LARGE {#1}}

\begin{centering}
\begin{tabular}{ccc}
\begin{tabular}{c}
\psset{unit=0.5cm} 
\begin{pspicture}(-5,-4)(4,4)
\psset{fillstyle=solid,fillcolor=vlgray,opacity=0.80}
\rput(!2 -3 2 \TT){\rnode{a}{}}
\rput(!2 -3 -2 \TT){\rnode{b}{}}
\rput(!2 2 -2 \TT){\rnode{c}{}}
\rput(!2 2 1 \TT){\rnode{d}{}}
\rput(!2 0 2 \TT){\rnode{e}{}}
\rput(!-2 -3 2 \TT){\rnode{f}{}}
\rput(!-2 -3 -2 \TT){\rnode{g}{}}
\rput(!-2 2 -2 \TT){\rnode{h}{}}
\rput(!-2 2 1 \TT){\rnode{i}{}}
\rput(!-2 0 2 \TT){\rnode{j}{}}
\pspolygon(f)(g)(h)(i)(j) 
\psline(b)(g)
\rput(g){\lblshape{$g$}} 
\pspolygon(a)(b)(c)(d)(e)
\pspolygon(a)(e)(j)(f)
\pspolygon(d)(c)(h)(i)
\rput(a){\lblshape{$a$}} 
\rput(b){\lblshape{$b$}} 
\rput(c){\lblshape{$c$}} 
\rput(d){\lblshape{$d$}} 
\rput(e){\lblshape{$e$}} 
\rput(f){\lblshape{$f$}} 
\rput(h){\lblshape{$h$}} 
\rput(i){\lblshape{$i$}} 
\rput(j){\lblshape{$j$}} 
\end{pspicture} 
\end{tabular}
& \rule{3cm}{0cm} &
\begin{tabular}{c}
\psset{unit=0.5cm} 
\begin{pspicture}(-5,-4)(4,4)
\psset{fillstyle=solid,fillcolor=vlgray,opacity=0.80}
\rput(!2 -3 2 \TT){\rnode{a}{}}
\rput(!2 -3 -2 \TT){\rnode{b}{}}
\rput(!2 2 -2 \TT){\rnode{c}{}}
\rput(!2 2 1 \TT){\rnode{d}{}}
\rput(!2 0 2 \TT){\rnode{e}{}}
\rput(!-2 -3 2 \TT){\rnode{f}{}}
\rput(!-2 -3 -2 \TT){\rnode{g}{}}
\rput(!-2 2 -2 \TT){\rnode{h}{}}
\rput(!-2 2 1 \TT){\rnode{i}{}}
\rput(!-2 0 2 \TT){\rnode{j}{}}
\pspolygon(f)(g)(h)(i)(j) 
\psline(b)(g)
\rput(!-2 -0.5 0 \TT){\facelbl{2}} 
\rput(!0 -3 0 \TT){\facelbl{3}} 
\rput(!0 -0.5 -2 \TT){\facelbl{\pstilt{45}{5}}} 
\pspolygon(a)(b)(c)(d)(e)
\pspolygon(a)(e)(j)(f)
\pspolygon(d)(c)(h)(i)
\pspolygon(e)(d)(i)(j)
\rput(!2 -0.5 0 \TT){\facelbl{0}} 
\rput(!0 2 -0.5 \TT){\facelbl{1}} 
\rput(!0 -1.5 2 \TT){\facelbl{\pstilt{45}{4}}} 
\rput{45}(!0 1 1.5 \TT){\facelbl{\pstilt{110}{6}}} 
\end{pspicture} 
\end{tabular}\\
\small labels of the vertices & & \small labels of the faces \\
\multicolumn{3}{c}{\Fig\label{fig:P-labels}} \\
\end{tabular}\\
\end{centering}

For each $\mu\in C(L,1)$ the algorithm assigns a monomial ideal of $S$ to each vertex, edge, 
face, and $3$-cell of $P$.  Since we are interested in computing only $H^1$ of the complex, we will ignore the value 
on the unique $3$-cell (although it is always $S$, the entire Cox ring).
We will also only give the values on the edges where the ideal is nonzero. 

\psframebox{$\mu=-1$} :

\setlength{\extrarowheight}{0.1cm} 
\begin{tabular}{ccccc}
\begin{tabular}{ccccc}
\multicolumn{5}{c}{\und{Vertices}} \\
\lblshape{$a$} & $0$ & & \lblshape{$f$} & $0$ \\
\lblshape{$b$} & $0$ & & \lblshape{$g$} & $0$ \\
\lblshape{$c$} & $0$ & & \lblshape{$h$} & $0$ \\
\lblshape{$d$} & $(y_1^3)$ & & \lblshape{$i$} & $(y_1^3)$ \\
\lblshape{$e$} & $0$ & & \lblshape{$j$} & $0$ \\
\end{tabular}
& \rule{1.0cm}{0cm} & 
\begin{tabular}{ccccc}
\multicolumn{5}{c}{\und{Nonzero ed}g\und{es}} \\
\lblshape{$ab$} & $S$ & & \lblshape{$ij$} & $(y_1^3)$  \\
\lblshape{$cd$} & $S$ & & \lblshape{$di$} & $(y_1^3)$  \\
\lblshape{$hi$} & $S$ & & \lblshape{$de$} & $(y_1^3)$  \\
\lblshape{$fg$} & $S$  \\
\\
\end{tabular}
& \rule{1.0cm}{0cm} & 
\begin{tabular}{ccccc}
\multicolumn{5}{c}{\und{Faces}} \\
\lblshape{$0$} & $S$ & & \lblshape{$4$} & $0$ \\
\lblshape{$1$} & $S$ & & \lblshape{$5$} & $0$ \\
\lblshape{$2$} & $S$ & & \lblshape{$6$} & $(y_1^3)$ \\
\lblshape{$3$} & $S$  \\
\\
\end{tabular}
\end{tabular}

The resulting complex is, in each degree $i$, the sum over the terms associated to the $i$-dimensional cells
(e.g., in degree $0$ the sum over the vertex terms, in degree $1$ the sum over the edge terms, and so on).

The differentials of the complex are a combination of the usual signs for the cohomological complex associated to $P$
(i.e., the dual of the homological complex, where the signs are determined by choosing an orientation on each of the
$k$-cells), and the natural inclusion of ideals.  That is, if the map between two summands is nonzero, the ideal
of the source is always contained in the ideal of the target. 

In this case one computes that the $H^1$ of this complex is $S$. 

\psframebox{$\mu=0$} :

\setlength{\extrarowheight}{0.1cm} 
\begin{tabular}{ccccc}
\begin{tabular}{ccccc}
\multicolumn{5}{c}{\und{Vertices}} \\
\lblshape{$a$} & $S$ & & \lblshape{$f$} & $S$ \\
\lblshape{$b$} & $0$ & & \lblshape{$g$} & $0$ \\
\lblshape{$c$} & $0$ & & \lblshape{$h$} & $0$ \\
\lblshape{$d$} & $(y_1^2)$ & & \lblshape{$i$} & $(y_1^2)$ \\
\lblshape{$e$} & $(y_1^2)$ & & \lblshape{$j$} & $(y_1^2)$ \\
\end{tabular}
& \rule{0.75cm}{0cm} & 
\begin{tabular}{cccccccc}
\multicolumn{8}{c}{\und{Nonzero ed}g\und{es}} \\
\lblshape{$ab$} & $S$ & & \lblshape{$ij$} & $(y_1^2)$ & & \lblshape{$af$} & $S$ \\
\lblshape{$cd$} & $S$ & & \lblshape{$di$} & $(y_1^2)$ & & \lblshape{$ae$} & $S$ \\
\lblshape{$hi$} & $S$ & & \lblshape{$de$} & $(y_1^2)$ & & \lblshape{$fj$} & $S$ \\
\lblshape{$fg$} & $S$ & & \lblshape{$ej$} & $(y_1^2)$ &  \\
\\
\end{tabular}
& \rule{0.75cm}{0cm} & 
\begin{tabular}{ccccc}
\multicolumn{5}{c}{\und{Faces}} \\
\lblshape{$0$} & $S$ & & \lblshape{$4$} & $S$ \\
\lblshape{$1$} & $S$ & & \lblshape{$5$} & $0$ \\
\lblshape{$2$} & $S$ & & \lblshape{$6$} & $(y_1^2)$ \\
\lblshape{$3$} & $S$  \\
\\
\end{tabular}
\end{tabular}

In this case $H^1$ of the complex is $S/(y_1^2)$.

\psframebox{$\mu=1$} : In this case the complex is like that for $\mu=0$, but with every instance of $(y_1^2)$
replaced by $(y_1)$.  The first cohomology of the resulting complex is $S/(y_1)$. 

In this particular example all the shifts $E_{\mu}$ (step ({\em ii}) from \S\ref{sec:algorithm-description}) are zero.
Thus, by the theorem, $R^1\phi_{*}L$ is the sheafification of the $S$-module 
$S\oplus S/(y_1^2) \oplus S/(y_1)$, and therefore
$R^1\phi_{*}L = \Osh_{Y} \oplus \frac{\Osh_{Y}}{\Osh_{Y}(-2B)} \oplus \frac{\Osh_{Y}}{\Osh_{Y}(-B)}$.
It is not hard to calculate independently that this is the correct answer (for instance by writing
down short exact sequences including $L$, and using the resulting long exact sequences of higher direct images).

\bpoint{Remarks} 
(\Euler{1}) For a given $\mu\in C(L,i)$, the summand $(R^{i}\phi_{*}L)_{\mu}$ may be reducible 
as an $\Osh_{Y}$-module.  That is, the decomposition of \eqref{eqn:TK-splitting} is in general
not a maximal decomposition into summands, although this was the case in Example \ref{sec:X-to-F1-example} above.

As an instance of this, let $Y$ be a toric surface, $\phi\colon X\longrightarrow Y$ the blowup of $Y$ at two toric 
fixed points $p_1$ and $p_2$, with respective exceptional divisors $E_1$ and $E_2$.
Setting $L=\Osh_{X}(-2E_1-2E_2)$, then (as in Example \ref{ex:not-fg}) $R^{1}\phi_{*}L = \Osh_{p_1}\oplus \Osh_{p_2}$.

However, since $\phi$ is a birational map, the induced map $T_{X}\longrightarrow T_{Y}$ is an isomorphism, and so 
$T_{K}$ is reduced to a single point (i.e., the trivial group), with trivial character group.  Thus, in this
case $C(L,1)=\{0\}$, and the two summands of $R^1\phi_{*}L$ belong to the same term $(R^1\phi_{*}L)_{0}$. 

(\Euler{2}) As Example \ref{sec:X-to-F1-example} suggests, for a given $\mu$, once one knows the ideals of $S$
associated to each maximal face of $P$, the ideal for any smaller face is obtained
by intersecting the ideals associated to all the maximal faces which contain it. 

For instance, in the complex for $\mu=0$ (in Example \ref{sec:X-to-F1-example}), 
the edge $ij$ is contained in faces $2$ and $6$, and
the intersection of the ideals $S$ and $(y_1^2)$ associated to those faces is the ideal $(y_1^2)$ associated
to $ij$.  Similarly, the vertex $a$ is contained in faces $0$, $3$, and $4$, each of which are assigned the ideal 
$S$ and which therefore has intersection $S$, the ideal assigned to $a$.   Finally, since face $5$ is assigned
the ideal $0$, any vertex or edge contained in that face will also be assigned the ideal $0$.

(\Euler{3})  
For each $\mu\in C(L,i)$ 
the complexes of monomial ideals is constructed by first constructing a complex 
$\CComp_{P}(\Ish_{D_{\mu}})$ of ideal sheaves on $Y$
(see the description of \S\ref{sec:construction-of-D-E-complex} below), 
along with an isomorphism\footnote{ Here, and for the rest of the paper, $\Hsh^{i}$ denotes the cohomology sheaf
of a complex of sheaves, and not the notation of \S\ref{sec:Temp-def-of-Hi}.}
$\Hsh^{i}(\CComp_{P}(\Ish_{D_{\mu}}))\otimes \Osh_{Y}(E_{\mu})\longrightarrow (R^{i}\phi_{*}L)_{\mu}$ 
Each term of the complex is a $T_Y$-equivariant ideal sheaf and hence corresponds to a monomial ideal of $S$.

In some ways this is our answer to the issue raised in \S\ref{sec:issue-finding-module}, of how to find an $S$-module
producing a given sheaf.  For an ideal sheaf, or an ideal sheaf twisted by a line bundle, it is easy to find
a corresponding module. The construction of \S\ref{sec:construction-of-D-E-complex} realizes the 
sheaf $(R^{i}\phi_{*}L)_{\mu}$ as the $i$-th cohomology 
group of a complex of ideal sheaves, possibly twisted by a line bundle.  Since the procedure of going from
modules to sheaves is an exact functor, to produce an $S$-module whose corresponding sheaf is $(R^{i}\phi_{*}L)_{\mu}$ 
it suffices to lift the terms of the complex to $S$-modules and then compute the $i$-th cohomology.

\bpoint{Outline of the paper} 

{\bf \S\ref{sec:cell-complex-toric-variety} :} This section contains a brief reminder of the general facts about cell 
complexes, 
as well as a discussion of the non-negative part, $X_{\geq 0}$ of a toric variety $X$.  It also states a result
we will need in constructing the a reduced \v{C}ech complex for a semi-proper toric variety, namely that
if $X$ is proper $X_{\geq 0}$ is homeomorphic to a ball (Theorem \ref{thm:top-structure-of-X}), 
and therefore has the cohomology of a point.

{\bf \S\ref{sec:reduced-cech-complexes} :} This section contains the axiomatic development of the idea of a 
reduced \v{C}ech 
complex.  Among the results are Theorem  \ref{thm:reduced-cech-computes-cohomology}, showing that when
the axioms (P1)--(P4) (see \S\ref{sec:Def-of-reduced-complex} and \S\ref{sec:compatibility-conditions}) are satisfied,
a reduced \v{C}ech complex based on $P$ computes the cohomology of any
sheaf of abelian groups $\Fsh$ whose higher cohomology vanishes on the elements of the cover, and
Theorem \ref{thm:reduced-cech-spectral-sequence}, giving a general `reduced \v{C}ech-to-cohomology' spectral
sequence.  This section also establishes that conditions (P1)--(P4) are necessary and sufficient for a
reduced \v{C}ech complex to compute cohomology. 

Example \ref{ex:axioms-dont-force-values-on-faces} gives an example of a cover satisfying the axioms showing that 
for a $\sigma\in P$, the open set $U_{\sigma}$ associated to $\sigma$ is not necessarily the intersection 
of the $U_{v}$  over those vertices $v\in \sigma$ (as is the case for the usual \v{C}ech cover, or the 
reduced \v{C}ech complex associated to a toric variety).  
Sections \ref{sec:comparison-argument-suggestion}--\ref{Ex:comparison-argument-example}
discuss the issue of whether one can deduce these results from the results about the usual \v{C}ech complex.

{\bf \S\ref{sec:covers-of-semi-proper} :} 
Theorem \ref{thm:construction-of-toric-reduced-cover} constructs a reduced \v{C}ech complex for any 
semi-proper toric variety.
This section also contains results on semi-proper maps, a discussion of toric morphisms, 
a proof that if $\phi\colon X\longrightarrow Y$ is a
toric morphism, and $X\longrightarrow Z\longrightarrow Y$ its Stein factorization, then $Z$ is itself
a toric variety (Lemma \ref{lem:toric-stein-factorization}), as well as intermediate results necessary for the proof 
of the theorem cited above.  In particular, Proposition \ref{prop:cohomology-of-W-pos} 
gives a useful result on the topology of the non-negative points of the fibre of a semi-proper toric map over a 
torus fixed point.

{\bf \S\ref{sec:construction-of-D-E-complex} :} 
Given a toric fibration $\phi\colon X\longrightarrow Y$ between smooth toric varieties, a $T_X$-linearized line 
bundle $L$ on $X$, and $\mu\in C(L,i)$, this section produces : two Cartier divisors $D_{\mu}$ on $X$
and $E_{\mu}$ on $Y$; a complex of ideal sheaves $\CComp_{P}(\Ish_{D_{\mu}})$ on $Y$; and a morphism 
$\Hsh^{i}\bigl(\CComp_{P}(\Ish_{D_{\mu}})\bigr)\otimes \Osh_{Y}(E_{\mu})\longrightarrow (R^{i}\phi_{*}L){\mu}$, which
is guaranteed to be an isomorphism. 

The description of the general construction of $\CComp_{P}(\Ish_{D_{\mu}})$, along with a morphism to 
$(R^{i}\phi_{*}L)_{\mu}$ is given in \S\ref{sec:basic-construction}--\S\ref{sec:variation-on-the-construction},
while the result constructing a pair $(D_{\mu},E_{\mu})$ such that the map to $(R^{i}\phi_{*}L)_{\mu}$ is
an isomorphism is Theorem \ref{thm:D-E-theorem}.

{\bf \S\ref{sec:appendix-algorithm} :} The above constructions depend on knowing the finite set $C(L,i)$ and 
on knowing local generators for $(R^{i}\phi_{*}L)$ as an $\Osh_{Y}$-module (e.g., over 
torus stable affines $V\subseteq Y$).  In this appendix, due to Sasha Zotine, the previous 
constructions are made effective, and turned into an algorithm, see \S\ref{sec:appendix-alg1} 
and \S\ref{sec:appendix-alg2}.  The resulting method has been implemented in {\it Macaulay2}, see \cite{Zo}. 

{\bf \S\ref{sec:appendix-examples} :} This section, also due to Sasha Zotine, contains further examples 
illustrating the algorithm, and exploring related questions.  
For instance, Example \ref{ex:convexity-of-weights} examines the question as to whether the set $C(L,i)$, 
considered as a subset of $\chi(T_{K})_{\RR}$, is the lattice points of a convex polytope.

\bpoint{Acknowledgements}  The idea to look for an algorithm to compute higher direct images of line bundles
(or in general, coherent sheaves) under toric morphisms is due to Gregory G.\ Smith.  The authors also thank 
Gregory Smith for describing the construction in the literature of a reduced \v{C}ech complex
on projective simplicial toric varieties (i.e., the construction described in \S\ref{sec:Oda-construction}) , 
and for encouragement while the arguments of this paper were being worked out.

\section{The cell complex associated to a proper toric variety} 
\label{sec:cell-complex-toric-variety}

\point 
As in the standard \v{C}ech complex, a reduced \v{C}ech complex on a topological space $X$ 
is a complex formed by the sections of a sheaf on a collection of open sets in $X$.  
In constructing the reduced \v{C}ech complex we need a way to keep track of the combinatorics of the restriction maps 
and the associated signs.  In this article we encode this data using a CW-complex, 
also called a cell complex in \cite[page 5]{H}.
In general the cell complex does not need to have any direct connection to the underlying topological space $X$, 
although this will be the case in our main application to toric varieties.

In this section we recall the definition of a cell complex, and, when $X$ is a proper toric variety, 
that the nonnegative part of $X$, $X_{\geq 0}$, is naturally a regular cell complex, homeomorphic to a sphere.
When $X$ is projective this cell complex is the standard polytope associated to $X$ under any projective 
embedding, and for this reason we use $P$ as the name of the cell complex throughout the article.

\bpoint{Cell complexes} 
\label{def:cell-complex}
We use the definition of CW-complex from \cite[IV. \S3]{M}.  
Specifically, a CW-complex is a Hausdorff
topological space $P$, along with an ascending sequence of closed subspaces
$$P^0 \subset P^1 \subset P^2 \subset \cdots $$
satisfying the conditions~:

\RomanList
\begin{enumerate}
\item $P^0$ has the discrete topology.
\item For each $n>0$, $P^n$ is obtained from $P^{n-1}$ by attaching a collection of $n$-cells.
\PauseEnumerate

We note that the topology of $P^n$ must be that constructed from the topology of $P^{n-1}$ and the attaching maps.  
See \cite[p.\ 519]{H} or \cite[IV.\ \S2]{M}.

In line with the notation for polytopes or cones, we use $\sigma$ (or similar letters) for the name of a closed cell. 
We emphasize that $\sigma$ denotes a closed cell, and not the interior.
We also use the familiar notation $\sigma_1\subseteq \sigma_2$ to denote inclusion of cells, analogous to
inclusions of faces in polytopes or cones, and use $|\sigma|$ for the dimension of a cell $\sigma$. 

This data is further required to satisfy~:

\ResumeEnumerate
\item $P$ is the union of the $P^n$, for $n\geq 0$.
\item The space $P$ and the subspaces $P^n$ all have the weak topology~: A subset $A$ is closed if an only if
$A\cap \sigma$ is closed for each $n$-cell $\sigma$  and for all $n\geq 0$. 
\end{enumerate}
\AlphaList

If $P$ is a finite cell complex (i.e., $P$ has only finitely many cells) then condition ({\em iv}) is automatic, see
\cite[p. 80]{M} or \cite[p.\ 519]{H}. 
This will be the case in our application to toric varieties.

\bpoint{Incidence functions} 
\label{sec:incidence-functions}
We recall that, given a cell complex $P$, for each pair of cells $\sigma$, $\sigma'$ of dimensions $i$ and $i+1$ 
respectively one may define an incidence number $\ep(\sigma,\sigma')\in \ZZ$, such that the collection of such
numbers satisfies 

\begin{enumerate}
\item Fixing an $i+1$-cell $\sigma'$, $\ep(\sigma,\sigma')=0$ for all but finitely many $i$-cells $\sigma$.
\item If $\sigma\not\subset \sigma'$, then $\ep(\sigma,\sigma')=0$. 
\item Fixing an $i$-cell $\sigma$, and $(i+2)$-cell $\sigma''$, 
$$\sum_{\sigma\subset \sigma' \subset \sigma'', |\sigma'|=i+1} \ep(\sigma,\sigma')\ep(\sigma',\sigma'')=0.$$
\end{enumerate}

See \cite[IV. \S5, Lemmas 5.1 and 5.2]{M}.

We call such a function $\ep$ an incidence function associated to $P$.  Such a function $\ep$ is not unique; its
construction depends on a choice of orientation for each of the cells. However, this non-uniqueness, which manifests
itself in the signs of the $\ep(\sigma,\sigma')$, will not affect the outcome of any of the constructions using $\ep$.

\bpoint{The non-negative part of a toric variety}  
\label{sec:non-neg-part-def}
Let $X$ be a complex toric variety, with fan $\Delta_{X} \subseteq N_{\RR}$, where $M$ and $N$ denote as usual the
lattice of torus characters and one-parameter subgroups respectively.  

As was first realized by Ehlers \cite[\S IV]{E}, since the transition data between the open affine
charts of $X$ are monomial maps, there is sense in asking for the points of $X$ with coordinates in 
$\RR_{\geq 0} := \left\{x\in \RR \st x\geq 0\right\}$. 

Intrinsically, for each cone $\sigma\in \Delta_{X}$ with associated semigroup $S_{\sigma}=\sigma^{\vee}\cap M$,
the $\RR_{\geq 0}$-points of the affine variety $U_{\sigma}\subseteq X$ may be identified with the semigroup
homomorphisms $\Hom_{\sg}(S_{\sigma},\RR_{\geq 0})$, where $\RR_{\geq 0}$ is thought of as a semigroup under 
multiplication. 
One may similarly look at the points of $X$ with coordinates in  $\RR_{>0} := \left\{x\in \RR \st x>0\right\}$.

We use $X_{\geq 0}$ and $X_{>0}$ respectively for the corresponding topological subspaces of $X(\CC)$.  

As elementary examples, we note that 
$$
\PP^{n}_{\geq 0} = \left\{ [x_0:x_1:\cdots : x_n] \st x_i\in \RR_{\geq 0} \right\}
= \left\{ [x_0:x_1:\cdots : x_n] \st x_i\in \RR_{\geq 0},\,\sum x_i = 1 \right\}
$$
is the $n$-dimensional simplex, and that when $X=(\PP^1)^{n}$,  $X_{\geq 0} = 
(\PP^1_{\geq 0})^n \cong [0,1]^n$ is the $n$-dimensional cube. 

\point 
\label{sec:cell-complex-structure}
As shown in \cite[Chap.\ I]{AMRT} via the exponential map one can identify
$N_{\RR}$ and $X_{>0}$.  Furthermore, via this map $N_{\RR}$ acts on $X_{\geq 0}$, breaking it into
finitely many orbits, each one corresponding to a cone of $\Delta_{X}$ in a way parallel to the usual decomposition
of $X$ under the action of $T_{X}$, see \cite[loc. cit.]{AMRT}. 

For each cone $\sigma\in \Delta_{X}$, the stabilizer of a point in the orbit corresponding to $\sigma$
is the subspace $\langle \sigma\rangle\subseteq N_{\RR}$, and so the orbit is isomorphic to 
$N_{\RR}/\langle \sigma\rangle$.  Thus $X_{\geq 0}$ decomposes naturally into open cells with the largest cell,
$X_{>0}$, identified with $N_{\RR}$ itself. 

We will need the result below in our construction of reduced \v{C}ech complexes for a semi-proper toric varieties.

\bpoint{Theorem} 
\label{thm:top-structure-of-X}
Let $X$ be a complete toric variety over $\CC$, and set $n=\dim(X)$.  Then 

\begin{enumerate}
\item The cells described above give $X_{\geq 0}$ the structure of a regular cell complex.

\smallskip
\item $X_{\geq 0}$ is homeomorphic to the $n$-dimensional ball.
\end{enumerate}

The key step is to show ({\em b}), which then easily implies ({\em a}).   

When $X$ is a smooth and proper toric variety the
fact that $X_{\geq 0}$ is homeomorphic to a ball was first shown by Ehlers, \cite[p. 143, Lemma 1]{E}.
The case of a projective, but possibly singular, toric variety was established by Jurkiewicz \cite{J}, see
also \cite[p. 81, Proposition]{F}.  In this case the topological space $X_{\geq 0}$ may be identified
with the image of the moment map, \cite[p.\ 82--83]{F}.
The result for the common generalization, that of a possibly singular proper toric variety, appears
as \cite[Theorem 2.15]{R}.

\point 
\label{sec:functoriality-initial-discussion}
The association $X\rightsquigarrow X_{\geq 0}$ is a functor from the category of toric varieties (with toric morphisms) 
to topological spaces.  If $\phi\colon X\longrightarrow Y$ is a toric morphism (\S\ref{sec:toric-morphisms})
then locally $\phi$ is given by monomial maps, and thus preserves the coordinates which are in $\RR_{\geq 0}$.
One may then ask if this functor is also a functor to the category of cell complexes.  That is, if the induced
map $\phi_{\geq 0}\colon X_{\geq 0}\longrightarrow Y_{\geq 0}$ is always a map of cell complexes, with the 
cell structure as given in \S\ref{sec:cell-complex-structure}.  

In general, the answer is no, although $\phi_{\geq 0}$ does respect this choice of cellular structure whenever 
$\phi$ is surjective, see \S\ref{sec:funct-of-non-neg-construction}.

\section{Reduced \v{C}ech complexes} 
\label{sec:reduced-cech-complexes}

\point 
\label{sec:Def-of-reduced-complex}
Let $X$ be a topological space, $P$ a cell complex as in \S\ref{def:cell-complex}, and let
$\ep$ be an incidence function for $P$ (\S\ref{sec:incidence-functions}).
We also assume given an assignment
$\sigma\mapsto \psi(\sigma)=:U_{\sigma}$ sending each cell $\sigma\in P$ to an open set $U_{\sigma}\subseteq X$.
To reduce notation we often omit reference to $\psi$ (e.g., naming groups $\HP$ instead of $\check{H}_{(P,\psi)}$).
We assume the pair $(P,\psi)$ satisfies the following conditions~:

\begin{itemize}
\item[(P1)] The assignment is order reversing with respect to inclusion~: if $\tau\subseteq \sigma$ then $U_{\sigma}\subseteq U_{\tau}$. 

\smallskip
\item[(P2)] The $U_{\sigma}$ cover $X$~: $\bigcup_{\sigma\in P}U_{\sigma} = X$.
\end{itemize}

Given (P1), (P2) is equivalent to the condition that $\bigcup_{v\in P^{0}} U_{v}=X$, where $v$ runs through the 
vertices (i.e., $0$-cells) of $P$.  Let $\Fsh$ be a sheaf of abelian groups on $X$.  For each $i\in \NN$ we set 
$$\CComp_{P}^{i}(\Fsh) := \prod_{\sigma\in P,\,|\sigma|=i} \Gamma(U_{\sigma},\Fsh),$$
where $|\sigma|$ denotes the dimension of $\sigma$. For each $i\in \NN$ we define a differential
$d^{i}\colon \CComp_{P}^{i}(\Fsh)\longrightarrow \CComp_{P}^{i+1}(\Fsh)$ by the following rule on the components
of the products.  For $\sigma$, $\sigma'\in P$, $|\sigma|=i$, $|\sigma'|=i+1$, we put 

\begin{equation}\label{eqn:def-of-differential}
d^{i}_{\sigma,\sigma'} : = \left\{
\begin{array}{cl}
0 & \mbox{if $\sigma\not\subset \sigma'$} \\
\ep(\sigma,\sigma')\rho_{U_{\sigma},U_{\sigma'}} & \mbox{if $\sigma\subset \sigma'$}. \rule{0cm}{0.4cm}  \\
\end{array}
\right.
\end{equation}

Here $\rho_{U_{\sigma},U_{\sigma'}}$ is the restriction map 
$\Gamma(U_{\sigma},\Fsh)\longrightarrow \Gamma(U_{\sigma'},\Fsh)$ induced by the inclusion
$U_{\sigma'}\subseteq U_{\sigma}$.  For fixed $\sigma'$ there are only finitely many $\sigma$ such that 
$d^{i}_{\sigma,\sigma'}\neq 0$, since there are only finitely many
$\sigma$ such that $\ep(\sigma,\sigma')\neq 0$. 
(In contrast, for fixed $\sigma$ there may be infinitely many $\sigma'$ with $\ep(\sigma,\sigma')\neq 0$.)
Thus, for each $\sigma'$ the map 
$$\sum_{\sigma} d^{i}_{\sigma,\sigma'}\colon \CComp_{P}^{i}(\Fsh) \longrightarrow \Gamma(U_{\sigma'},\Fsh)$$
makes sense, and induces a map 
$d^{i}\colon\CComp_{P}^{i}(\Fsh)\longrightarrow \CComp_{P}^{i+1}(\Fsh)$.
The functoriality of the restriction homomorphisms (i.e., under composite inclusions) and the axioms for
an incidence function show that $d^{i+1}\circ d^{i}=0$. 
We denote the resulting complex by $\CComplex_{P}(\Fsh)$, and call it the {\em reduced \v{C}ech complex} (based
on $(P,\psi)$ and $\Fsh$). 

We use the term reduced since in general this complex may be smaller than the usual \v{C}ech complex
associated to the cover.  See \S\ref{sec:first-example-of-reduced} for an example.

\bpoint{Functoriality} 
\label{sec:functoriality-of-complex} 
The construction in \S\ref{sec:Def-of-reduced-complex} is functorial in $\Fsh$.  If 
$\varphi\colon \Fsh\longrightarrow \Gsh$ is a map of sheaves of abelian groups, then $\varphi$ induces maps
$\Gamma(U_{\sigma},\Fsh)\longrightarrow \Gamma(U_{\sigma},\Gsh)$ for each $\sigma\in P$, and hence maps
$\CComp_{P}^{i}(\Fsh)\longrightarrow \CComp_{P}^{i}(\Gsh)$ for each $i$.
The morphisms on global sections over each $U_{\sigma}$ commute with restriction, and hence with the 
differential constructed from \eqref{eqn:def-of-differential}. This gives rise to a morphism of complexes 
$\CComplex_{P}(\varphi)\colon \CComplex_{P}(\Fsh)\longrightarrow \CComplex_{P}(\Gsh)$.  
Given another homomorphism $\varphi'\colon \Gsh\longrightarrow \Hsh$ of sheaves of abelian groups, we have
$\CComplex_{P}(\varphi'\circ\varphi) = \CComplex_{P}(\varphi')\circ \CComplex_{P}(\varphi)$, 
since the morphisms on global sections are similarly functorial.

\bpoint{Definition of $\HP^{i}(\Fsh)$} 
\label{sec:Def-of-HP}
For any sheaf $\Fsh$ of abelian groups on $X$, and for each $i\geq 0$ we put 
$\HP^{i}(\Fsh):=H^{i}(\CComplex_{P}(\Fsh))$. 
Given a morphism $\varphi\colon \Fsh\longrightarrow \Gsh$ of sheaves, the morphism of complexes from
\S\ref{sec:functoriality-of-complex} induces morphisms 
$\HP^{i}(\varphi)\colon \HP^{i}(\Fsh)\longrightarrow \HP^{i}(\Gsh)$ for each $i$, the construction again being
functorial in $\Fsh$.

\point {\bf Examples}. 
\label{sec:first-example-of-reduced}
(\Euler{1}) Let $X$ be a topological space, and $\{U_{\alpha}\}_{\alpha\in I}$ an open
cover of $X$.  To the cover we associate the nerve complex $P$, the simplicial complex with vertices the 
$\alpha\in I$, and $\sigma=\{\alpha_0,\alpha_1,\ldots, \alpha_r\}\in P$ if $\cap_{i=0}^{r} U_{\alpha_i}\neq\emptyset$.
If we give $I$ a total order, then there is a standard incidence function on $P$. If 
$\sigma'=\{\alpha_0,\alpha_1,\ldots, \alpha_{i+1}\}$, with $\alpha_0 < \alpha_1<\cdots < \alpha_{i+1}$, and 
$\sigma=\{\alpha_0,\ldots, \alpha_{j-1},\hat{\alpha}_{j},\alpha_{j+1},\cdots,\alpha_{i+1}\}$, 
then we set $\ep(\sigma,\sigma')=(-1)^{j}$.

The simplicial complex $P$, along with the assignment $\sigma\mapsto U_{\sigma}:=\cap_{\alpha\in \sigma} U_{\alpha}$ 
satisfies (P1) and (P2) above, and for a given sheaf of abelian groups $\Fsh$, the complex
$\CComplex_{P}(\Fsh)$ is the usual \v{C}ech complex associated to the cover. 

(\Euler{2}) To give a concrete example of a smaller complex, let $X=\PP^1\times \PP^1$.  
We make $X$ into a toric variety by letting $T=(\CC^{*})^2$ act on $X$ in the standard way, 
$$(t_1,t_2)\cdot ([x_1:y_1],[x_2,y_2])= ([t_1x_1:y_1], [t_2x_2:y_2]).$$
There are $9$ orbits of $T$ on $X$ under this action.

\newgray{vlgray}{0.9}
\hfill
\begin{tabular}{c}
\psset{unit=2.0cm}
\begin{pspicture}(-1,-1)(1,1)
\pspolygon[fillstyle=solid,fillcolor=vlgray,linecolor=lightgray](-1,-1)(-1,1)(1,1)(1,-1)
\psline(-0.75,-0.85)(-0.75,0.85)
\psline(0.75,-0.85)(0.75,0.85)
\psline(-0.85,-0.75)(0.85,-0.75)
\psline(-0.85,0.75)(0.85,0.75)
\psset{fillstyle=solid,fillcolor=white}
\pscircle(-0.75,-0.75){0.03}
\pscircle(-0.75,0.75){0.03}
\pscircle(0.75,-0.75){0.03}
\pscircle(0.75,0.75){0.03}
\rput(0,0){\tiny $\sigma$}
\rput(0,0.88){\tiny $\tau_{12}$} 
\rput(-0.88,0){\tiny $\tau_{13}$}
\rput(0,-0.88){\tiny $\tau_{34}$} 
\rput(0.88,0){\tiny $\tau_{24}$}
\rput(-0.88,0.88){\tiny $v_1$} 
\rput(0.88,0.88){\tiny $v_2$} 
\rput(-0.88,-0.88){\tiny $v_3$} 
\rput(0.88,-0.88){\tiny $v_4$} 
\rput(0,-1.45){$(\PP^1\times\PP^1)_{\geq 0}$}
\end{pspicture}\\
\Fig\label{fig:P1xP1}
\end{tabular}

\vspace{-4cm}
\parshape 1 0cm 12cm
Four  $0$-dimensional orbits~:

\smallskip
\begin{tabular}{ll}
$O_{v_1} = \{([0:1],[0:1])\}$, & $O_{v_2} = \{([1:0],[0:1])\}$ \\
$O_{v_3} = \{([0:1],[1:0])\}$, & and $O_{v_4} = \{([1:0],[1:0])\}$; \rule{0cm}{0.5cm}\\
\end{tabular}

four  $1$-dimensional orbits~:

\smallskip
\begin{tabular}{ll}
$O_{\tau_{12}} = \{([*:1],[0:1])\}$, & $O_{\tau_{34}} = \{([*:1],[1:0])\}$ \\
$O_{\tau_{13}} = \{([0:1],[*:1])\}$, & and $O_{\tau_{24}} = \{([1:0],[*:1])\}$; \rule{0cm}{0.5cm}\\
\end{tabular}

and one $2$-dimensional orbit : 

\smallskip
$O_{\sigma} = \{([*:1],[*:1])\}$. 

\parshape 1 0cm \textwidth 
Here and above $*$ denotes any non-zero element of $\CC$ (not necessarily the same element in each coordinate in the 
case of $O_{\sigma}$). 

The cell complex $P$ we associate to this situation is the square $X_{\geq 0}$ (\S\ref{sec:non-neg-part-def}),
with vertices $v_1$,\ldots, $v_4$, 
edges $\tau_{12}$, \ldots, $\tau_{34}$, and $2$-dimensional face $\sigma$.  A vertex $v_{k}$ is contained in an 
edge $\tau_{ij}$ if $O_{v_{k}}\subseteq \overline{O}_{\tau_{ij}}$ (here denoted in the labelling by the fact 
that $k\in \{i,\,j\}$), and where $\sigma$ contains all of the edges and vertices.   We give $P$ any valid incidence
function (e.g., one given by choosing orientations on the edges, and on $\sigma$).

To each face $\gamma\in P$ we associate the set 

\begin{equation}\label{eqn:square-example}
U_{\gamma} := \bigcup_{O_{\gamma}\subseteq \overline{O}_{\gamma'}} O_{\gamma'}.
\end{equation}

In this case one can check that each $U_{v_i}$ is an $\AA^2$, each $U_{\tau_{ij}}$ is a $\CC^{*}\times \CC$,
and $U_{\sigma}\cong (\CC^{*})^2$.  In particular, each $U_{\gamma}$ is an open affine subset of $X$.  (On a
toric variety this is a general feature when taking the union of torus orbits whose closure contains a fixed orbit, 
see \cite[p.\ 24, proof of theorem 6]{TEI}.)

Condition (P1) is satisfied by virtue of \eqref{eqn:square-example}, and (P2) is satisfied since
each orbit appears at least once, and so the $U_{\gamma}$ cover $X$.

The $f$-vector of $P$ is $(4,4,1)$, and thus the resulting reduced \v{C}ech complex has $4$ summands in degree $0$, 
$4$ in degree $1$, and $1$ in degree $2$.  On the other hand, starting with the $4$ open sets $U_{v_1}$,\ldots,
$U_{v_4}$, if we construct the total \v{C}ech complex, the nerve complex is the $3$-dimensional simplex with
$f$-vector $(4,6,4,1)$ (including a term in degree $3$, larger than the dimension of $X$). 
More generally, applying a similar construction to $(\PP^1)^n$, the polytope for the reduced \v{C}ech complex is the
$n$-dimensional cube, while the total \v{C}ech complex has nerve complex the $(2^{n}-1)$-dimensional simplex. 

\point
Since the cover in the previous example is by open affines, it is well-known that the total \v{C}ech complex 
computes the cohomology of quasi-coherent sheaves on $X$.   
We seek conditions so that the reduced \v{C}ech complex will compute
cohomology of an appropriate class of sheaves, both here, and in the general case of an arbitrary topological space.

In the case of this example, or the similar extension of this idea to semi-proper toric varieties (see 
\S\ref{sec:covers-of-semi-proper}) the class will include quasi-coherent sheaves.  
For proper toric varieties this gives a way of computing cohomology where
the combinatorics of the complex are determined by the ``polytope of $X$'' 
(i.e, $X_{\geq 0}$ from \S\ref{sec:non-neg-part-def}).

\point 
\label{sec:compatibility-conditions}
The conditions we need are compatibility conditions between 
the topology of the polytope $P$ and the combinatorics of the cover. 
For any $x\in X$ we put
$$P_{x} := \left\{\sigma\in P \st x\in U_{\sigma}\right\}.$$

By (P2), $P_{x}$ is always nonempty.  By (P1), if $\sigma\in P_{x}$ and $\tau\subseteq \sigma$ then $\tau\in P_{x}$.  
Thus, $P_{x}$ is a sub cellular complex of $P$.  

The required compatibility condition is that each $P_{x}$ should have the cohomology of a point.
For convenience in Lemmas \ref{lem:reduced-H0-correct} and \ref{lem:Godement-zero-in-reduced-cech}
below, we break this into two conditions.

\begin{itemize}
\item[(P3)] For each $x\in X$, $H^0(P_{x},\ZZ)=\ZZ$.

\smallskip
\item[(P4)] For each $x\in X$, $H^i(P_{x},\ZZ)=0$ for all $i\geq 1$.
\end{itemize}

In the case of the total \v{C}ech complex of a cover with associated nerve complex $P$, 
each $P_{x}$ is a simplex (possibly infinite), and so satisfies (P3) and (P4) above. 

\newpage
\bpoint{Lemma}\label{lem:reduced-H0-correct}
Condition (P3) is a necessary and sufficient condition so that $\HP^{0}(\Fsh)=H^0(X,\Fsh)$ for all sheaves of
abelian groups $\Fsh$ on $X$.

\bpf Let $P(0)$ denote the set of vertices of $P$.  The restriction maps 
$$
\begin{array}{rcll}
H^{0}(X,\Fsh) & \longrightarrow & \prod_{v\in P(0)} \Fsh(U_{v}) & (= \CComp_{P}^{0}(\Fsh)) \\
s & \xmapsto{\rule{0.5cm}{0cm}} & \prod s|_{U_v}  \rule{0cm}{0.6cm} \\
\end{array}
$$
induce a map $H^{0}(X,\Fsh)\longrightarrow \HP^{0}(\Fsh)$.  This map is injective since, by (P1) and (P2), the $U_{v}$
are a cover of $X$.  

To show that the map is surjective, suppose that $\prod s_{v}$ ($s_{v}\in \Fsh(U_{v})$)
is an element of $\CComp^{0}_{P}(\Fsh)$ which goes to zero under the differential map.   We will show that the $s_{v}$
patch together to a global section $s$.
To do that, it is sufficient to show that for each $x\in X$, the stalks $s_{v,x}$, for those $v$ such that
$x\in U_{v}$, are all equal in $\Fsh_{x}$.

So, fix $x\in X$.  The vertices $v$ such that $x\in U_{v}$ are the vertices of $P_{x}$.  Let $v$ and $w$ be
any two vertices of $P_{x}$.  We assume that (P3) holds. Then, since $H^0(P_{x},\ZZ)=\ZZ$, there is a finite list
of vertices $v=v_0$, $v_1$, $v_2$, \ldots, $v_m=w$,  and edges $e_1$, $e_2$, \ldots, $e_m$, all in $P_{x}$ 
with $e_i$ connecting $v_{i-1}$ to $v_{i}$ for $i=1$,\ldots, $m$. 

The cocycle condition associated to each $e_{i}$ is $s_{v_{i-1}}|_{U_{e_i}} = s_{v_{i}}|_{U_{e_i}}$. 
Since $x\in U_{e_i}$ ($e_i$ is an edge of $P_{x}$), this gives $s_{v_{i-1},x}=s_{v_{i},x}$ for $i=1$, \ldots, $m$, 
and so $s_{v,x}=s_{w,x}$.   Thus $H^0(X,\Fsh) = \HP^{0}(\Fsh)$. 

To see that (P3) is a necessary condition, let $x$ be a point of $X$ and $i\colon \{x\}\hookrightarrow X$ the 
canonical injection. We consider the sheaf $i_{*}\ZZ$ on $X$. 

Clearly $H^{0}(X,i_{*}\ZZ)=H^0(\{x\},\ZZ)=\ZZ$, and since $\CComplex_{P}(i_{*}\ZZ)$ is the cohomology complex for 
$P_{x}$ with coefficients in $\ZZ$, we also have $\HP^{0}(i_{*}\ZZ)=H^0(P_{x},\ZZ)$.  Therefore in order to have 
$\HP^{0}(i_{*}\Fsh)=H^0(X,i_{*}\ZZ)$, we must have $H^0(P_{x},\ZZ)=\ZZ$.
\epf

\bpoint{Lemma}\label{lem:Godement-zero-in-reduced-cech} Let $\Fsh$ be a sheaf of abelian groups on $X$, 
and $I_{\Fsh}$ the corresponding flasque Godement sheaf, i.e., the sheaf defined \cite[\S4.3, p.\ 167]{G} for
all open sets $U\subseteq X$ by 

\begin{equation}\label{eqn:Godement-def}
I_{\Fsh}(U) = \prod_{x\in U} \Fsh_{x},
\end{equation}

where $\Fsh_{x}$ denotes the stalk of $\Fsh$ at $x$. 
Assume that $(P,\psi)$ satisfies (P4).  Then $$\HP^{i}(I_{\Fsh})=0\,\mbox{for all $i\geq 1$}.$$

\bpf By \eqref{eqn:Godement-def} for each $\sigma\in P$ we have 
$I_{\Fsh}(U_{\sigma}) = \prod_{x\in U_{\sigma}} \Fsh_{x}$.  
All the maps in $\CComplex_{P}(I_{\Fsh})$ act componentwise with respect this product structure. 
It follows that 

\begin{equation}\label{eqn:cohomology-of-Godement-sheaf}
\HP^{i}(I_{\Fsh}) = \prod_{x\in X} H^{i}(P_{x},\Fsh_{x}),
\end{equation}

where $H^{i}(P_{x},\Fsh_{x})$ denotes the cohomology of the topological space $P_{x}$ with coefficients in the 
abelian group $\Fsh_{x}$.  If $P$ satisfies (P4) then by the universal coefficient theorem the right hand side of 
\eqref{eqn:cohomology-of-Godement-sheaf} is zero for all $i\geq 1$.  \epf

\bpoint{Definition of $\HP^{i,j}(\Fsh)$}  
\label{sec:def-of-Hij}
Extending the construction in \S\ref{sec:Def-of-reduced-complex}, given a sheaf $\Fsh$
of abelian groups on $X$, for each $i$, $j\in \NN$ we put 
$$\CComp_{P}^{i,j}(\Fsh) := \prod_{\sigma\in P,\,|\sigma|=i} H^{j}(U_{\sigma},\Fsh).$$
If $\sigma\subseteq \sigma'$, then $U_{\sigma'} \subseteq U_{\sigma}$, and we have maps 
$H^{j}(U_{\sigma},\Fsh)\longrightarrow H^{j}(U_{\sigma'},\Fsh)$ induced by restriction.  
For each fixed $j$  the formula for 
$d_{\sigma,\sigma'}$ in \eqref{eqn:def-of-differential}, where we interpret $\rho_{U_{\sigma},U_{\sigma'}}$ as
the map on cohomology above, allows us as before to construct maps 
$d^{i}\colon \CComp_{P}^{i,j}(\Fsh) \longrightarrow \CComp_{P}^{i+1,j}(\Fsh)$ 
giving rise to a complex $\CComp_{P}^{\smallbullet,j}(\Fsh)$.
We then define cohomology groups
$\HP^{i,j}(\Fsh):=H^{i}(\CComp_{P}^{\smallbullet,j}(\Fsh))$. 
The constructions of the complex and the groups are again functorial in $\Fsh$. 

To compare with the previous notation, $\CComp_{P}^{i}(\Fsh)=\CComp_{P}^{i,0}(\Fsh)$, 
$\CComplex_{P}(\Fsh)=\CComp_{P}^{\smallbullet,0}$, and $\HP^{i}(\Fsh)=\HP^{i,0}(\Fsh)$. 
Since we are mostly interested in the case $j=0$, we will largely use the simpler notation $\HP^{i}(\Fsh)$.
Our main use of the other groups is in connection with the result below.

\bpoint{Theorem}
\label{thm:reduced-cech-spectral-sequence}
Let $X$ be a topological space with the data of a cell complex $P$ and assignment $\psi$ as 
in \S\ref{sec:Def-of-reduced-complex}, and assume that $(P,\psi)$ satisfies (P1)--(P4).
Then for any sheaf $\Fsh$ of abelian groups on $X$ there is a spectral sequence 

\begin{equation}\label{eqn:spectral-sequence}
E_2^{i,j}:=\HP^{i,j}(\Fsh) \implies H^{i+j}(X,\Fsh).
\end{equation}

I.e., a first-quadrant spectral sequence with $E_2$ term $E_2^{i,j}=\HP^{i,j}(\Fsh)$, abutting to the
sheaf cohomology of $\Fsh$.

\bpf 
Let 
$$0 \longrightarrow \Fsh \longrightarrow \Ish^{0} 
\longrightarrow \Ish^{1} \longrightarrow \Ish^{2} \longrightarrow \Ish^{3} \longrightarrow \cdots $$
be the canonical resolution of $\Fsh$ by Godement sheaves \cite[\S4.3, p.168]{G}
(beginning with $\Ish^{0}=I_{\Fsh}$ and injection $\Fsh\hookrightarrow \Ish^{0}$, one sets
$\Ish^{1}=I_{\Ish^{0}/\Fsh}$, and continues inductively) and consider the complex

\begin{equation}\label{eqn:Godement-complex}
0 \longrightarrow \Ish^{0} \longrightarrow \Ish^{1} \longrightarrow \Ish^{2} \longrightarrow \Ish^{3} 
\longrightarrow \cdots
\end{equation}

with each $\Ish^{j}$ placed in degree $j$.  The crucial property of \eqref{eqn:Godement-complex} is
that for any open $U\subseteq X$, the corresponding complex of global sections 

\begin{equation}\label{eqn:Godement-complex-U}
0 \longrightarrow \Gamma(U,\Ish^{0}) \longrightarrow \Gamma(U,\Ish^{1}) \longrightarrow \Gamma(U,\Ish^{2}) 
\longrightarrow \Gamma(U,\Ish^{3}) \longrightarrow \cdots
\end{equation}

computes $H^{\smallbullet}(U,\Fsh)$. (In \cite[\S4.4, p.\ 173]{G} the cohomology of \eqref{eqn:Godement-complex-U}
is taken as the {definition} of $H^{\smallbullet}(U,\Fsh)$. For the construction of sheaf cohomology by
other approaches, the equality follows since \eqref{eqn:Godement-complex} is a flasque resolution of $\Fsh$.)
Applying $\CComp_{P}^{\smallbullet}(\cdot)$ to 
\eqref{eqn:Godement-complex}, we obtain, by functoriality of the construction (\S\ref{sec:functoriality-of-complex}),
a complex of complexes, i.e., a double complex~:

\neweqnum
\hspace{1cm}
\label{eqn:double-complex}
\begin{tabular}{c}
\psset{xunit=2.5cm, yunit=1.5cm,nodesep=0.2cm}
\begin{pspicture}(-0.5,-0.5)(4.4,3.5)
\psset{linestyle=dashed,linecolor=lightgray}
\psline(-0.5,-0.5)(4.3,-0.5)
\psline(-0.5,0.5)(4.3,0.5)
\psline(-0.5,1.5)(4.3,1.5)
\psline(-0.5,2.5)(4.3,2.5)
\psline(-0.5,-0.5)(-0.5,3.3)
\psline(0.5,-0.5)(0.5,3.3)
\psline(1.5,-0.5)(1.5,3.3)
\psline(2.5,-0.5)(2.5,3.3)
\psline(3.5,-0.5)(3.5,3.3)
\psset{linestyle=solid,linecolor=black}
\rput(0,0){\rnode{00}{$\CComp_{P}^{0}(\Ish^{0})$}}
\rput(1,0){\rnode{10}{$\CComp_{P}^{0}(\Ish^{1})$}}
\rput(2,0){\rnode{20}{$\CComp_{P}^{0}(\Ish^{2})$}}
\rput(3,0){\rnode{30}{$\CComp_{P}^{0}(\Ish^{3})$}}
\rput(0,1){\rnode{01}{$\CComp_{P}^{1}(\Ish^{0})$}}
\rput(1,1){\rnode{11}{$\CComp_{P}^{1}(\Ish^{1})$}}
\rput(2,1){\rnode{21}{$\CComp_{P}^{1}(\Ish^{2})$}}
\rput(3,1){\rnode{31}{$\CComp_{P}^{1}(\Ish^{3})$}}
\rput(0,2){\rnode{02}{$\CComp_{P}^{2}(\Ish^{0})$}}
\rput(1,2){\rnode{12}{$\CComp_{P}^{2}(\Ish^{1})$}}
\rput(2,2){\rnode{22}{$\CComp_{P}^{2}(\Ish^{2})$}}
\rput(3,2){\rnode{32}{$\CComp_{P}^{2}(\Ish^{3})$}}
\rput(0,3){\rnode{03}{$\vdots$}}
\rput(1,3){\rnode{13}{$\vdots$}}
\rput(2,3){\rnode{23}{$\vdots$}}
\rput(3,3){\rnode{33}{$\vdots$}}
\rput(3.9,0){\rnode{40}{$\cdots$}}
\rput(3.9,1){\rnode{41}{$\cdots$}}
\rput(3.9,2){\rnode{42}{$\cdots$}}
\psset{arrows=->}
\ncline{00}{10}
\ncline{10}{20}
\ncline{20}{30}
\ncline{30}{40}
\ncline{01}{11}
\ncline{11}{21}
\ncline{21}{31}
\ncline{31}{41}
\ncline{02}{12}
\ncline{12}{22}
\ncline{22}{32}
\ncline{32}{42}
\ncline{00}{01}
\ncline{01}{02}
\ncline{02}{03}
\ncline{10}{11}
\ncline{11}{12}
\ncline{12}{13}
\ncline{20}{21}
\ncline{21}{22}
\ncline{22}{23}
\ncline{30}{31}
\ncline{31}{32}
\ncline{32}{33}
\end{pspicture}
\end{tabular}.

\vspace{0.25cm}
Comparing the two spectral sequences for the associated total complex gives the result.
We now check the details. 

If we start by first taking the cohomology of the vertical maps, then by Lemmas 
\ref{lem:reduced-H0-correct} and \ref{lem:Godement-zero-in-reduced-cech} the double complex is reduced to 
$$
0 \longrightarrow \Gamma(X,\Ish^{0}) \longrightarrow \Gamma(X,\Ish^{1}) \longrightarrow \Gamma(X,\Ish^{2}) 
\longrightarrow \Gamma(X,\Ish^{3}) \longrightarrow \cdots.
$$
I.e., we are reduced to the setting of \eqref{eqn:Godement-complex-U} with $U=X$.  The cohomology of this complex
is therefore $H^{\smallbullet}(X,\Fsh)$, showing that the double complex computes the cohomology of $X$ with 
coefficients in $\Fsh$. 

On the other hand, the horizontal maps in \eqref{eqn:double-complex} are, by the construction in 
\S\ref{sec:functoriality-of-complex} the product of maps of the form 
$$\cdots\longrightarrow 
\Gamma(U_{\sigma},\Ish^{j-1})\longrightarrow 
\Gamma(U_{\sigma},\Ish^{j})\longrightarrow 
\Gamma(U_{\sigma},\Ish^{j+1})\longrightarrow \cdots
$$
i.e., products of maps of the type \eqref{eqn:Godement-complex-U} with $U=U_{\sigma}$. 
Taking cohomology, the result for a fixed row $i$ is the product, over those $\sigma$ with $|\sigma|=i$, 
of the cohomology groups $H^{j}(U_{\sigma},\Fsh)$. This product is $\CComp_{P}^{i,j}(\Fsh)$, by definition
of the latter (\S\ref{sec:def-of-Hij}).  
The $E_{2}$ term of this spectral sequence, resulting from taking the cohomology of the columns (i.e., with $j$
fixed) is, in row $i$ and column $j$, the group $\HP^{i,j}(\Fsh)$,
again by the definition in \S\ref{sec:def-of-Hij}.  
This completes the proof of the theorem.  \epf

%

\bpoint{Definition of $\Ab_{P}$} 
\label{sec:Def-of-AbP}
Given a topological space $X$, and given $(P,\psi)$ as in 
\S\ref{sec:Def-of-reduced-complex} we set $\Ab_{P}$ to be full subcategory of the category of sheaves of abelian 
groups on $X$ consisting of those sheaves $\Fsh$ so that 

\begin{equation}\label{eqn:Def-of-Ab}
H^{j}(U_{\sigma},\Fsh)=0\rule{0.25cm}{0cm}\mbox{for all $\sigma\in P$ and all $j\geq 1$.}
\end{equation}

If $X$ is a scheme, and $\{U_{\sigma}\}_{\sigma\in P}$ a cover of $X$ by affine open sets, then all 
quasi-coherent sheaves on $X$ are in $\Ab_{P}$. 

\bpoint{Theorem} 
\label{thm:reduced-cech-computes-cohomology}
Let $X$ be a topological space with the data of a cell complex $P$ and assignment $\psi$ as described
in \S\ref{sec:Def-of-reduced-complex}, and assume that $(P,\psi)$ satisfies (P1)--(P4).
Then for every sheaf $\Fsh$ in $\Ab_{P}$, the reduced \v{C}ech complex associated to $P$ computes the cohomology 
of $\Fsh$.  I.e., for every $\Fsh\in \Ab_{P}$, $\HP^{i}(\Fsh)\cong H^{i}(X,\Fsh)$ for all $i\geq 0$. 

\bpf
By the definition of $\Ab_{P}$ and of $\CComp_{P}^{i,j}(\Fsh)$, if $\Fsh\in \Ab_{P}$ then 
$\CComp_{P}^{i,j}(\Fsh)=0$ for all $i$ and all $j\geq 1$.
Therefore $\HP^{i,j}(\Fsh)=0$ for all $i$ and all $j\geq 1$, and so the spectral sequence \eqref{eqn:spectral-sequence}
in Theorem \ref{thm:reduced-cech-spectral-sequence} degenerates immediately, giving
$H^{i}(X,\Fsh)=\HP^{i,0}(\Fsh)=\HP^{i}(\Fsh)$ for all $i$.
\epf

\bpoint{Corollary} 
\label{cor:reduced-cech-complex-on-a-subspace}
Let $Z\subseteq X$ be a subset.  For each $\sigma\in P$ set $V_{\sigma} := Z\cap U_{\sigma}$.
Then the mapping $\sigma\to V_{\sigma}$ satisfies (P1)--(P4) as a covering of $Z$.   Thus by 
Theorem \ref{thm:reduced-cech-computes-cohomology}, for all $\Fsh\in \Ab_{P,Z}$,
$\HP^{i}(\Fsh)\cong H^{i}(Z,\Fsh)$ for all $i\geq 0$.  Here by $\Ab_{P,Z}$, we mean the category
of sheaves as in \S\ref{sec:Def-of-AbP}, but with respect to the covering $\{V_{\sigma}\}_{\sigma\in P}$ of $Z$. 

\bpf
For each $x\in Z$
the cell-complexes $P_{x}$ are the same whether using the covering 
$\{U_{\sigma}\}_{\sigma \in P}$ of $X$ or the covering
$\{V_{\sigma}\}_{\sigma \in P}$ of $Z$.  Thus (P3) and (P4) for $X$ imply (P3) and (P4) for $Y$. 
The implication that (P1) and (P2) for $X$ imply the corresponding conditions for $Z$ is also elementary. 
\epf

\point {\bf Remark}. 
\label{sec:remark-P4-necessary}
Theorem \ref{thm:reduced-cech-computes-cohomology} shows that (P1)--(P4) are sufficient conditions so that
$\HP^{i}(\Fsh)= H^{i}(X,\Fsh)$ for all $\Fsh\in \Ab_{P}$ and all $i\geq 0$.  Conditions (P1) and (P2) are necessary
to even define the groups $\HP^i$.  Lemma \ref{lem:reduced-H0-correct} shows that (P3) is a necessary condition
to have equality of $\HP^{0}$.  

To see that (P4) is also a necessary condition, we consider $I_{\ZZsheaf}$, the Godement flasque sheaf associated
to the locally constant sheaf $\Fsh=\ZZsheaf$.   Since $I_{\ZZsheaf}$ is flasque, $H^{i}(U,I_{\ZZsheaf})=0$
for all $i\geq 1$, and all open subsets $U\subseteq X$.  In particular, $I_{\ZZsheaf}\in \Ab_{P}$, and 
$H^{i}(X,I_{\ZZsheaf})=0$ for all $i\geq 0$. 

For this sheaf, \eqref{eqn:cohomology-of-Godement-sheaf} becomes 
$$\HP^{i}(I_{\ZZsheaf}) = \prod_{x\in X} H^{i}(P_{x},\ZZ).$$

Thus, in order to have $\HP^{i}(I_{\ZZsheaf})=H^{i}(X,I_{\ZZsheaf})=0$ for all $i\geq 0$, we must have 
$H^{i}(P_{x},\ZZ)=0$ for all $i\geq 0$ and all $x\in X$. 

\bpoint{Computing with a sub cell complex} 
\label{sec:sub-cell-complex}
Let $X$ and $(P,\psi)$ be as in Theorem \ref{thm:reduced-cech-computes-cohomology}, with $(P,\psi)$ satisfying
(P1)--(P4).  Let $Q\subseteq P$ be a sub cell complex, and 
let $\psi|_{Q}$ denote the restriction of the assignment $\psi$ to $Q$. 

Since $(P,\psi)$ satisfies (P1), $(Q,\psi|_{Q})$ does as well.  If $(Q,\psi|_{Q})$ also satisfies (P2), 
then for each sheaf $\Fsh$ of abelian groups on $X$ we may construct 
groups $\HQ^{i}(\Fsh)$ as in \S\ref{sec:Def-of-reduced-complex}--\ref{sec:Def-of-HP}. 
Furthermore, $\CComp_{Q}^{\smallbullet}(\Fsh)$ is then a quotient complex
of $\CComp_{P}^{\smallbullet}(\Fsh)$, and so for each $i$ we have a map $\HP^{i}(\Fsh)\longrightarrow \HQ^{i}(\Fsh)$. 

The cover $\{U_{\sigma}\}_{\sigma\in Q}$ is a subcover of $\{U_{\sigma}\}_{\sigma\in P}$, and thus 
$\Ab_{P}\subseteq \Ab_{Q}$.   If $(Q,\psi|_{Q})$ also satisfies (P3) and (P4), 
then by Theorem \ref{thm:reduced-cech-computes-cohomology} for each $\Fsh\in \Ab_{P}$ and each $i\geq 0$ the groups
$\HP^{i}(\Fsh)$ and $\HQ^{i}(\Fsh)$ are isomorphic, since each is isomorphic to $H^{i}(X,\Fsh)$. 

\bpoint{Lemma} \label{lem:pass-to-subpolytope} 
Let $X$ and $(P,\psi)$ be as in Theorem \ref{thm:reduced-cech-computes-cohomology}, with $(P,\psi)$ satisfying
(P1)--(P4).  Let $Q\subseteq P$ be a sub cell complex, and assume that $(Q,\psi|_{Q})$ also satisfies
(P1)--(P4).   Then for every $\Fsh\in \Ab_{P}$, and all $i\geq 0$, the map 
$\HP^{i}(\Fsh)\longrightarrow \HQ^{i}(\Fsh)$ induced by the quotient map 
$\CComp_{P}^{\smallbullet}(\Fsh)\longrightarrow \CComp_{Q}^{\smallbullet}(\Fsh)$ is an isomorphism.

\bpf We first consider the case that $\Ish=I_{\Gsh}$ is the flasque Godement sheaf associated to a sheaf 
of abelian groups $\Gsh$ on $X$.   Then by
Lemma \ref{lem:Godement-zero-in-reduced-cech} both $\HP^{i}(\Ish)$
and $\HQ^{i}(\Ish)$ are zero when $i\geq 1$.  Using \eqref{eqn:cohomology-of-Godement-sheaf} the map on 
$\check{H}^{0}$ induced by the quotient map on complexes is the componentwise map 

\vspace{-0.2cm}
\begin{equation}
\label{eqn:H0-subcomplex-restriction}
\HP^{0}(\Ish) = \prod_{x\in X} H^0(P_{x},\Gsh_{x})\longrightarrow 
\prod_{x\in X} H^0(Q_{x},\Gsh_{x}) = \HQ^{0}(\Ish).
\end{equation}

\vspace{-0.2cm}
Thus the issue is to show that each map $H^0(P_{x},\Gsh_{x})\longrightarrow H^0(Q_{x},\Gsh_{x})$ is
an isomorphism.   Elements of $H^0(P_{x},\Gsh_{x})$ are functions on the vertices of $P_{x}$, 
taking values in $\Gsh_{x}$, and having the same values on vertices connected by edges of $P_{x}$. 
A similar description holds for $H^0(Q_{x},\Gsh_{x})$.  

Since the assignment function for $Q$ is the restriction of that for $P$, 
$\sigma\in Q_{x}$ implies that $\sigma\in P_{x}$, and thus $Q_{x}\subseteq P_{x}$. 
The map 
$H^0(P_{x},\Gsh_{x})\longrightarrow H^0(Q_{x},\Gsh_{x})$ is the restriction of a function on the vertices of $P_{x}$
to the vertices of $Q_{x}$.

By (P3), any two vertices of $P_{x}$ are connected by a chain of edges. 
Thus $H^0(P_x,\Gsh_{x})\cong \Gsh_{x}$, with the isomorphism being given by evaluating on any vertex of $P_{x}$. 
Similarly $H^0(Q_x,\Gsh_{x})\cong \Gsh_{x}$.  Since $Q_{x}\subseteq P_{x}$, in both cases we may choose to evaluate
on a vertex in $Q_{x}$.  Thus the map $H^0(P_{x},\Gsh_{x})\longrightarrow H^0(Q_{x},\Gsh_{x})$ is an 
isomorphism.

In the general case, we take the complex \eqref{eqn:Godement-complex} of Godement sheaves associated to $\Fsh$, and
obtain double complexes, one each for $P$ and $Q$ as in \eqref{eqn:double-complex}.  The quotient map 
$\CComp_{P}^{\smallbullet}\longrightarrow \CComp_{Q}^{\smallbullet}$ induces a map of double complexes, 
and a corresponding map of total complexes. 

The total complexes have two associated spectral sequences.  In the second (in the order used in the proof of 
Theorem \ref{thm:reduced-cech-spectral-sequence}) the map on $E_2$ terms is the collection of maps
$\HP^{i}(\Fsh)\longrightarrow\HQ^{i}(\Fsh)$ we are interested in.   
This second spectral sequence degenerates at the $E_2$ term when $\Fsh\in \Ab_{P}$, and so it suffices to show
that the map on total complexes induces an isomorphism in cohomology. 

For that we may use the first spectral sequence. The map on $E_1$ terms of that sequence is 

\begin{equation}
\label{eqn:gamma-identity-diagram}
\rule{1cm}{0cm} 
\begin{array}{ccccccccccc}
0 & \longrightarrow & \HP^{0}(\Ish^{0}) & \longrightarrow & \HP^{0}(\Ish^{1}) & \longrightarrow & \HP^{0}(\Ish^{2}) 
& \longrightarrow & \HP^{0}(\Ish^{3}) & \longrightarrow & \cdots \\
& & \xdownarrow{0.35cm} 
& & \xdownarrow{0.35cm} 
& & \xdownarrow{0.35cm} 
& & \xdownarrow{0.35cm} 
\\
0 & \longrightarrow & \HQ^{0}(\Ish^{0}) & \longrightarrow & \HQ^{0}(\Ish^{1}) & \longrightarrow & \HQ^{0}(\Ish^{2}) 
& \longrightarrow & \HQ^{0}(\Ish^{3}) & \longrightarrow & \cdots \\
\end{array}
\end{equation}

\vspace{-0.15cm}
and we have already
shown that these maps are isomorphisms when $\Ish$ is a flasque Godement sheaf. 
This finishes the proof of the lemma. \epf

\point {\bf Remark}. The proof of the lemma also shows the more precise statement that, under the identifications 
$\HP^{i}(\Fsh)=H^{i}(X,\Fsh)$ and $\HQ^{i}(\Fsh)=H^{i}(X,\Fsh)$ provided by 
Theorem \ref{thm:reduced-cech-computes-cohomology}, the map $\HP^{i}(\Fsh)\longrightarrow \HQ^{i}(\Fsh)$ 
is the identity.

To see this, we first note that the diagram 

\begin{centering}
\begin{tabular}{c}
\begin{pspicture}(0,0)(3,2)
\psset{nodesep=2mm,arrows=->}
\rput(0,1.5){\rnode{HP}{$H^{0}(P_{x},\Gsh_{x})$}}
\rput(2.5,1.5){\rnode{GP}{$\Gsh_{x}$}}
\rput(0,0){\rnode{HQ}{$H^{0}(Q_{x},\Gsh_{x})$}}
\rput(2.5,0){\rnode{GQ}{$\Gsh_{x}$}}
\ncline{HP}{HQ}
\ncline{HP}{GP}
\ncline{HQ}{GQ}
\rput(2.5,0.75){\rotateright{$\xeq{\rule{0.75cm}{0cm}}$}}
\rput(1.55,1.7){\tiny$\sim$}
\rput(1.55,0.2){\tiny$\sim$}
\rput(0.2,0.75){\tiny\rotateright{$\sim$}}
\rput(3.0,0.75){\tiny $\operatorname{Id}_{\Gsh_x}$} 
\end{pspicture}
\end{tabular}\\
\end{centering} 

\vspace{0.25cm} 
commutes since the horizontal isomorphisms are given by evaluating at a vertex in $P_x$ or $Q_x$ respectively, and
the vertical left hand isomorphism is given by restricting a function on the vertices of $P_{x}$ to the vertices
of $Q_{x}$. 

Thus (under the identification of each $\check{H}^{0}$ with $\prod \Gsh_x$)  the map in 
\eqref{eqn:H0-subcomplex-restriction} is the identity, and
so the vertical maps in \eqref{eqn:gamma-identity-diagram} are also the identity, and so induce the identity map
$H^{i}(X,\Fsh)\longrightarrow H^{i}(X,\Fsh)$ on cohomology.

But, the maps $\HP^{i}(\Fsh)\longrightarrow \HQ^{i}(\Fsh)$ are (under the identifications with $H^{i}(X,\Fsh)$)
exactly these maps. (The
identification of $\HP^{i}(\Fsh)$ and $H^{i}(X,\Fsh)$ is through the degeneration of the two spectral sequences
associated to the double complexes, and \eqref{eqn:gamma-identity-diagram} is the map on $E_1$ terms of the spectral
sequence whose $E_2$ term is the $H^{i}(X,\Fsh)$.) \epf

We will not need the result below, but record it for completeness. 

\bpoint{Lemma} \label{lem:Ab-lemma}
Let 
$$0\longrightarrow \Fsh \longrightarrow \Gsh \longrightarrow \Hsh \longrightarrow 0$$
be a short exact sequence of sheaves of abelian groups on $X$.

\begin{enumerate}
\item If $\Fsh\in \Ab_{P}$ then we have a long exact sequence of the $\HP^{\smallbullet}$~:
$$ \cdots \longrightarrow \HP^{i-1}(\Hsh) \longrightarrow \HP^{i}(\Fsh)\longrightarrow \HP^{i}(\Gsh)\longrightarrow 
\HP^{i}(\Hsh)\longrightarrow \HP^{i+1}(\Fsh)\longrightarrow \cdots$$
\item If $\Fsh$, $\Hsh\in \Ab_{P}$ then $\Gsh\in \Ab_{P}$.
\item If $\Fsh$, $\Gsh\in \Ab_{P}$ then $\Hsh\in \Ab_{P}$.
\end{enumerate}

\bpf
({\em a}) If $\Fsh\in \Ab_{P}$ then for each $\sigma\in P$, since $H^{1}(U_{\sigma},\Fsh)=0$ we have a short
exact sequence 
$$0\longrightarrow \Gamma(U_{\sigma},\Fsh) \longrightarrow \Gamma(U_{\sigma},U_{\sigma}) \longrightarrow 
\Gamma(U_{\sigma},\Hsh) \longrightarrow 0,$$
and therefore a short exact sequence of complexes 
$$0\longrightarrow \CComplex_{P}(\Fsh) \longrightarrow \CComplex_{P}(\Gsh) \longrightarrow \CComplex_{P}(\Hsh)
\longrightarrow 0,$$
giving the long exact sequence in cohomology. 
Parts (b) and (c) are clear from the long exact sequence of (usual) cohomology for each $U_{\sigma}$.
\epf

\point {\bf Remark}. If, in part (a) of the previous lemma, $\Fsh$, $\Gsh$, and $\Hsh$ are all in 
$\Ab_{P}$, then by Theorem \ref{thm:reduced-cech-computes-cohomology} the reduced \v{C}ech cohomology groups
and usual cohomology groups coincide.  Under those isomorphisms, the long exact sequence in (a) is, of course, the
long exact sequence of usual cohomology groups.   This can be seen by starting with a short exact
sequence of the corresponding double complexes, and considering the two associated spectral sequences. 

We finish off this section with examples addressing two natural questions associated to this construction.

\bpoint{Example} 
\label{ex:axioms-dont-force-values-on-faces}
In the standard \v{C}ech complex, and in the reduced \v{C}ech complexes 
constructed in Theorem \ref{thm:construction-of-toric-reduced-cover}, one feature of the assignment
$\psi$ is that for each $\sigma\in P$, $U_{\sigma}=\bigcap_{v\in \sigma} U_{v}$.  
I.e., each $U_{\sigma}$ is the intersection of the $U_{v}$ as $v$ runs over the vertices of $\sigma$.  One may
ask if this is forced by the conditions (P1)--(P4).  Here is a simple example to show that it is not, even
with the requirement that each $U_{\sigma}$ be affine, so that $\Ab_{P}$ contains the category of quasi-coherent 
sheaves. 

Let $X=\PP(\Osh_{\PP^1}\oplus \Osh_{\PP^1}(1))$ be the first Hirzebruch surface, and $\pi\colon X\longrightarrow \PP^2$
the blowdown map to $\PP^2$. 
We use the cell complex $P:=X_{\geq 0}$, which again, as in the case of $\PP^1\times\PP^1$,  is a square, with cells 
corresponding to the torus orbits on $X$.
We use $v_1$,\ldots, $v_4$ for the four torus-fixed points on $X$, $\tau_{ij}$ for the one dimensional orbits, 
with the indices $ij$ indicating that $v_i$ and $v_j$ are in the closure of these orbits, 
and $\sigma$ for the largest torus orbit, $T_{X}$.   
We choose the indices so that the closure of $O_{\tau_{12}}$ is the exceptional divisor of the blowdown map $\pi$. 

As in \S\ref{sec:first-example-of-reduced} (and Theorem \ref{thm:construction-of-toric-reduced-cover})
to each face $\gamma\in P$ we associate the affine open set 
$$U'_{\gamma} := \bigcup_{O_{\gamma}\subseteq \overline{O}_{\gamma'}} O_{\gamma'}.$$

With this cover, there are three possibilities for the topological type of $P_{x}$, depending on $x$.
If $x\in X$ is one of the four torus fixed points, $P_{x}$ is a single point. If $x$ is on one of the one-dimensional
torus orbits $P_{x}$ is a segment, homeomorphic to $[0,1]$.   If $x$ is in the two-dimensional orbit, i.e., in $T_{X}$,
then $P_{x}$ is all of $P$, and so is a square.  
This cover satisfies the condition that for all $\gamma\in P$, $U'_{\gamma}=\bigcap_{v\in \gamma} U'_{v}$.

We modify the cover as follows.  Fix a line $\ell\subset\PP^2$ not passing through any of the 
coordinate points.   For $\gamma\in P$, $\gamma$ a vertex, or $\gamma=\tau_{ij}$ with $\tau_{ij}\neq \tau_{12}$, we
set $U_{\gamma}=U'_{\gamma}$, and then set $U_{\tau_{12}}=U'_{\tau_{12}}\setminus\pi^{-1}(\ell)$ and
$U_{\sigma}=U'_{\sigma}\setminus\pi^{-1}(\ell)$.  Since $\ell$ does not pass through any of the coordinate points,
$\pi^{-1}(\ell)$ does not meet the exceptional divisor of the blowdown (nor any of the torus fixed points),
and so no point on the exceptional divisor is removed by this process.  When $\gamma$ is $\tau_{12}$ or $\sigma$, 
we have removed a Cartier divisor from $U'_{\gamma}$, so the result is still affine. Thus, this is a cover
by affine opens.  This cover again satisfies (P1) and (P2).

To see that the cover satisfies (P3) and (P4), we note that for $x$ a torus fixed point and $x$ in a one-dimensional
torus orbit, the description of $P_{x}$ is as above.  For $x\in T_{X}$ there are two possibilities,
$x\in T_{X}\cap\pi^{-1}(\ell)$ or $x\in T_{X}\setminus \pi^{-1}(\ell)$. In the first case, $P_{x}$ is the boundary
of the square, minus the edge corresponding to the exceptional divisor, and in the second $P_{x}$ is again a square.
These possibilities are illustrated in Figure \ref{fig:possible-Pxs} below.

\vspace{0.5cm}
\begin{centering} 
\begin{tabular}{|c|c|c|c|}
\hline
& & & \\
\begin{tabular}{c}
\begin{pspicture}(-0.5,-0.5)(0.5,0.5)
\psset{fillstyle=solid,linecolor=black}
\psset{fillcolor=white}
\pscircle(-0.5,-0.5){0.05}
\end{pspicture}
\end{tabular}
&
\begin{tabular}{c}
\begin{pspicture}(-0.5,-0.5)(0.5,0.5)
\psset{fillstyle=solid,linecolor=black}
\psline(-0.5,-0.5)(0.5,-0.5)
\psset{fillcolor=white}
\pscircle(-0.5,-0.5){0.05}
\pscircle(0.5,-0.5){0.05}
\end{pspicture}
\end{tabular}
&
\begin{tabular}{c}
\begin{pspicture}(-0.5,-0.5)(0.5,0.5)
\psset{fillstyle=solid,linecolor=black}
\psline[fillstyle=none](-0.5,0.5)(-0.5,-0.5)(0.5,-0.5)(0.5,0.5)
\psset{fillcolor=white}
\pscircle(-0.5,-0.5){0.05}
\pscircle(0.5,-0.5){0.05}
\pscircle(0.5,0.5){0.05}
\pscircle(-0.5,0.5){0.05}
\end{pspicture}
\end{tabular}
&
\begin{tabular}{c}
\begin{pspicture}(-0.5,-0.5)(0.5,0.5)
\psset{fillstyle=solid,linecolor=black}
\pspolygon[fillcolor=vlgray](-0.5,-0.5)(-0.5,0.5)(0.5,0.5)(0.5,-0.5)
\psset{fillcolor=white}
\pscircle(-0.5,-0.5){0.05}
\pscircle(0.5,-0.5){0.05}
\pscircle(0.5,0.5){0.05}
\pscircle(-0.5,0.5){0.05}
\end{pspicture}
\end{tabular} \\
\small $x=v_i$ & \small $x\in O_{\tau_{ij}}$  & \small $x\in\pi^{-1}(\ell)\cap T_{X}$ 
& \small $x\in T_{X}\setminus\pi^{-1}(\ell)$ \rule[-0.3cm]{0cm}{0.8cm}\\
\hline
\multicolumn{4}{c}{} \\
\multicolumn{4}{c}{\Fig\label{fig:possible-Pxs} Possibilities for $P_x$} \\
\end{tabular}\\
\end{centering}

\vspace{0.5cm}
In all cases $P_{x}$ has the cohomology of a point, and thus this cover satisfies (P3) and (P4). 
For $\gamma=\tau_{12}$ or $\gamma=\sigma$, $U_{\gamma}$ is a proper subset of $\bigcap_{v\in \gamma} U_{\gamma}$.


\point 
\label{sec:comparison-argument-suggestion}
Given a topological space $X$, cell complex $P$, and an assignment $\psi$ satisfying (P1)--(P4), one might
ask if it is possible to prove the result of Theorem \ref{thm:reduced-cech-computes-cohomology} 
by using the total \v{C}ech complex $\CComp^{\smallbullet}(\Fsh)$ associated to the 
cover $\{U_{v}\}_{v\in P(0)}$, finding a morphism between the complexes 
$\CComp^{\smallbullet}(\Fsh)$ and $\CComp^{\smallbullet}_{P}(\Fsh)$  and proving that the cohomology of the 
mapping cone is zero.

The author does not see how to show that this is impossible, but also does not see how to carry out such a plan.
Among the difficulties are : (\Euler{1}) It is not clear how to give a general construction of a map of complexes
between $\CComp(\Fsh)$ and $\CComp_{P}(\Fsh)$.   (\Euler{2}) It is also not clear how one would prove that the
cohomology of the resulting mapping cone is zero.  (\Euler{3}) Implicit in this question
is that the cohomology of the total \v{C}ech complex should compute the actual sheaf cohomology of $\Fsh$.  For this
to happen, the usual requirement is that the higher cohomology of $\Fsh$ should be zero on all the open sets
which appear in the total \v{C}ech complex.  However, Theorem \ref{thm:reduced-cech-computes-cohomology}
applies to sheaves $\Fsh\in\Ab_{P}$, i.e., to the category of sheaves whose higher cohomology vanishes on the
open sets $U_{\sigma}$, $\sigma\in P$, and this may be a different category than that of the sheaves whose
cohomology is computed by the total \v{C}ech complex.

The following is a simple example of the third issue.

\bpoint{Example} 
\label{Ex:comparison-argument-example}
Let $X$ be the affine plane with doubled origin, i.e,. the nonseparated scheme obtained by
taking two copies of $\AA^2$ and glueing them together by the identity map everywhere but the origin of each of the
two planes.

Let $P$ be a square, with vertices $v_{1}$,\ldots, $v_{4}$, edges $\tau_{ij}$, and $2$-cell $\sigma$, organized
as in Figure \ref{fig:P1xP1}.  To $v_1$ we associate the open subset $\AA^2$ obtained by removing one origin,
to $v_4$ the open subset obtained by removing the other origin, again equal to $\AA^2$, to $v_2$ we associate
$\AA^2$ minus the $x$-axis (an open subset of both the original copies of $\AA^2$) and to $v_3$ we associate
$\AA^2$ minus the $y$-axis.  For the other $\gamma\in P$ we define $U_{\gamma}$ by $U_{\gamma}=\bigcup_{v\in \gamma} U_{v}$. 

By this rule both $U_{\tau_{12}}$ and $U_{\tau_{24}}$ are equal to $U_{v_2}$ ($\AA^2\setminus\{\mbox{$x$-axis}\}$),
similarly $U_{\tau_{13}}=U_{\tau_{34}}=U_{v_{3}}=\AA^2\setminus\{\mbox{$y$-axis}\}$, and finally 
$U_{\sigma}$ is $\AA^2$ minus both the $x$- and $y$-axes.  Thus, all $U_{\gamma}$ are affine varieties.

This assignment satisfies (P1)--(P4).  Conditions (P1) and (P2) are easy to see.  For (P3) and (P4), if $x$ is
one of the two copies of the origin, $P_{x}$ is a single point (the corresponding vertex of $P$), if $x$ is on
the $x$- or $y$-axis, but not on an origin $P_{x}$ is an L, consisting of three vertices and two connecting edges.  
For instance, if $x$ is on the $y$-axis, but not at an origin $P_{x}$ consists of the edges $\tau_{13}$ and $\tau_{34}$,
and corresponding vertices $v_1$, $v_3$, and $v_4$.  For $x$ off the $x$- and $y$-axes $P_{x}$ is the square.  In
all cases $P_{x}$ has the cohomology of a point. 

By Theorem \ref{thm:reduced-cech-computes-cohomology}, the reduced \v{C}ech complex associated to the cover
computes the cohomology of all quasi-coherent sheaves on $X$.

But, if we construct the associated total \v{C}ech complex, we will need to include the open set 
$U_{v_1}\cap U_{v_4}$ which is $\AA^2$ minus the origin.  Since this set is not an affine scheme, there
are coherent sheaves with nonzero higher cohomology on this open set. Thus, the associated total \v{C}ech cover
does not automatically compute the cohomology of coherent sheaves on $X$.  In particular, it is even less 
clear how a comparison argument of the type suggested in \S\ref{sec:comparison-argument-suggestion} would work
in this example. 

In the case of the reduced \v{C}ech complex associated to a simplicial projective toric variety $X$, there
is such a conceptual explanation.  Let $I$ be the ideal of the GIT unstable locus $Z\subseteq \AA^{s}$
(in the notation of \S\ref{sec:Oda-construction}), and
$S$ the coordinate ring of $\AA^{s}$, i.e., the Cox ring of $X$.  
By \cite[Proposition 6.10]{Mi1} the reduced complex and the standard complex computing the resolution of $S/I$ 
are quasi-isomorphic.  
From this one can deduce that the reduced \v{C}ech complex computes (for coherent sheaves) the same groups
as the standard \v{C}ech complex.

\section{Torus orbits and reduced \v{C}ech covers of semi-proper toric varieties}
\label{sec:covers-of-semi-proper}

\point 
In this section we show, for a semi-proper toric variety $X$,  
how to construct a pair $(P,\psi)$ satisfying (P1)--(P4) of \S\ref{sec:reduced-cech-complexes} such that

\RomanList
\begin{enumerate}
\item each $U_{\sigma}\subseteq X$ is a torus-stable open affine, and 
\item $P$ is a cell complex of dimension $\cdim(X)$.

\end{enumerate}
\AlphaList

Here $\cdim(X)$ is the cohomological dimension of $X$, 
the smallest $n$ so that $H^{i}(X,\Fsh)=0$ for all coherent sheaves $\Fsh$ on $X$, and all $i>n$.

A semi-proper variety is one for which the natural morphism $X\longrightarrow \Spec(\Gamma(X,\Osh_{X}))$ is proper.
This class includes proper varieties, and more generally varieties of the form $\phi^{-1}(V)$ where 
$\phi\colon X\longrightarrow Y$ is a proper morphism and $V\subseteq Y$ is affine.

If one is willing to drop condition ({\em ii}) then there is an easy solution~:
Take a cover of $X$ by torus-invariant affines and use the total \v{C}ech complex.  However (see example 2 in 
\S\ref{sec:first-example-of-reduced}) for $X=(\PP^1)^{n}$ the resulting $P$ has dimension $2^{n}-1$, while $\cdim(X)=n$.

On the other hand, for an $n$-dimensional projective $X$ (so $\cdim(X)=n$ again), if one is willing to drop 
condition ({\em i}), but still require the 
cover to be affine, then there is also an easy solution.  Embedding $X$ in projective space, and taking the 
affine cover of $X$ obtained as the complement of ${n+1}$ general hyperplanes,
the total \v{C}ech complex has associated polytope the $n$-dimensional simplex, the 
smallest $n$-dimensional polytope possible.
The disadvantage is that now the open sets are not torus-stable, and so one cannot use the torus action to organize
the computations. 

The case of semi-proper toric $X$ is one where we can achieve both goals, and where, moreover, 
the cell complex $P$ is directly tied to the structure of torus orbits on $X$.
This last feature will be useful in constructing an algorithm to compute higher direct images. 

In \S\ref{sec:orbits-and-affines}--\S\ref{prop:lining-up-of-torus-orbits} 
we review results about torus orbits, open affines on toric varieties, and toric morphisms, 
and prove some preliminary results.  The main result of the section is Theorem 
\ref{thm:construction-of-toric-reduced-cover}, constructing the pair $(P,\psi)$
when $X$ is toric and semi-proper (satisfying, of course, ({\em i}) and ({\em ii}) above, as well as (P1)--(P4)).


We begin by justifying one of the statements in the introductory paragraphs. 

\bpoint{Lemma} 
\label{lem:semi-proper-facts}

\begin{enumerate}
\item 
Suppose that $\phi\colon X\longrightarrow V$ is a proper morphism of schemes, with $V$ affine.
Then $X$ is semi-proper.
\item If $X$ is semi-proper and $W\subseteq X$ a closed subscheme, then $W$ is semi-proper.
\end{enumerate}

\bpf
({\em a}) 
Let $X\stackrel{\pi}{\longrightarrow} Z\stackrel{\phi'}{\longrightarrow} V$ be the Stein factorization of $\phi$. 
Then $\pi$ is proper, $\pi_{*}\Osh_{X}=\Osh_{Z}$ and $\phi'$ is an integral morphism. 
Since $V$ is affine and $\phi'$ an affine morphism, $Z$ is affine.  Therefore $Z=\Spec(A)$ with 
$$A = \Gamma(Z,\Osh_{Z})=\Gamma(Z,\pi_{*}\Osh_{X}) = \Gamma(X,\Osh_{X}).$$
Thus the morphism $X\longrightarrow \Spec(\Gamma(X,\Osh_{X}))=Z$ is proper, and so $X$ is semi-proper. 

({\em b}) The composition of the closed immersion $W\hookrightarrow X$ and the natural morphism
$X\longrightarrow \Spec(\Gamma(X,\Osh_{X}))$ is proper, with the target being an affine scheme. Thus
by ({\em a}) $W$ is semi-proper.
\epf

\bpoint{Orbits and open affines on toric varieties} 
\label{sec:orbits-and-affines}
Let $X$ be a toric variety with torus $T_{X}$ and fan $\Sigma_{X}$. 
For each cone $\sigma\in \Sigma_{X}$ we use $O(\sigma)$ for the corresponding $T_{X}$-orbit, and $U_{\sigma}$ for
the union
$$U_{\sigma} := \bigcup_{O(\sigma)\subseteq \overline{O(\sigma')}} O(\sigma').$$
We recall that each $U_{\sigma}$ is a $T_{X}$-stable open affine subset of $X$, and that this procedure gives
a three-way bijection 
$$
\left\{\mbox{Cones $\sigma\in \Sigma_{X}$}\rule{0cm}{0.4cm}\right\} \longleftrightarrow
\left\{\mbox{$T_X$-orbits $O(\sigma)$ on $X$}\rule{0cm}{0.4cm}\right\} \longleftrightarrow
\left\{\mbox{$T_X$-stable open affine subsets of $X$}\rule{0cm}{0.4cm}\right\},
$$
see for example \cite[p.\ 24, theorem 6]{TEI}.
These bijections have the following connections with the natural partial orders on each set.  
For $\sigma$, $\sigma'\in \Sigma_{X}$, 
$$\sigma\subseteq \sigma'  \iff \overline{O(\sigma)}\supseteq O(\sigma') \iff U_{\sigma}\subseteq U_{\sigma'}.$$

\bpoint{Proposition}
\label{prop:torus-stable-covers} 
Let $U\subseteq X$ be a $T_X$-stable open subset, and $\sigma_1$,\ldots, $\sigma_r$ the cones of $\Sigma_{X}$
such that each $O(\sigma_i)$ is closed in $U$.  Then $U_{\sigma_1}$,\ldots, $U_{\sigma_r}$ is the minimal cover
of $U$ by $T_X$-stable open affine subsets. 

\bpf
Since $O(\sigma_i)\subseteq U$, and since open sets are stable under generalization, it follows that each 
$U_{\sigma_i}\subseteq U$.  Fix $i$ and let $U_{\sigma'}\subseteq U$ be a $T_X$-stable open affine containing
$O(\sigma_i)$.  Then $O(\sigma')\subseteq \overline{O(\sigma_i)}$ (where the closure denotes the closure in $X$).  
Since $O(\sigma_i)$ is closed in $U$, 
if $\sigma'\neq \sigma_i$, then $O(\sigma')\not\subseteq U$, 
contradicting the assumption that $U_{\sigma'}\subseteq U$.
Thus, the only $U_{\sigma'}\subseteq U$ containing $O(\sigma_i)$ is $U_{\sigma_i}$, and so 
each $U_{\sigma_i}$ must be part of any $T_X$-stable open affine cover of $U$. 

Finally let us check that $\bigcup_{i=1}^{r} U_{\sigma_i} = U$.   Let $O_{\sigma'}\subseteq U$ be a $T_X$-orbit.
If $O_{\sigma'}$ is closed in $U$, then $\sigma'=\sigma_i$ for some $i$, and $O(\sigma')\subseteq U_{\sigma_i}$. 
If $O_{\sigma'}$ is not closed in $U$, then $\overline{O(\sigma')}\cap U$ contains a $T_X$-orbit $O(\sigma'')$ 
closed in $U$. Then again $\sigma''=\sigma_i$ for some $i$, and $O(\sigma')\subseteq U_{\sigma_i}$.  Thus the 
$U_{\sigma_i}$ cover $U$. \epf

\bpoint{Toric morphisms} 
\label{sec:toric-morphisms}
Let $X$ and $Y$ be toric varieties, with corresponding tori $T_{X}$ and $T_{Y}$,
and inclusions $T_{X}\hookrightarrow X$ and $T_{Y}\hookrightarrow Y$.  Following \cite[\S2, p.\ 138]{deCMM}
a {\em toric morphism} is a morphism $\phi\colon X\longrightarrow Y$ such that $\phi(T_{X})\subseteq T_{Y}$, and 
such that the induced morphism $\phi|_{T_{X}}\colon T_{X}\longrightarrow T_{Y}$ is a morphism of algebraic groups. 
Thus we have a commutative diagram

\vspace{-0.75cm}
\begin{centering}
\begin{tabular}{c}
\begin{pspicture}(0,0)(3,3)
\psset{arrows=->,nodesep=2mm}
\rput(0,2){\rnode{X}{$X$}}
\rput(2,2){\rnode{TX}{$T_{X}$}} 
\rput(0,0){\rnode{Y}{$Y$}}
\rput(2,0){\rnode{TY}{$T_{Y}$}} 
\ncline{X}{Y}
\bput(0.5){$\phi$} 
\ncline{TX}{TY}
\aput(0.5){$\phi|_{T_{X}}$} 
\psset{arrows=H->,hookwidth=1.5mm}
\ncline{TX}{X}
\ncline{TY}{Y}
\end{pspicture}
\end{tabular}\\
\end{centering} 

\vspace{0.5cm}
with $\varphi|_{T_{X}}$ a group homomorphism.   This automatically implies that $\phi$ is $T_{X}$-equivariant 
with respect to the $T_{X}$-action on $Y$ induced by $\varphi|_{T_{X}}$.  For instance, the $T_{X}$-equivariance can 
be expressed as the commutativity of the diagram 

\begin{centering}
\begin{tabular}{c}
\begin{pspicture}(-0.3,-0.3)(4.3,2.5)
\psset{nodesep=2mm,arrows=->} 
\rput(0,2){\rnode{TXX}{$T_{X}\times X$}} 
\rput(0,0){\rnode{TYY}{$T_{Y}\times Y$}} 
\rput(3.5,2){\rnode{X}{$X$}}
\rput(3.5,0){\rnode{Y}{$Y$}}
\ncline{TXX}{TYY}
\bput(0.5){$\phi|_{T_X}\times \phi$}
\ncline{X}{Y}
\aput(0.5){$\phi$} 
\ncline{TXX}{X}
\aput(0.45){$\alpha_X$} 
\ncline{TYY}{Y}
\aput(0.45){$\alpha_Y$}
\end{pspicture},
\end{tabular} \\
\end{centering}

\vspace{0.25cm}
where $\alpha_{X}$ and $\alpha_{Y}$ are the morphisms giving the action of (respectively) $T_X$ on $X$ and $T_Y$ 
on $Y$.  
On the open dense set $T_X\times T_X\subseteq T_X\times X$ the commutativity of the diagram amounts to the fact
that $\varphi|_{T_X}$ is a group homomorphism. 
Since $Y$ is separated, the locus in $T_{X}\times X$ where the morphisms $\phi\circ \alpha_{X}$ and
$\alpha_{Y}\circ (\phi|_{T_X}\times\phi)$ agree is closed, and therefore this locus is equal to all of $T_{X}\times X$. 

\bpoint{Example}
\label{sec:toric-maps-not-rigid-example}
In toric geometry we expect compatibility with the torus action to ``rigidify'' the geometry, making it
essentially combinatorial.  This is not quite the case for a toric morphism $\phi\colon X\longrightarrow Y$ when
$\dim(X)<\dim(Y)$, or more generally when $\dim(\phi(X))< \dim(Y)$. 

For instance, taking $X=\PP^1$, and $Y=\PP^2$ (with torus actions $t\cdot [x:y]=[tx:y]$ and $(s_0,s_1)\cdot 
[z_0:z_1:z_2] =[s_0z_0:s_1z_1:z_2]$ and inclusions $t\mapsto [t:1]$, $(s_0,s_1)\mapsto [s_0:s_1:1]$) 
the morphism $\phi([x:y])=[x:x:y]$ is a toric morphism.

\begin{centering}
\begin{tabular}{ccc}
\begin{tabular}{c}
\begin{pspicture}(-1.5,-2.2)(1.5,2.5)
\psline(-1.0,0)(1.0,0)
\psset{fillstyle=solid,fillcolor=white}
\pscircle(0.7,0){0.045}
\pscircle(-0.7,0){0.045}
\rput(-0.7,0.3){\tiny $[0:1]$}
\rput(0.7,0.3){\tiny $[1:0]$}
\rput(0,-2.0){$\PP^1$}
\end{pspicture}
\end{tabular}
& 
\rule{0.75cm}{0cm}
$\xrightarrow{\rule{0.75cm}{0cm}\phi\rule{0.75cm}{0cm}}$
\rule{0.75cm}{0cm}
&
\begin{tabular}{c}
\begin{pspicture}(-1.5,-2.2)(1.5,2.5)
\psset{fillstyle=solid,fillcolor=vlgray,linecolor=lightgray}
\pscircle(0,0){1.5}
\psset{linestyle=dashed,linecolor=gray, dash=1.5mm 1mm}
\parametricplot{-1.2}{1.2}{t -0.55}
\parametricplot{-1.2}{1.2}{t -0.55 120 \Rot}
\parametricplot{-1.2}{1.2}{t -0.55 -120 \Rot}
\psset{linecolor=black,fillstyle=none,linestyle=solid}
\psline(1.3;210)(1.3;30)
\psset{fillcolor=white,fillstyle=solid}
\pscircle(1.1;90){0.06}
\pscircle(1.1;210){0.06}
\pscircle(1.1;-30){0.06}
\pscircle(0.55;30){0.045}
\rput(2.2;210){\tiny $[0:0:1]$} 
\rput(2.2;-30){\tiny $[0:1:0]$} 
\rput(1.9;90){\tiny $[1:0:0]$} 
\rput(0,-2.0){$\PP^2$}
\end{pspicture}
\end{tabular}\\
\multicolumn{3}{c}{\Fig}
\end{tabular}\\
\end{centering} 

\vspace{0.25cm}
In this example the image of $X$ is not any of the boundary $\PP^1$'s in $Y$, but rather passes through the
largest open torus orbit in $Y$.  Indeed this must happen, since by the group homomorphism condition, the origin
of $T_X$ (in the open orbit $T_{X}\subseteq X$) must be sent to the origin of $T_Y$ 
(in the open orbit $T_{Y}\subseteq Y$). 

\bpoint{Stein factorization of toric morphisms} Let $\phi\colon X\longrightarrow Y$ be a proper toric morphism and
$X\stackrel{\pi}{\longrightarrow} Z\stackrel{\phi'}{\longrightarrow} Y$ the Stein factorization of $\phi$. 

In \cite[Remark 2.3, p.\ 140]{deCMM} the authors show that in this case the Stein factorization may be further 
refined,  introducing a factorization of $\phi'$, 
$Z\longrightarrow W\stackrel{w}{\longrightarrow} Y$
such that $W$ is a toric variety, the induced maps $X\longrightarrow W$ and $W\longrightarrow Y$ are toric morphisms,
and such that $w$ induces a closed immersion\footnote{ The morphism $\phi'$ is finite, and may map with degree $>1$ 
onto its image.  The construction of $W$ absorbs this finite-to-one behaviour, at least generically. }
$T_{W}\hookrightarrow T_{Y}$. 

It is more convenient for us to work directly with $Z$. 
Below we give a proof, similar in spirit to that of \cite[loc.\ cit.]{deCMM}, 
that $Z$ is a toric variety, and that $\pi$ and $\phi'$ are toric morphisms.

\bpoint{Lemma} 
\label{lem:toric-stein-factorization}
Let $\phi\colon X\longrightarrow Y$ be a proper toric morphism, and 
$X\stackrel{\pi}{\longrightarrow} Z\stackrel{\phi'}{\longrightarrow} Y$ the Stein factorization of $\phi$. 
Then $Z$ may be given the structure of a toric variety in such a way that the maps $\pi$ and $\phi'$ are toric
morphisms.  Moreover this structure is unique. 

\bpf
Recall that, in the proof of Stein factorization, $Z=\Spec(\phi_{*}\Osh_{X})$ is constructed as the relative affine 
scheme associated to the sheaf $\phi_{*}\Osh_{X}$ of $\Osh_{Y}$-algebras on $Y$, 
a sheaf supported on the closed set $\phi(X)$.

If $V\subseteq Y$ is an open subset such that $V\cap \phi(X)\neq\emptyset$, then $\phi^{-1}(V)$ is a nonempty open 
subset of the normal irreducible variety $X$, and therefore $(\phi_{*}\Osh_{X})(V) = \Osh_{X}(\phi^{-1}(V))$ is
a normal domain.   It follows that $Z$ is an irreducible normal variety. 

Since $\phi$ is a $T_{X}$-equivariant morphism $\phi_{*}\Osh_{X}$, and so $Z$, inherits a $T_{X}$-action, under which
the morphisms $\pi$ and $\phi'$ are $T_X$-equivariant. 

Let $1_{X}\in X$ denote the image of the identity of $T_{X}$ under the inclusion into $X$, and 
set $1_{Z}:=\pi(1_{X})\in Z$.  Let $H\subseteq T_{Z}$ be the closed subgroup scheme stabilizing $1_{Z}$ under
the $T_{X}$-action.  Since $T_{X}$ is abelian, $H$ is a normal subgroup of $T_{X}$, and we may form the quotient
$T_{Z}:=T_{X}/H$.  The quotient is a torus, since all quotients of tori are again tori\footnote{ For instance,
by \cite[p.\ 114, \S8.5, Proposition]{B} a linear algebraic group is a torus if and only if it is connected and 
diagonalizable.  Connectedness and diagonalizability are inherited by a quotient group.}.

By the orbit-stabilizer theorem, the map $T_{X}\longrightarrow Z$ given by 
$t\mapsto t\cdot 1_{Z}$ ($=\pi(t)$) identifies the orbit of $1_{Z}$ with the quotient variety $T_{X}/H$, and hence with $T_{Z}$.  
It follows from the $T_X$-equivariance of $\pi$ that the induced map $\pi|_{T_X}\colon T_{X}\longrightarrow T_{Z}$
is a group homomorphism. 

By the theorem of Chevalley \cite[Th\'{e}or\`{e}me 1.8.4]{EGA-IV.1}, the $T_X$-orbit $T_Z\subseteq Z$ is open in its closure.  
But $T_Z=\pi(T_X)$, and $T_X$ is dense in $X$.
Since $\pi$ is proper and surjective, we conclude that the closure of $T_{Z}$ is all of $Z$, and 
thus that $T_Z$ is an open dense subset of $Z$.

Consider the map 
$$
\begin{array}{rcccc}
\alpha_{Z} & \colon & H\times Z & \longrightarrow & Z\times Z. \\[0.1cm]
& & (h,z) & \xmapsto{\rule{0.75cm}{0cm}} & (h\cdot z,z) \\
\end{array}
$$
Since $T_{X}$ acts on $T_{Z}\subseteq Z$ through the quotient map, $H$ fixes
all points of $T_{Z}$, and therefore $H\times T_{Z}$ is contained in the closed set $\alpha_{Z}^{-1}(\Delta_{Z})$, 
where $\Delta_{Z}$ denotes the diagonal of $Z\times Z$. Thus, since $T_Z$ is dense in $Z$, 
$\alpha_{Z}^{-1}(\Delta_{Z})=H\times Z$.  I.e., $H$ fixes all points of $Z$, and therefore the $T_X$-action on $Z$
factors through the quotient map to $T_Z$. 

In summary, $Z$ is a normal irreducible variety, with an action of a torus $T_Z$, and compatible open immersion
$T_{Z}\hookrightarrow Z$.  Therefore $Z$ is a toric variety with torus $T_Z$.

The fact that $\phi'|_{T_Z}\colon T_Z\longrightarrow T_{Y}$ is a group homomorphism follows immediately from the 
fact that $\pi|_{T_{X}}\colon T_{X}\longrightarrow T_{Z}$ is a surjective group homomorphism, and that
$\phi|_{T_X}=\phi'|_{T_Z}\circ\pi|_{T_X}$ is a group homomorphism.

Finally, the uniqueness of the toric structure on $Z$ (compatible with $\pi$ and $\phi'$) is clear from the proof, since
the $T_{X}$-action on $Z$ determines $T_Z$, the action of $T_Z$ on $Z$, and the open immersion 
$T_Z\hookrightarrow Z$.  \epf

\point Given a proper toric morphism $\phi\colon X\longrightarrow Y$, and a line bundle $L$ on $X$ our goal 
is to give an algorithm to compute the higher direct images sheaves $R^{i}\phi_{*}L$. 

Letting $X\stackrel{\pi}{\longrightarrow} Z\stackrel{\phi'}{\longrightarrow} Y$ be the Stein factorization of $\phi$
as in Lemma \ref{lem:toric-stein-factorization}, since $\phi'$ is finite we have 
$R^{i}\phi_{*}L = \phi'_{*}\left(R^{i}\pi_{*}L\right)$ for all $i\geq 0$. 
``Compute'', in the sense of this paper, means to produce a finitely generated module over the appropriate
Cox ring which sheafifies to the higher direct image sheaf desired.
The factorization 
$R^{i}\phi_{*}L = \phi'_{*}\left(R^{i}\pi_{*}L\right)$ allows us to split the problem into two pieces, 
first finding a module computing $R^{i}\pi_{*}L$ and then finding a module which computes the result of 
applying $\phi_{*}'$.  

In this paper we concentrate solely on the first problem, and leave the second for future work.
We note that the second problem is not solved by applying the tensor product functor 
$S_{Y}\otimes_{S_{Z}}\underline{\rule{0.4cm}{0cm}}$ , where $S_{Y}$ and $S_{Z}$ denote the Cox rings of $Z$ and $Y$ 
respectively.  For instance, even though $\phi'\colon Z\longrightarrow Y$ is a finite map, this does not imply
that $S_{Z}$ is a finite module over $S_{Y}$.  As a simple example, set $Z=\PP^1\times\PP^1$, $Y=\PP^3$, and
let $\phi'\colon Z\hookrightarrow Y$ be the standard embedding (which, after making compatible choices of toric
structure on $Z$ and $Y$, is a toric morphism).  Then $S_{Y}$ is a polynomial ring in four variables, which
we label $Y_{ij}$ with $i$, $j\in \{1,2\}$, and $S_{Z}$ is also a polynomial ring in four variables, which
we label $Z_1$, $Z_2$, $W_1$, $W_2$.  The homomorphism
$$(\phi')^{*}\colon S_{Y}=\CC[Y_{11},Y_{12},Y_{21},Y_{22}]\longrightarrow \CC[Z_1,Z_2,W_1,W_2]=S_{Z}$$
is given  by $Y_{ij}\mapsto Z_iW_j$ for all $i$, $j\in \{1,2\}$, and thus $S_{Z}$ is not a finite $S_Y$-module.

Returning to the first step, that of computing $R^{i}\pi_{*}L$, after renaming $Z$ and $\pi$ 
we are thus in the case that $\phi\colon X\longrightarrow Y$ is a 
proper {\em toric fibration}, i.e., a proper toric map such that $\phi_{*}\Osh_{X}=\Osh_{Y}$.

\point 
\label{sec:surjectivity-consequence}
The hypothesis that $\phi$ is a toric fibration implies that $\phi$ is surjective with connected fibres.  
The surjectivity, combined with $T_X$-equivariance of $\phi$ means that each $T_X$-orbit in $X$ maps surjectively 
onto a $T_Y$-orbit in $Y$. 
This lining up of torus orbits makes the map essentially combinatorial, in contrast to the kind of behaviour 
in the example of \S\ref{sec:toric-maps-not-rigid-example}.

\point {\bf Proposition} 
\label{prop:lining-up-of-torus-orbits} 
(Combinatorial consequences of the surjectivity hypothesis).
Let $\phi\colon X\longrightarrow Y$ be a proper surjective toric morphism, 
$\tau\in \Sigma_{Y}$ a cone, and $V_{\tau}\subseteq Y$ the corresponding $T_Y$-stable affine 
open\footnote{ We use $V_{\tau}$ instead of $U_{\tau}$ to keep track of the fact that this is an open subset of $Y$ 
and not $X$} in $Y$.  We consider cones $\sigma\in \Sigma_{X}$ with $O(\sigma)\subseteq \phi^{-1}(V_{\tau})$. 

\begin{enumerate}
\item If $O(\sigma)$ is closed in $\phi^{-1}(V_{\tau})$ then $\phi(O(\sigma))=O(\tau)$, the morphism
\PauseEnumerate

\begin{equation}
\label{eqn:O-sigma-to-O-tau-map}
\phi|_{O(\sigma)}\colon O(\sigma)\longrightarrow O(\tau)
\end{equation}

\parshape 1 1.3cm 12cm
is a finite surjective map, and $\dim(O(\sigma))=\dim(O(\tau))$.

\ResumeEnumerate
\item If $O(\sigma)\subseteq \phi^{-1}(V_{\tau})$ then $\dim(O(\sigma))\geq \dim(O(\tau))$.

\medskip
\item If $O(\sigma)\subseteq \phi^{-1}(V_{\tau})$ then $O(\sigma)$ is closed in 
$\phi^{-1}(V_{\tau})$ if and only if $\dim(O(\sigma))=\dim(O(\tau))$.
\end{enumerate}

\parshape 1 0cm \textwidth

\bpf
Set $U:=\phi^{-1}(V_{\tau})$.  We first prove ({\em a}). 
If $O(\sigma)$ is closed in $U$, then since $\phi|_{U}\colon U\longrightarrow V_{\tau}$
is proper, $\phi(O(\sigma))$ is closed in $V_{\tau}$.  By the remark in \S\ref{sec:surjectivity-consequence},
$\phi(O(\sigma))$ is a $T_{Y}$-orbit.  Since $O(\tau)$ is the unique closed $T_{Y}$-orbit in $V_{\tau}$, we conclude 
that $\phi(O(\sigma))=O(\tau)$, i.e., \eqref{eqn:O-sigma-to-O-tau-map} is a surjective map.

Both $O(\sigma)$ and $O(\tau)$ are affine varieties (specifically 
$O(\sigma)\cong (\CC^{*})^{k}$ and $O(\tau)\cong (\CC^{*})^{\ell}$ with $k=\codim(\sigma,N_{\RR,X})$ and
$\ell=\codim(\tau,N_{\RR,Y})$), so $\phi|_{O(\sigma)}$ is a map between affine varieties and hence an affine map. 
But $\phi|_{O(\sigma)}$ is also the composition of the closed inclusion $O(\sigma)\hookrightarrow U$ and
the proper map $\phi|_{U}\colon U\longrightarrow V_{\tau}$ and so proper.

Since a proper affine map between finite-type schemes is finite, \eqref{eqn:O-sigma-to-O-tau-map} is a finite map, 
and we have already seen that it is surjective. 
Thus \eqref{eqn:O-sigma-to-O-tau-map} is a finite surjective map, and so $\dim(O(\sigma))=\dim(O(\tau))$. 

For ({\em b}), if $O(\sigma)$ is not closed in $U$ then the closure of $O(\sigma)$ in $U$ contains a $T_X$-orbit
of strictly smaller dimension which {\em is} closed in $U$.  By part (a), this smaller orbit has dimension equal to
$\dim(O(\tau))$.  Thus, in this case, $\dim(O(\sigma))> \dim(O(\tau))$. 

On the other hand, if $O(\sigma)$ is closed in $U$ then $\dim(O(\sigma))=\dim(O(\tau))$ by ({\em a}).  In both cases
we have $\dim(O(\sigma))\geq \dim(O(\tau))$. 

For ({\em c}), the direction $O(\sigma)$ is closed in $U$ $\implies$ $\dim(O(\sigma))=\dim(O(\tau))$ 
is part of ({\em a}).  Now assume that $\dim(O(\sigma))=\dim(O(\tau))$.  If $O(\sigma)$ were not closed in $U$ then,
as in the first part of the argument for ({\em b}) we would conclude that $\dim(O(\sigma))> \dim(O(\tau))$, 
in contradiction to our assumption. \epf

\bpoint{Functoriality of the $\geq 0$ construction and cell structure}
\label{sec:funct-of-non-neg-construction}
As noted in \S\ref{sec:functoriality-initial-discussion} a toric morphism $\phi\colon X\longrightarrow Y$
induces a map $\phi_{\geq 0}\colon X_{\geq 0}\longrightarrow Y_{\geq 0}$ of topological spaces, and this association
is functorial. 
Giving $X_{\geq 0}$ and $Y_{\geq 0}$ the structure of a cell complex as in \S\ref{sec:cell-complex-structure},
the map $\phi_{\geq 0}$ is, in general, not a map of cell complexes.

Recall that, by definition, a continuous map $f\colon Z\longrightarrow W$ 
between two topological spaces with the structure of cell complexes is a 
cellular map if and only if each $k$-cell of $Z$ is mapped into the union of the cells of dimension $\leq k$ of $W$.

Thus, the map $\phi_{\geq 0}$ corresponding to the example in \S\ref{sec:toric-maps-not-rigid-example} is not
a map of cell complexes, since the $1$-cell of $X_{\geq 0}$ (corresponding to the open torus orbit in $\PP^1$) 
is taken to the $2$-cell of $Y_{\geq 0}$ (corresponding to the open torus orbit in $\PP^2$).

On the other hand, if $\phi$ is proper and surjective, then $\phi_{\geq 0}$ is a morphism of cell complexes.
The cells of $X_{\geq 0}$ and $Y_{\geq 0}$ are in one-to-one correspondence with
the torus orbits on $X$ and $Y$ respectively, a correspondence compatible with inclusions, closures, and 
dimensions (the real dimension of a cell is equal to the complex dimension of the corresponding orbit), 
and also compatible with the maps $\phi_{\geq 0}$ and $\phi$.  
When $\phi$ is proper and surjective the fact that a $k$-cell of
$X_{\geq 0}$ maps into the cells of dimension $\leq k$ of $Y_{\geq 0}$ is therefore a consequence of
Proposition \ref{prop:lining-up-of-torus-orbits}(b).

\point 
\label{sec:discussion-of-W-non-neg}
In the proof of Theorem \ref{thm:construction-of-toric-reduced-cover} below we will need to construct
a particular cell complex, and know that it has the cohomology of a point. 

To set this up, suppose that $X$ is a semi-proper toric variety, set $Z=\Spec(\Gamma(X,\Osh_{X}))$ and let 
$\phi\colon X\longrightarrow Z$ be the natural morphism.  By
Lemma \ref{lem:toric-stein-factorization}, $Z$ has the structure of a toric variety, and with that structure
$\phi$ is a toric morphism. 

Let $\Sigma_{Z}$ be the fan of $Z$. Since $Z$ is affine, $\Sigma_{Z}$ has a unique maximal cone $\tau$,
corresponding to the unique closed torus-stable orbit in $Z$.   The key case is when this orbit $O(\tau)$
is zero dimensional, i.e., when $z_0:=O(\tau)$ is a closed point of $Z$. 

In this case set $W=\phi^{-1}(z_0)$, with its reduced subscheme structure. Since $\phi$ is $T_{X}$-equivariant,  
and $z_0$ fixed, $W$ is stable under the action of $T_{X}$, and thus a union of $T_{X}$-orbits.  
Moreover, $W$ is proper, and thus each component of $W$ is a proper toric variety.

Since $W$ is given locally in $X$ by the vanishing of monomial conditions, $W_{\geq 0}$ makes sense (it
is the subset of $X_{\geq 0}$ given by the same monomial conditions) and is the union
of the $(W_i)_{\geq 0}$ where $W_i$, $i=1$,\ldots, $s$ are the components of $W$.  Moreover, $W_{\geq 0}$ has the 
structure of a cell complex, given by the cells of each $(W_i)_{\geq 0}$ 
(as in \S\ref{sec:cell-complex-structure}).  The closure of a cell in $(W_i)_{\geq 0}$ is a union of cells 
in $(W_i)_{\geq 0}$, and so contained in $W_{\geq 0}$.  
Alternatively, since $\phi$ is a surjective map, it induces a map of cell complexes
$X_{\geq 0}\longrightarrow Z_{\geq 0}$, and $W_{\geq 0}$ is the sub cellular complex of $X_{\geq 0}$ which is the
fibre over the unique $0$-cell of $Z_{\geq 0}$ (the cell corresponding to $z_0$).

\point
\label{prop:cohomology-of-W-pos}
{\bf Proposition} (Cohomology type of $W_{\geq 0}$).
As in \S\ref{sec:discussion-of-W-non-neg} let $X$ be a semi-proper toric variety such that the unique
closed torus orbit of $Z:=\Spec(\Gamma(X,\Osh_{X}))$ is a point $z_0$, and set $W:=\phi^{-1}(z_0)$,
where $\phi\colon X\longrightarrow Z$ is the natural map.  Then 
$H^{0}(W_{\geq 0},\ZZ)=\ZZ$ and $H^{i}(W_{\geq 0},\ZZ)=0$ for all $i\geq 1$. 

\bpf
If $Z$ is itself a point, $Z=z_0$, then $W=X$, and therefore by Theorem \ref{thm:top-structure-of-X}(b)
$W_{\geq 0}=X_{\geq 0}$ is homeomorphic to a ball.  
We may therefore suppose that $Z$ is not a point.   

If $W$ is irreducible, then since $W$ is a proper toric variety, 
Theorem \ref{thm:top-structure-of-X}(b) again gives that $W_{\geq 0}$, is homeomorphic to a ball. 
Otherwise, $W$ is reducible, with irreducible components $W_i$, $i=1$, \ldots, $s$, and $W_{\geq 0}$ is the union
of the closed subsets $(W_i)_{\geq 0}$, $i=1$,\ldots, $s$. 

Recall that a cell complex is called regular if the inclusion of each (relatively open) 
cell of dimension $k$ extends to a homeomorphism of the $k$-dimensional ball onto the closure of that cell.  
By Theorem \ref{thm:top-structure-of-X}(a) each $(W_i)_{\geq 0}$ is a regular cell complex.  Since each
$(W_i)_{\geq 0}$ is closed in $W_{\geq 0}$, it follows that $W_{\geq 0}$ is also a regular cell complex. 
The cohomology of a regular cell complex is completely combinatorial, determined solely by the incidence 
relations among the cells (see for example \cite[Proposition 4.7.8 and Corollary 4.7.9]{OM}).  To compute the
cohomology of $W_{\geq 0}$ we may therefore use any combinatorial model keeping track of the cells and their 
incidence relations.  We will use the (order and dimension reversing) correspondence between cells of $W_{\geq 0}$
and cones in $\Sigma_{X}$, parallel to the correspondence between $T_{X}$-orbits on $X$ and cones
in $\Sigma_{X}$.  

Let $\Sigma_{Z}\subseteq N_{Z,\RR}$ be the fan of $Z$, and $\tau\in \Sigma_{Z}$ the cone corresponding to $z_0$.  
Since $\dim(Z)>0$, $\dim(N_{Z,\RR})>0$, and since the unique closed orbit of $Z$ is a point, 
$|\tau|$ is a pointed convex subset\footnote{ In this subsection, in contrast to the rest of the paper, we use
$|\tau|$ to denote the support of a cone, rather than its dimension.} of $N_{Z,\RR}$. 
Since the map $\phi\colon X\longrightarrow Z$ is surjective, the induced map on tori is surjective, and hence
the induced map $\phi_{N}\colon N_{X,\RR}\longrightarrow N_{Z,\RR}$ is also surjective.   Let $|\tau|^{\circ}$
denote the interior of the support of $\tau$.   The cones corresponding to the cells in $W_{\geq 0}$ (equivalently
to the $T_{X}$-orbits in $W$) are those whose relative interiors lie in $\phi^{-1}_{N}(|\tau|^{\circ})$ i.e., 
those cones of $\Sigma_{X}$ whose image under $\phi_{N}$ is not contained in the boundary of $|\tau|$. 
We also note that since $Z$ is affine, and $\tau$ the maximal cone of $\Sigma_{Z}$, all cones in $\Sigma_{X}$
have image in $|\tau|$. 

Set $n=\dim(X)$, let $S^{n-1}\subset N_{X,\RR}$ denote the unit sphere,
and set $U:=S^{n-1}\cap \phi^{-1}_{N}(|\tau|^{\circ})$.  The cones of $\sigma_{X}$ give a 
decomposition of $\overline{U}$, and it is the complex built out of those cones which intersect $U$ which
we are interested in.  This complex is naturally a quotient complex $\mathcal{Q}_{\smallbullet}$ of a complex
computing the homology of $S^{n-1}$, a quotient complex which we now construct.

The pullback of $|\tau|^{\circ}$ is an open convex subset of $N_{X,\RR}$, cut out by the pullbacks of the strict
inequalities defining $|\tau|^{\circ}$ in $N_{Z,\RR}$.   It follows that $U$ is contractible to a point within
$S^{n-1}$.  Let $K$ be the closed complement of $U$ in $S^{n-1}$. The previous remark shows that $K$ has the homotopy
type of $S^{n-1}$ minus a point, and thus that $K$ has the homology of a point. 
Choose a decomposition of $K$ into cells compatible on the boundary of $U$ with the decomposition of $\overline{U}$ 
given by $\Sigma_{X}$.   We thus obtain a decomposition of $S^{n-1}$ into cells, compatible with the 
decomposition of $K$ and $U$, and a short exact sequence of homology complexes~:

\begin{equation}
\label{eqn:homology-complex}
0 \longrightarrow \mathcal{C}_{\smallbullet}(K,\ZZ) \longrightarrow \mathcal{C}_{\smallbullet}(S^{n-1},\ZZ)
\longrightarrow \mathcal{Q}_{\smallbullet}\longrightarrow 0.
\end{equation}

where $\mathcal{C}_{\smallbullet}(K,\ZZ)$ and $\mathcal{C}_{\smallbullet}(S^{n-1}\ZZ)$ are the homology complexes 
coming from the respective decompositions of $K$ and $S^{n-1}$, and where $\mathcal{Q}_{\smallbullet}$ denotes 
the quotient complex. 

In the correspondence between cells of $W_{\geq 0}$ and cones of $\Sigma_{X}$, a $k$-dimensional cell of $W_{\geq 0}$
corresponds to an $(n-k)$-dimensional cone of $\Sigma_{X}$, and hence induces an $((n-1)-k)$-dimensional cell
in $U$.  The correspondence between orbits and cones is order reversing.  It follows that the complex which computes
the cohomology of $W_{\geq 0}$ is the same as the homology complex $\mathcal{Q}_{\smallbullet}$ above, and that
$H^{i}(W_{\geq 0},\ZZ) = H_{n-1-i}(\mathcal{Q}_{\smallbullet})$ for all $i$.  
The statement we want to prove about the cohomology 
of $W_{\geq 0}$ becomes the statement that $H_{n-1}(\mathcal{Q}_{\smallbullet})=\ZZ$, 
and that $H_{j}(\mathcal{Q}_{\smallbullet})=0$ for all $j\neq n-1$.  
Since $K$ has the homology of a point, this follows immediately from the long exact sequence in homology associated to 
\eqref{eqn:homology-complex}. \epf

\point
\label{lem:splitting-lemma}
{\bf Lemma} (Product lemma).
Let $X$ be a semi-proper toric variety and $Z=\Spec(\Gamma(X,\Osh_X))$. 
Suppose that the unique closed torus orbit in $Z$ is not a point, and thus is of the form $(\CC^{*})^{k}$ for
some $k\geq 1$.
Then there is a semi-proper toric variety $X_1$ and map $\phi_1\colon X_1\longrightarrow Z_1$, with 
$Z_1=\Spec(\Gamma(X_1,\Osh_{X_1}))$ such that the 
natural map $\phi\colon X\longrightarrow Z$ is the product of $\phi_1$ with $(\CC^{*})^{k}$, and such that the
unique closed torus orbit of $Z_1$ is a point. 

\bpf
Let $\Sigma_{X}\subset N_{X,\RR}$ and $\Sigma_{Z}\subset N_{Z,\RR}$ be the fans of $X$ and $Z$ respectively.
Since the unique closed orbit of $Z$ is isomorphic to $(\CC^{*})^{k}$, $\Sigma_{Z}$ spans a subspace of
codimension $k$ in $N_{Z,\RR}$.  As in the proof of Proposition \ref{prop:cohomology-of-W-pos}, the fact that 
$\phi$ is surjective implies that $\phi_{N}\colon N_{X,\RR}\longrightarrow N_{Z,\RR}$ is surjective.
Thus, $\phi^{-1}_{N}(\Sigma_{Z})$ spans a subspace of codimension $k$ in $N_{X,\RR}$, one which contains
the fan $\Sigma_{X}$, since $X$ maps to $Z$. 

Let $X_1$ be the toric variety defined by the fan $\Sigma_{X}$ thought of as a subspace of the span of 
$\phi^{-1}_{N}(\Sigma_{Z})$ (a rationally defined subspace of $N_{X,\RR}$), $Z_1$ the toric variety defined 
by the fan $\Sigma_{Z}$ thought of as a subspace of the span of $\Sigma_{Z}$, and 
$\phi_1\colon X_1\longrightarrow Z_1$ the map induced by restricting $\phi_{N}$ to those subspaces. 

Then : ({\em i}) $\phi_1$ is proper, because this may be checked on the map of fans, and the map of fans for 
$\phi_1$ is the same as the map of fans for $\phi$. ({\em ii}) $Z_1$ is affine since $\Sigma_{Z}$ has a unique
maximal cone. ({\em iii}) The unique closed torus orbit in $Z_1$ is a point since $\Sigma_{Z_1}$ (=$\Sigma_{Z}$) 
spans $N_{Z_1,\RR}$. ({\em iv}) $\phi$ is the product of $\phi_1$ with $(\CC^{*})^{k}$ (clear from the construction). 
Finally ({\em v}) $Z_1=\Spec(\Gamma(X_1,\Osh_{X_1}))$, since we have a natural map
$\Gamma(Z_1,\Osh_{Z_1})\longrightarrow \Gamma(X_1,\Osh_{X_1})$ which (by ({\em iv})) after tensoring over $\CC$
with $\CC[z^{\pm}_{1},\ldots, z^{\pm}_{k}]$ becomes an isomorphism, and therefore must have been an isomorphism
to begin with, and furthermore, $Z_1$ is affine and hence 
$Z_1=\Spec(\Gamma(Z_1,\Osh_{Z_1}))=\Spec(\Gamma(X_1,\Osh_{X_1}))$. \epf

\point
\label{thm:construction-of-toric-reduced-cover}
{\bf Theorem} (Construction of a reduced \v{C}ech complex for a semi-proper toric variety).
Let $X$ be a semi-proper toric variety.  Then there exists a reduced \v{C}ech complex $(P,\psi)$ such that
each $U_{\sigma}$ is a torus-stable affine and such that $\dim(P)=\cdim(X)$. 

\bpf
Let $Z=\Spec(\Gamma(X,\Osh_{X}))$ and let $\phi\colon X\longrightarrow Z$ be the natural morphism.  By
Lemma \ref{lem:toric-stein-factorization} $Z$ has the structure of a toric variety, and with that structure
$\phi$ is a toric morphism. 

Let $\Sigma_{Z}$ be the fan of $Z$. Since $Z$ is affine, $\Sigma_{Z}$ has a unique maximal cone $\tau$,
corresponding to the unique closed torus-stable orbit in $Z$. We first consider the case
that $O(\tau)$ is zero dimensional, i.e., when $z_0:=O(\tau)$ is a point of $Z$.

In this case let $W=\phi^{-1}(z_0)$, with its reduced subscheme structure, and set $P:=W_{\geq 0}$. 
Using the correspondence between closed cells of $W_{\geq 0}$ and the closures of $T_X$-orbits on $W$, we label
the closed cells of $P$ by the corresponding $\sigma\in\Sigma_{X}$.   For each $\sigma\in P$, we let $U_{\sigma}$ be
the $T_X$-stable open affine subset of $X$ corresponding to $\sigma$ under the standard correspondence between 
cones in $\Sigma_{X}$ and torus stable open affines (as in \S\ref{sec:orbits-and-affines}).
Since the correspondence between closures of orbits and open affines is order reversing,
this correspondence satisfies (P1).  The vertices of $P$ correspond to the zero-dimensional $T_X$-orbits in $W$, and
by Proposition \ref{prop:lining-up-of-torus-orbits}(c) these are the only closed $T_X$-orbits in $X$.
By Proposition \ref{prop:torus-stable-covers} the open sets $U_{v}$, for $v\in P$ a vertex, cover $X$.  
Thus this cover satisfies (P2). 

We now check that this cover satisfies (P3) and (P4).   Let $x$ be a point of $X$.   First suppose that
$x$ is in the largest torus orbit, i.e., that $x\in T_{X}$.  By definition $P_{x}$ consists of the cells $\sigma\in P$
such that $x\in U_{\sigma}$.  Since each $U_{\sigma}$ contains $T_{X}$, when $x\in T_{X}$, $P_{x}=P=W_{\geq 0}$.
By Proposition \ref{prop:cohomology-of-W-pos}, $W_{\geq 0}$ has the cohomology of a point.

If $x\not\in T_{X}$, let $\gamma\in \Sigma_{X}$ be the cone so that $x\in O(\gamma)$,
and set $X'=\overline{O}({\gamma})$.  
Since $X'$ is a closed subset of a semi-proper variety, $X'$ is also semi-proper 
(Lemma \ref{lem:semi-proper-facts}(b)), and of course $X'$ is a toric variety.  
Moreover, setting $Z'=\Spec(\Gamma(X',\Osh_{X'}))$, 
Lemma \ref{lem:toric-stein-factorization} and the proof of 
Lemma \ref{lem:semi-proper-facts}(a) shows that we have a commutative square of $T_X$-equivariant morphisms

\begin{centering}
\label{eqn:commutative-stein}
\begin{tabular}{c}
\begin{pspicture}(-0.3,-0.3)(4.3,2.5)
\psset{nodesep=2mm}
\rput(0,2){\rnode{Xp}{$X'$}} 
\rput(0,0){\rnode{Zp}{$Z'$}} 
\rput(3.5,2){\rnode{X}{$X$}}
\rput(3.5,0){\rnode{Z}{$Z$}}
\ncline{Xp}{Zp}
\bput(0.5){$\phi'$}
\ncline{X}{Z}
\aput(0.5){$\phi$} 
\psset{arrows=H->,hookwidth=-1.5mm} 
\ncline{Xp}{X}
\psset{arrows=->}
\ncline{Zp}{Z}
\aput(0.45){$j$}
\end{pspicture} 
\end{tabular},\\
\end{centering}

with the map $j\colon Z'\longrightarrow Z$ finite.   It is not clear to the author 
whether $j$ is a closed immersion, but since $j$ is finite and $T_X$-equivariant, the following argument
shows that $j^{-1}(z_0)$ consists of a single point.   First, since $j$ is finite $j^{-1}(z_0)$ consists of a finite 
number of points, and since $z_0$ is $T_X$-stable, the same is true of the subset $j^{-1}(z_0)$.  But
$T_{X}$ is connected and so must fix each point of $j^{-1}(z_0)$.  Since  $Z'$ is an affine toric variety it
has a unique closed $T_{X}$-orbit, and thus $j^{-1}(z_0)$ consists of a single point, which we name $z'_0$.

Set $W':=(\phi')^{-1}(z'_0)$.  By the commutativity of the diagram and the definition of $W$ we have the equality 
$W'= W\cap X'$ of closed subsets of $X$.  Thus, a $T_X$-orbit is in $W'$ if and only if it is in $W$ and in $X'$.  
The same is therefore true of the corresponding cells of the cell decomposition,  and so 
$$W'_{\geq 0} = W_{\geq 0} \cap X'_{\geq 0}.$$

For any $\sigma\in \Sigma_{X}$, $x\in U_{\sigma}$ if and only if $O(\gamma)\subset U_{\sigma}$, if and only if
$O(\sigma)\subseteq \overline{O}(\gamma)=X'$.    By definition of $P_{x}$, a cell $\sigma$ in $P=W_{\geq 0}$ 
is contained in $P_{x}$ if and only if $x\in U_{\sigma}$. 

Thus, $P_{x}$ consists of those cells $\sigma\in W_{\geq 0}$ such that
$O(\sigma)\subseteq X'$ (i.e., so that the cell $\sigma$ is contained in $X_{\geq 0}$), 
and therefore 
$$P_{x}=W_{\geq 0}\cap X'_{\geq 0}=W'_{\geq 0}.$$  
Since $X'$ is semi-proper, and the unique closed torus orbit of $Z'$ is zero-dimensional, we may
apply Proposition \ref{prop:cohomology-of-W-pos} to conclude that $P_{x}=W'_{\geq 0}$ has the cohomology
of a point.

This finishes the argument that the cover satisfies (P1)--(P4) in the case when the closed orbit
of $Z$ is zero-dimensional.  Continuing in that case, we next show that $\dim(P)=\cdim(X)$.

Let $d=\dim(W)$, and let $W_k$ be a component of $W$ of dimension $d$.  Since the real dimensions of the cells of $P$
are equal to the complex dimension of the corresponding $T_X$-orbits, $\dim(P)=d$.  
By \cite[Main Theorem]{K} there is a coherent sheaf $\Fsh$ on $W_k$ such that $H^{d}(W_k,\Fsh)\neq 0$.  Pushing forward
$\Fsh$ by the closed immersion $W_k\hookrightarrow X$, we conclude that $\cdim(X)\geq \dim(P)$.  

Since the cover constructed above satisfies (P1)--(P4), and each element of the cover is affine, by 
Theorem \ref{thm:reduced-cech-computes-cohomology} for any quasi-coherent sheaf $\Fsh$ on $X$  the reduced 
\v{C}ech complex $\CComplex_{P}(\Fsh)$ computes the cohomology of $X$. Therefore $\cdim(X)\leq \dim(P)$ and
thus $\dim(P)=\cdim(X)$.   This proves Theorem \ref{thm:construction-of-toric-reduced-cover}
in the case that the unique closed orbit in $Z$ is zero-dimensional.

In the general case, if the closed orbit in $Z$ has dimension $k$, we use the construction above and take the product
with $(\CC^{*})^k$. Specifically, we apply Lemma \ref{lem:splitting-lemma} to obtain a semi-proper toric variety $X_1$,
with map $\phi_1\colon X_1\longrightarrow Z_1:=\Spec(\Gamma(X_1,\Osh_{X_1}))$ such that the closed orbit of $Z_1$ is
zero-dimensional, and such that $\phi$ is obtained from $\phi_1$ by taking the product with $(\CC^{*})^k$ (i.e., 
by base-changing with the map from $(\CC^{*})^k$ to a point), and such that $T_{X} = T_{X_1}\times (\CC^{*})^k$.

Set $P:=(X_1)_{\geq 0}$.  Each $T_X$-orbit in $X$ is the product with $(\CC^{*})^k$ of a $T_{X_1}$-orbit
in $X_1$.  Thus, if we take the affine cover $\{V_{\sigma}\}_{\sigma \in P}$ of $X_1$ produced by the argument 
in the case above, and take the product with
$(\CC^{*})^k$, we obtain an affine cover $\{U_{\sigma}\}_{\sigma\in P}$ of $X$ by the $T_X$-stable affines 
$U_{\sigma}=V_{\sigma}\times (\CC^{*})^k$.  
Each point $x$ of $X$ may be written as a pair $(x_1,t_2)$, with $x_1\in X_1$ and $t_2\in (\CC^{*})^k$. 
Then, the subcomplex $P_{x}$ with respect to the cover $\{U_{\sigma}\}_{\sigma\in P}$ of $X$ is 
$P_{x_1}$ with respect to the cover $\{V_{\sigma}\}_{\sigma\in P}$ of $X_1$.
Thus, the fact that the cover of $X$ by the $U_{\sigma}$ satisfies (P3) and (P4) follows from the fact that the cover
of $X_1$ by the $V_{\sigma}$ satisfies those conditions.  The conditions (P1) and (P2) are also clear. 

Finally, to check that $\cdim(X)=\dim(P)$, since the cover satisfies (P1)--(P4) the reduced
\v{C}ech complexes $\CComplex_{P}(\Fsh)$ compute the cohomology of any quasi-coherent sheaf $\Fsh$, and 
we again get the inequality $\cdim(X)\leq \dim(P)$.  The inclusion of $X_1$ into $X$ as a closed subscheme
(for instance, as $X_1\times \{1\}\hookrightarrow X_1\times (\CC^{*})^k=X$, with $1$ denoting the unit of 
$(\CC^{*})^{k}$) gives the opposite inequality $\dim(P)=\cdim(X_1)\leq \cdim(X)$.

This finishes the proof of Theorem \ref{thm:construction-of-toric-reduced-cover}. \epf

\bpoint{Example}
\label{ex:semi-proper-example} 
Here is an illustration of Theorem \ref{thm:construction-of-toric-reduced-cover}. 
Let us start with the setup of Example \ref{sec:X-to-F1-example}, a map $\phi\colon X\longrightarrow Y$ with $Y$ 
the first Hirzebruch surface, and $B\subset Y$ the curve of self-intersection $(-1)$. Let $V\subset Y$ be a 
torus-stable open affine corresponding to one of the $T_Y$-fixed points on $B$, and let us apply 
the theorem to the semi-proper variety $\phi^{-1}(V)$. 

\vspace{-0.2cm}
\hfill
\begin{tabular}{c}
\psset{unit=0.5cm} 
\begin{pspicture}(-5,-7)(4,4)
\psset{fillstyle=solid,fillcolor=white,opacity=0.85,linecolor=gray}
\rput(!2 -3 2 \TT){\rnode{a}{}}
\rput(!2 -3 -2 \TT){\rnode{b}{}}
\rput(!2 2 -2 \TT){\rnode{c}{}}
\rput(!2 2 1 \TT){\rnode{d}{}}
\rput(!2 0 2 \TT){\rnode{e}{}}
\rput(!-2 -3 2 \TT){\rnode{f}{}}
\rput(!-2 -3 -2 \TT){\rnode{g}{}}
\rput(!-2 2 -2 \TT){\rnode{h}{}}
\rput(!-2 2 1 \TT){\rnode{i}{}}
\rput(!-2 0 2 \TT){\rnode{j}{}}
\pspolygon(f)(g)(h)(i)(j)
\psline(b)(g)
\pspolygon(a)(b)(c)(d)(e)
\pspolygon(a)(e)(j)(f)
\pspolygon(d)(c)(h)(i)
\rput(!2 -3 -6 \TT){\rnode{p}{}}
\rput(!2 2 -6 \TT){\rnode{q}{}}
\rput(!-2 2 -6 \TT){\rnode{r}{}}
\rput(!-2 -3 -6 \TT){\rnode{s}{}}
\pspolygon(p)(q)(r)(s)
\psline[arrows=->](!0 -0.5 -3.5 \TT)(!0 -0.5 -5.5 \TT)
\psset{linecolor=black,linewidth=0.4mm}
\psline(e)(d)(c)
\psset{fillstyle=solid,fillcolor=black} 
\pscircle(e){0.14}
\pscircle(d){0.14}
\pscircle(c){0.14}
\pscircle(q){0.14}
\rput(q){\rule{0.80cm}{0cm}$v_0$}  
\end{pspicture} \\
\Fig\label{fig:semi-proper} \\
\end{tabular}

\vspace{-6.1cm}
\parshape 1 0cm 11.5cm
The open set $V$ is isomorphic to $\AA^2$, and $B\cap V$ is one of the axes.  Following the construction in
Example \ref{sec:X-to-F1-example}, $\phi^{-1}(V)$ is then the blowup of $V\times \PP^1$ along a torus-fixed
curve lying over $V\cap B$. 

\parshape 1 0cm 11.5cm
The affine open set $V$ has a unique torus fixed point, namely the origin $z_0$ of $V$, which corresponds to a 
vertex $v_0$ of $Y_{\geq 0}$.  Following the construction in Theorem \ref{thm:construction-of-toric-reduced-cover}, 
we set $W=\phi^{-1}(z_0)$, and use $P=W_{\geq 0}$ for constructing the reduced \v{C}ech cover. 

\parshape 1 0cm 11.5cm
Because passing to the non-negative points is functorial, we may also realize $P$ as $P=\phi_{\geq 0}^{-1}(v)$.
Figure \ref{fig:semi-proper} shows the fibre of $\phi_{\geq 0}$ over $v_0$ (as well as the rest of $X_{\geq 0}$
mapping to $Y_{\geq 0}$ as in Figure \ref{fig:X-to-Y}). 

\parshape 2 0cm 11.5cm 0cm \textwidth
As the figure shows, $P$ is the cell complex consisting of three vertices and two edges.  
We may also compute that $W=\phi^{-1}(z_0)$ is the union of two $\PP^1$'s, meeting in a node, with the node
an end point of each of the intervals $(\PP^1)_{\geq 0}$.   As guaranteed by 
Proposition \ref{prop:cohomology-of-W-pos}, $P$ has the cohomology of a point.

The proof of Theorem \ref{thm:construction-of-toric-reduced-cover} further guarantees that for each $x\in X$ the
subcomplex $P_{x}\subseteq P$ also has the cohomology of a point (for instance, it is impossible that $P_x$
be one of the edges of $P$ and an isolated vertex). 

Given a coherent sheaf on $\phi^{-1}(V)$, the resulting reduced \v{C}ech complex has three terms in 
degree $0$, and two terms in degree $1$.  The variety $\phi^{-1}(V)$ has cohomological dimension $1$, 
the dimension of $P$.

\bpoint{Remark} 
Even if $X$ is a semi-proper toric variety defined over a field other than $\CC$ (for instance, a field of 
characteristic $p$) we may still use the method above to construct a reduced \v{C}ech complex on $X$ 
satisfying the conditions of Theorem \ref{thm:construction-of-toric-reduced-cover}. 

Letting $X_{\CC}$ be the corresponding toric variety over $\CC$
(i.e., using the same fan), we construct a reduced \v{C}ech complex on $X$ by using the cell complex
$P:=(X_{\CC})_{\geq 0}$, and for $\sigma\in P$, the corresponding open affine $U_{\sigma}$ of $X$. 
The properties (P1)--(P4) only depend on the combinatorics of the torus orbits and their closures, and
how the cones in $\Sigma_{X}$ map to the cones in $\Sigma_{Z}$.  This information is the same for $X$ as it is
for $X_{\CC}$.  The arguments to show that $\dim(P)=\cdim(X)$ also carry over since (in the case that the
unique closed orbit of $Z$ is not a point) Lemma \ref{lem:splitting-lemma} is valid, with the same proof
(and with $\CC^{*}$ replaced by $k^{*}$), over fields other than $\CC$. 

%


\section{Construction of the complex and related results}
\label{sec:construction-of-D-E-complex}

Throughout this section $\phi\colon X\longrightarrow Y$ denotes a {\em toric fibration}, i.e., a
proper toric map between toric varieties, which in addition satisfies $\phi_{*}\Osh_{X}=\Osh_{Y}$.   

\bpoint{Lemma} 
\label{lem:kernel-of-toric-fibration}
Let $\phi\colon X\longrightarrow Y$ be a toric fibration.  Then the kernel of the induced
homomorphism of tori $\phi|_{T_{X}}\colon T_{X}\longrightarrow T_{Y}$ is connected, and thus the kernel is
itself a torus.

\bpf Base changing, we may assume that $Y=T_{Y}$.  By \cite[Th\'{e}or\`{e}m (6.9.1)]{EGA-IV.2}
there is a nonempty open set $U\subseteq Y$ over which $\phi$ is flat.  By equivariance of the morphism, and
the fact that $Y=T_{Y}$, we conclude that $\phi$ is flat over all of $Y$.

Let $\eta$ be the generic point of $Y$.  Then the fibre $X_{\eta}$, being the localization of the normal
variety $X$, is also normal.  Since we are in characteristic zero, the corresponding geometric fibre
$X_{\overline{\eta}}$ (i.e,. base changing to an algebraic closure of $\kappa(\eta)$) remains normal. 
By \cite[Th\'{e}or\`{e}m (12.2.4)(iv)]{EGA-IV.3}, the set $U\subseteq Y$ where the fibres are geometrically
normal is open.  The previous argument shows that this set is nonempty, and by equivariance we conclude that
each fibre of $\phi$ is geometrically normal.  Since each fibre of $\phi$ is also connected, each fibre of 
$\phi$ (over $T_{Y}$) is therefore irreducible.

In particular, the fibre $\phi^{-1}(1_{Y})$ is irreducible.  Since the kernel of 
$\phi|_{T_{X}}\colon T_{X}\longrightarrow T_{Y}$ is a nonempty subset of this fibre, the kernel is also
irreducible, and hence connected. 
By \cite[p. 112, \S8.2, Proposition ({\em c})]{B} every algebraic subgroup of a torus is diagonalizable,
and using \cite[p.\ 114, \S8.5, Proposition]{B} again, a diagonalizable and connected linear algebraic group
is a torus.  Thus the kernel is a torus.
\epf


\bpoint{Torus linearizations}
\label{sec:TK-decomposition}
Recall that on a toric variety $X$ there is an exact sequence

\begin{equation}
\label{eqn:linearization-exact-sequence}
0\longrightarrow M \longrightarrow \Pic^{T_{X}}(X) \longrightarrow \Pic(X)\longrightarrow 0,
\end{equation}

where $\Pic^{T_{X}}(X)$ is the group of $T_{X}$-linearized line bundles on $X$, i.e., the group of pairs $(L,\Phi)$, 
with $L$ a line bundle on $X$ and $\Phi$ a lift of the $T_{X}$ action on $X$ to the total space of $L$.

As above, let 
$\phi\colon X\longrightarrow Y$ be a toric fibration, and $T_{K}$ the kernel of 
$\phi|_{T_{X}}\colon T_{X}\longrightarrow T_{Y}$.   
The morphism $\phi$ is $T_{K}$-equivariant, with $T_{K}$ acting trivially on $Y$. 
By Lemma \ref{lem:kernel-of-toric-fibration}, $T_{K}$ is again a torus, of dimension $\dim(T_{K})=\dim(X)-\dim(Y)$.  

Suppose that $L$ is a line bundle on $X$, with a choice of $T_{X}$-linearization. Since $T_{K}$ is a subgroup
of $T_{X}$, we have by restriction an action of $T_{K}$ on $L$. Fix $i\geq 0$.  Since
$\phi$ is $T_{K}$-equivariant, with $T_{K}$ acting trivially on $Y$, and since $L$ has a $T_{K}$-action, 
$T_{K}$ also acts on the sheaf $R^{i}\phi_{*}L$.
It follows that $R^{i}\phi_{*}L$ decomposes as a direct sum of $T_{K}$-eigensheaves~:

\begin{equation}
R^{i}\phi_{*}L = \opp_{\mu \in C(L,i)} \left(R^{i}\phi_{*}L\right)_{\mu},
\end{equation}

where $(R^{i}\phi_{*}L)_{\mu}$ denotes the subsheaf of $R^{i}\phi_{*}L$ on which $T_{K}$ acts via the 
character $\mu\in \chi(T_{K})$,
and $C(L,i)\subseteq \chi(T_{K})$ the set of characters such that $(R^{i}\phi_{*}L)_{\mu}\neq 0$.   
Since $\phi$ is proper, $R^{i}\phi_{*}L$ is a coherent sheaf on $Y$, and it follows that $C(L,i)$ is a finite set.

In \S\ref{sec:basic-construction}--\S\ref{sec:variation-on-the-construction} 
below we give, for each $\mu\in C(L,i)$, a construction of a finite complex of subsheaves of a line bundle on $Y$, 
depending on a pair $(D,E)$, with $D$ a Cartier divisor on $X$ and $E$ a Cartier divisor on $Y$,
such that the $i$-th cohomology of the complex maps naturally
to $\left(R^{i}\phi_{*}L\right)_{\mu}$.  
In Theorem \ref{thm:D-E-theorem} we show that one can choose $(D,E)$ so that this map is an isomorphism of sheaves.

Before giving the construction of this complex, we turn to another useful piece of information accompanying the
choice of a $T_{X}$-linearization.

\bpoint{Torus characters associated to Cartier divisors}
\label{sec:characters-associated-to-Cartier-divisors}
Let $L$ be a bundle on $X$, and $D$ a $T_X$-stable Cartier divisor on $X$, not
necessarily effective, such that as line bundles $\Osh_{X}(D)\cong L$.
To $D$ is a associated a rational section $s_{D}$ of $L$, well-defined up to nonzero scalar, 
such that $\div(s_{D})=D$.   

Assume that we have also chosen a $T_{X}$-linearization of $L$. 
Since $D$ is $T_X$-stable, under the $T_{X}$-linearization $s_{D}$ corresponds to a character $\mu'\in \chi(T_{X})=M$.
This is easiest to see when $D$ is effective, since then $s_{D}\in H^0(X,L)$, and $s_{D}$ is an eigenvector
for the induced $T_{X}$-action on $H^0(X,L)$.   In the case that $D$ is not effective, we can write $D=D_1-D_2$, with
each of $D_1$ and $D_2$ effective and $T_{X}$-fixed.  This gives corresponding line bundles $L_i=\Osh_{X}(D_i)$, $i=1$, $2$, such that $L=L_1\otimes L_2^{*}$.  Choose $T_{X}$-linearizations on $L_1$ and $L_2$ so that the induced
$T_{X}$-linearization on $L_1\otimes L_{2}^{*}$ is that of the original linearization of $L$.  Then $s_{D_1}$ and
$s_{D_2}$ are eigensections of $H^0(X,L_1)$ and $H^0(X,L_2)$, with respective $T_{X}$-characters $\mu'_1$ and $\mu'_2$,
and $s_{D}$ corresponds to the $T_{X}$-character $\mu':=\mu'_1-\mu'_2$.  This difference is independent of the
$T_{X}$-linearizations chosen on the $L_i$, as long as the condition on $L_1\otimes L_2^{*}$ is met. 
By restricting $\mu'$ to $T_{K}$, we obtain a $T_{K}$-character $\mu$.

As in \cite[Chap. I]{TEI} for each $\alpha\in M$ we use $\charx^{\alpha}$ for the corresponding
homomorphism $\charx^{\alpha}\colon T_{X}\longrightarrow \CC^{*}$, which, via 
$\Gamma(T_{X},\Osh_{T_{X}})\subset \CC(T_{X})=\CC(X)$ we also consider to be a rational function on $X$.

Given a $D$ as above, 
for any $\alpha\in M$ the Cartier divisor $D+\div(\charx^{\alpha})$ is $T_{X}$-stable, associated to the same line 
bundle $L$, and has $T_{X}$-character $\mu'+\alpha$.  
Conversely, if $D'$ is a $T_{X}$-stable Cartier divisor such that $\Osh_{X}(D')=L$, and such that $s_{D'}$ has
$T_{X}$-character $\mu'+\alpha$, then $D'=D+\div(\charx^{\alpha})$.
Thus, once a $T_X$-linearization on $L$ has been fixed, 
the above procedure gives a one-to-one correspondence between
$T_X$-fixed Cartier divisors whose corresponding line bundle is $L$, and characters of $T_{X}$. 
We will use this in the proof of Theorem \ref{thm:D-E-theorem}.

\bpoint{Basic Construction of the complex} 
\label{sec:basic-construction} 
Let $\phi\colon X\longrightarrow Y$ be a toric fibration, $L$ a $T_{X}$-linearized line bundle
on $X$ and, as above, let $D$ be a $T_{X}$-stable Cartier divisor such that $\Osh_{X}(D)\cong L$, with corresponding
$T_{K}$-character $\mu$.

Let $P$, $\sigma\to U_{\sigma}$ be the reduced \v{C}ech complex of $X$ given by 
Theorem \ref{thm:construction-of-toric-reduced-cover}.   For each $\sigma\in P$ we define an ideal sheaf
$\Ish_{D,\sigma}$ on $Y$ by setting, for each open $V\subseteq Y$, 

\vspace{-0.5cm}

\begin{equation}
\label{eqn:sigma-ideal-def}
\rule{0.75cm}{0cm}
\raisebox{-0.5cm}{
\begin{array}[t]{rcl}
\Ish_{D,\sigma}(V) & = & \left\{ f\in \Osh_{Y}(V) \st 
\phi^{*}(f) s_{D}|_{U_{\sigma}\cap \phi^{-1}(V)}\, \mbox{ is a regular section of $L|_{U_{\sigma}\cap \phi^{-1}(V)}$}
\right\} \\[3mm]
& = & \left\{f\in \Osh_{Y}(V) \st (D+\div(\phi^{*}(f)))|_{U_{\sigma}\cap \phi^{-1}(V)} \geq 0\right\}. \\
\end{array}}
\end{equation}


We define the complex $C^{\smallbullet}_{P}(\Ish_{D})$ (which, despite notationally looking very similar to 
$\CComplex_{P}$, is a complex of ideal sheaves on $Y$) by setting, for each $i\geq 0$, 
$$C^{i}_{P}(\Ish_{D}) := \prod_{\sigma\in P,\, |\sigma|=i} \Ish_{D,\sigma},$$
and with the differentials $d^{i}$ obtained by combining the sheaf restriction maps and 
incidence function $\ep(\sigma,\sigma')$ as in \S\ref{sec:Def-of-reduced-complex}.

We note that, unlike the notation $\CComplex_{P}$, where we are using a fixed sheaf $\Fsh$, the ideal sheaf 
$\Ish_{D,\sigma}$ is different for each $\sigma\in P$.

By virtue of the definition of $\Ish_{D,\sigma}$, for each $\sigma\in P$ and each open $V\subseteq Y$ we have a natural 
inclusion of abelian groups~:

\begin{equation}
\label{eqn:I_sigma-natural-map} 
\begin{array}[t]{rcl}
\Ish_{D,\sigma}(V) & \hookrightarrow & \Gamma(U_{\sigma}\cap \phi^{-1}(V), L) \\[2mm]
 f & \xmapsto{\rule{0.5cm}{0cm}} & \phi^{*}(f)\cdot s_{D}|_{U_{\sigma}\cap \phi^{-1}(V)}  \\
 \end{array}
\end{equation}

Since $T_{K}$ acts trivially on $Y$, and via the character $\mu$ on $s_{D}$, the image of this map lies
in $\Gamma(U_{\sigma}\cap \phi^{-1}(V),L)_{\mu}$ (a notation which makes sense, since $U_{\sigma}$ and
$\phi^{-1}(V)$ are both stable under the action of $T_{K}$, and thus
$\Gamma(U_{\sigma}\cap \phi^{-1}(V),L)$ also has a $T_{K}$ action).

For fixed $V$, the groups $\Gamma(U_{\sigma}\cap \phi^{-1}(V), L)$, $\sigma\in P$ themselves form a complex. 
This is the complex $\CComplex_{P}(L|_{\phi^{-1}(V)})$ of \S\ref{sec:Def-of-reduced-complex} attached to the
topological space $\phi^{-1}(V)$, the cell complex $P$, the assignment $\sigma\mapsto U_{\sigma}\cap V$, and
the sheaf $L|_{\phi^{-1}(V)}$.

Although this assignment satisfies (P1)--(P4) (see e.g., Corollary \ref{cor:reduced-cech-complex-on-a-subspace}),
it is not necessarily true that $\CComplex_{P}(L|_{\phi^{-1}(V)})$ computes the 
cohomology $H^{\smallbullet}(\phi^{-1}(V),L|_{\phi^{-1}(V)})$.
This is because, since $V\subseteq Y$ may be an arbitrary open set, the sets $U_{\sigma}\cap \phi^{-1}(V)$ do not have 
to be affine, and therefore the line bundle $L$ may not be in the category $\Ab_{P}$ (relative to this cover)
needed to apply Theorem \ref{thm:reduced-cech-computes-cohomology}.  

But, by Theorem \ref{thm:reduced-cech-spectral-sequence} there is a spectral sequence, whose $E_{2}^{i,0}$ term
is $H^{i}\left(\CComplex_{P}(L|_{\phi^{-1}(V)})\right)$, and which abuts to 
$H^{\smallbullet}(\phi^{-1}(V),L|_{\phi^{-1}(V)})$.  The edge morphism of this spectral sequence gives, for
each $i\geq 0$, a natural map

\begin{equation}
\label{eqn:edge-morphism-map}
H^{i}\left(\CComplex_{P}(L|_{\phi^{-1}(V)})\right) \longrightarrow H^{i}(\phi^{-1}(V),L|_{\phi^{-1}(V)}).
\end{equation}

For each open $V\subseteq Y$, evaluating the complex of ideal sheaves $C^{\smallbullet}_{P}(\Ish_{D})$ on $V$
gives a complex of abelian groups, whose $i$-th cohomology we denote by $H^{i}(C^{\smallbullet}_{P}(\Ish_{D})(V))$.
For each $i$, the assignment $V\mapsto H^{i}(C^{\smallbullet}_{P}(\Ish_{D})(V))$ defines a presheaf of abelian groups 
on $Y$, and we use $\Hsh^{i}\left(C^{\smallbullet}_{P}(\Ish_{D})\right)$ for its sheafification. 

For each $V$, the maps \eqref{eqn:I_sigma-natural-map} give a map of complexes 
$C^{\smallbullet}_{P}(\Ish_{D})(V)\longrightarrow \CComplex_{P}(L|_{\phi^{-1}(V)})$.  Taking cohomology, and using 
\eqref{eqn:edge-morphism-map}, we therefore get a map of presheaves which for each $i$, and each $V$ is a map 

\begin{equation}
\label{eqn:map-of-presheaves}
H^{i}(C^{\smallbullet}_{P}(\Ish_{D})(V))\longrightarrow H^{i}(\phi^{-1}(V),L|_{\phi^{-1}(V)}).
\end{equation}

Since $R^{i}\phi_{*}L$ is the sheafification of the presheaf $V\mapsto H^{i}(\phi^{-1}(V),L|_{\phi^{-1}(V)})$, 
the maps \eqref{eqn:map-of-presheaves} sheafify to give a map of sheaves 
$\Hsh^{i}\left(C^{\smallbullet}_{P}(\Ish_{D})\right)\longrightarrow R^{i}\phi_{*}L$.
Moreover, since $D$ has $T_{K}$-character $\mu$, the image of this map lies in $(R^{i}\phi_{*}L)_{\mu}$.

This is the ``basic construction'' attached to the data above, including the choice of $D$.

\bpoint{Variation on the construction} 
\label{sec:variation-on-the-construction} 
Suppose, as in \S\ref{sec:basic-construction} that $\phi\colon X\longrightarrow Y$ is a toric fibration, and
$L$ a $T_{X}$-linearized line bundle on $X$.  In this variation we assume we are given 
a pair $(D,E)$, with $D$ a $T_{X}$-stable Cartier divisor on $X$, and $E$ a $T_Y$-stable Cartier divisor on $Y$.  
Setting $\Msh:=\Osh_{Y}(E)$, we additionally assume that $\Osh_{X}(D)\cong L\otimes\phi^{*}\Msh^{*}$, 
and that $D$ has $T_{K}$-character $\mu$.

To determine the $T_{K}$-character of $D$, we must give a $T_{X}$-action on $L\otimes \phi^{*}\Msh^{*}$, and thus 
a $T_{X}$-action on $\phi^{*}\Msh$.  We do this by pulling back any $T_{Y}$-linearization of $\Msh$, and letting
$T_{X}$ act through the quotient map.  Since $T_{K}$ is in the kernel of the quotient map, the $T_{K}$-character is
independent of the $T_{Y}$-linearization chosen.   We will use this independence implicitly below. 

Following the construction in \S\ref{sec:basic-construction}, we then obtain a map of sheaves
$$\Hsh^{i}\left(C^{\smallbullet}_{P}(\Ish_{D})\right)\longrightarrow 
R^{i}\phi_{*}(L\otimes\phi^{*}\Msh^{*}) = \left(R^{i}\phi_{*}L\right)\otimes\Msh^{*}. $$
Tensoring with $\Msh$, this becomes a map 
$\Hsh^{i}\left(C^{\smallbullet}_{P}(\Ish_{D})\right)\otimes\Msh \longrightarrow R^{i}\phi_{*}L.$
As above, the image of this map lies in the $\mu$-component of $R^{i}\phi_{*}L$, so that, for each $i\geq 0$, 
this variation really defines a map of sheaves

\begin{equation}
\label{eqn:D-E-sheaf-map}
\Hsh^{i}\left(C^{\smallbullet}_{P}(\Ish_{D})\right)\otimes\Msh \longrightarrow \left(R^{i}\phi_{*}L\right)_{\mu}.
\end{equation}

\bpoint{Theorem ($(D,E)$ theorem)} 
\label{thm:D-E-theorem}
Let $\phi\colon X\longrightarrow Y$ be a toric fibration between smooth proper
toric varieties, and $T_{K}$ the kernel of the induced map on tori.  Let $L$ be a $T_{X}$-linearized line bundle
on $X$, and fix an $i\geq 0$. As in \S\ref{sec:TK-decomposition} we have a decomposition 
$$
R^{i}\phi_{*}L = \opp_{\mu \in C(L,i)} \left(R^{i}\phi_{*}L\right)_{\mu},
$$
with $C(L,i)$ a finite set of $T_{K}$-characters.  Then, for each $\mu\in C(L,i)$, there exists a pair $(D,E)$ 
as in \S\ref{sec:variation-on-the-construction} 
such that the map \eqref{eqn:D-E-sheaf-map} is an isomorphism.
Moreover, if $(D,E)$ is such a pair, then for any effective $T_{Y}$-stable Cartier divisor $E'$ on $Y$, 
$(D-\phi^{*}(E'),E+E')$ is also such a pair.

\bpf
As in \S\ref{sec:basic-construction} let 
$P$, $\sigma\to U_{\sigma}$ be the reduced \v{C}ech complex of $X$ given by 
Theorem \ref{thm:construction-of-toric-reduced-cover};  in particular, $P$ is the cell-complex $X_{\geq 0}$.  
Let $V_1$,\ldots, $V_k$ be a minimal cover of $Y$ by torus-stable open affines.   Thus each $V_j$ corresponds
to a $T_Y$-fixed point of $Y$, or equivalently, a $0$-cell $v_j$ of the cell complex $Y_{\geq 0}$. 

For each $j$, let $(Q_j,\psi_j)$ be the reduced \v{C}ech cover produced by 
Theorem \ref{thm:construction-of-toric-reduced-cover} applied to the semi-proper toric variety $\phi^{-1}(V_j)$
(the variety is semi-proper by Lemma \ref{lem:semi-proper-facts}({\em a}))
\footnote{ As in Example \ref{ex:semi-proper-example} each $Q_j$ is the inverse image of $v_j$ under the map
$X_{\geq 0}\longrightarrow Y_{\geq 0}$.}.   

We fix $i$, and $\mu\in C(L,i)$ (i.e., such that $(R^{i}\phi_{*}L)_{\mu}\neq 0$).
The proof of the theorem proceeds in two steps.

{\sc Step I} : For each $j$, find a divisor $D_j$ as in the basic construction in \S\ref{sec:basic-construction},
with respect to the reduced \v{C}ech complex $(Q_j,\psi_j)$ on $\phi^{-1}(V_j)$, such that the map 
$$\Hsh^{i}\left(C^{\smallbullet}_{Q_j}(\Ish_{D_j})\right)\bigr|_{V_j}
\longrightarrow (R^{i}\phi_{*}L|_{V_j})_{\mu}|_{V_j}$$
is an isomorphism.

{\sc Step II} : Patch these local choices together to get a pair $(D,E)$ as in 
\S\ref{sec:variation-on-the-construction} such that the map \eqref{eqn:D-E-sheaf-map} is an isomorphism. 

{\em Construction in Step I}. Fix $j$, and the pair $(Q_j,\psi_j)$ giving a cover of $\phi^{-1}(V_j)$.  For
each $\sigma\in Q_j$, we use $U_{\sigma}$ for the open set of $\phi^{-1}(V_j)$ provided by the theorem. 
By Theorem \ref{thm:construction-of-toric-reduced-cover} the complex $\CComplex_{Q_j}(L)$ computes the cohomology
of $L$ on $\phi^{-1}(V_j)$.  In degrees $i-1$, $i$, and $i+1$ this complex is 

\begin{equation}
\label{eqn:total-Q-complex}
\begin{array}{ccccccccc}
&& \rule{1.70cm}{0cm} \und{i-1} 
&& \rule{1.20cm}{0cm} \und{i} 
&& \rule{1.70cm}{0cm} \und{i+1} \\[2mm]
\cdots 
& \longrightarrow & \displaystyle \opp_{\sigma\in Q_j, |\sigma|=i-1} \Gamma(U_{\sigma},L)
& \stackrel{d^{i-1}}{\longrightarrow} & \displaystyle \opp_{\sigma\in Q_j, |\sigma|=i} \Gamma(U_{\sigma},L)
& \stackrel{d^{i}}{\longrightarrow} & \displaystyle \opp_{\sigma\in Q_j, |\sigma|=i+1} \Gamma(U_{\sigma},L)
& \longrightarrow & \cdots \\
\end{array}.
\end{equation}

Each of the sets $U_{\sigma}$ is stable under the action of $T_{X}$, and thus each $\Gamma(U_{\sigma},L)$ has
a $T_X$-action.  Since the differentials are $T_{X}$-equivariant (each is induced via the natural restriction maps),
the complex above splits into a direct sum of complexes of abelian groups, one for each $T_{X}$-character 
$\alpha\in M$~:

\begin{equation-w-subscript}{\alpha}
\label{eqn:alpha-Q-complex}
\begin{array}{ccccccccc}
&& \rule{1.70cm}{0cm} \und{i-1} 
&& \rule{1.20cm}{0cm} \und{i} 
&& \rule{1.70cm}{0cm} \und{i+1} \\[2mm]
\cdots 
& \longrightarrow & \displaystyle \opp_{\sigma\in Q_j, |\sigma|=i-1} \Gamma(U_{\sigma},L)_{\alpha}
& \stackrel{d^{i-1}}{\longrightarrow} & \displaystyle \opp_{\sigma\in Q_j, |\sigma|=i} \Gamma(U_{\sigma},L)_{\alpha}
& \stackrel{d^{i}}{\longrightarrow} & \displaystyle \opp_{\sigma\in Q_j, |\sigma|=i+1} \Gamma(U_{\sigma},L)_{\alpha}
& \longrightarrow & \cdots \\
\end{array}.
\end{equation-w-subscript}

Because the splitting is as a direct sum of complexes, for any subset $S\subseteq M$, if one takes the direct sum 
of the $\mbox{\eqref{eqn:alpha-Q-complex}}_{\alpha}$ for all $\alpha\in S$, 
the inclusion of the resulting complex into \eqref{eqn:total-Q-complex} induces an injection at the level of 
cohomology. 

Since $(R^{i}\phi_{*}L)_{\mu}$ is a coherent sheaf, and $V_j$ affine, $(R^{i}\phi_{*}L)_{\mu}(V_j)\subseteq 
H^{i}(\phi^{-1}(V_i),L)$ is a finitely generated module over $\Gamma(V_j,\Osh_{Y})$.  
Let $m_1$,\ldots, $m_s$ be a finite set of $T_{X}$-eigenvectors (of $H^{i}(\phi^{-1}(V_i),L)$) which generate
$(R^{i}\phi_{*}L)_{\mu}(V_j)$, and for each $t$, $t=1$,\ldots, $s$, let $\mu'_t$ be the corresponding $T_X$-character.
Since these generators all have $T_{K}$-character $\mu$, for each $t$ we have $\mu'_t|_{T_K}=\mu$.

By the discussion in \S\ref{sec:characters-associated-to-Cartier-divisors}, for each $\mu_{t}'\in M$ 
there is a unique $T_X$-stable divisor $D_{j,t}$ such that $\Osh_{X}(D_{j,t})=L$ and such that
$D_{j,t}$ has $T_{X}$-character $\mu'_{t}$. 
To connect this with the complexes $\mbox{\eqref{eqn:alpha-Q-complex}}_{\alpha}$ above we note that, for any 
$\sigma\in Q_{j}$, if $D_{j,t}$ has no poles in $U_{\sigma}$ then the
the vector space $\Gamma(U_{\sigma},L)_{\mu'_{t}}$ is one dimensional,
spanned by $s_{D_{j,t}}|_{U_{\sigma}}$.  Otherwise, if $D_{j,t}$ does have poles in $U_{\sigma}$,  then
$\Gamma(U_{\sigma},L)_{\mu'_{t}}$ is zero.

Since each $V_j$ corresponds to a $T_Y$-fixed point of $Y$, for each $j$ we have 
$V_{j}=\Spec(\CC[u_1,\ldots, u_m])$, with $m=\dim(Y)$, and where the $u_{k}$ are eigenfunctions for $T_Y$.

Given any $t=2$,\ldots, $s$, $s_{D_{j,t}}/s_{D_{j,1}}$ is a rational section of $\Osh_{X}$ with associated $T_X$-stable
divisor $D_{j,t}-D_{j,1}$.  Therefore $D_{j,t}-D_{j,1}=\div(f_t)$ for some $T_X$-eigenfunction, $f_t\in \CC(X)$.   
But, both $D_{j,1}$ and $D_{j,t}$ have the same $T_K$-character $\mu$, and therefore $f_t$ has trivial $T_K$-character. 

The exact sequence
$$0\longrightarrow T_{K} \longrightarrow T_{X}\longrightarrow T_{Y}\longrightarrow 0$$
shows that the pullback map $\phi^{*}\colon \CC(T_{Y})\longrightarrow \CC(T_{X})$ induces an isomorphism

\begin{equation}
\label{eqn:T_K-invt-functions}
\CC(Y)=\CC(T_{Y})\stackrel{\sim}{\longrightarrow} \CC(T_{X})^{T_K}=\CC(X)^{T_{K}}.
\end{equation}

Thus, $f_t=\phi^{*}(u_1^{a_{1,t}}u_2^{a_{2,t}}\cdots u_m^{a_{m,t}})$ for unique $a_{\ell,t}\in \ZZ$, and so
$D_{j,t} = D_{1,t} + \div\bigl(\phi^{*}(u_1^{a_{1,t}}\cdots u_{m}^{a_{m,t}})\bigr)$, and
(up to scalar) $s_{D_{j,t}} = \phi^{*}(u_1^{a_{1,t}}\cdots u_{m}^{a_{m,t}})\cdot s_{D_{j,1}}$. 
For completeness we also define the exponents $(a_{1,t},\ldots, a_{m,t})$ in the case $t=1$ by setting 
$(a_{1,1},\ldots, a_{m,1}) = (0,\,0,\,\ldots, 0)$. 

For each $\ell=1$,\ldots, $m$, set $b_{\ell}=\min(a_{k,1},a_{i,2},\ldots, a_{k,s})$, i.e., set 
$b_{\ell}$ to be the minimum
of the exponents of $u_{\ell}$ appearing above as $t$ goes from $1$ to $s$.  We note that all $b_{\ell}\leq 0$. 
Finally, set $$D_j=D_{j,1}+\div\bigl(\phi^{*}(u_1^{b_1}\cdots u_{m}^{b_m})\bigr).$$

By construction, for each $t$,
$$D_{j,t} = D_{j}+\div\bigl(\phi^{*}(u_1^{a_{1,t}-b_1}u_2^{a_{2,t}-b_{t}}\cdots u_m^{a_{m,t}-b_m})\bigr),$$
and so also (up to scalar)

\begin{equation}
\label{eqn:s-D-j-span} 
s_{D_{j,t}} = \phi^{*}(u_1^{a_{1,t}-b_1}u_2^{a_{2,t}-b_{t}}\cdots u_m^{a_{m,t}-b_m})\cdot s_{D_j},
\end{equation}

with each exponent $a_{\ell,t}-b_{\ell}\geq 0$.   Since $D_j$ only differs from $D_{j,1}$ by the divisor of a rational
function pulled back from $Y$, $\Osh_{X}(D_j)\cong L$, and $D_j$ still has $T_{K}$-character $\mu$. 

We apply the basic construction of \S\ref{sec:basic-construction} to the morphism $\phi^{-1}(V_j)\longrightarrow V_j$
using the divisor $D_j$. For $\sigma\in Q_j$, and $V=V_j$, \eqref{eqn:sigma-ideal-def} is then
$$\Ish_{D_j,\sigma}(V_j) = 
\left\{f\in \Osh_{Y}(V_j) \st (D_j+\div(\phi^{*}(f)))|_{U_{\sigma}} \geq 0\right\}. \\
$$
Let $\beta_{\ell}$, $\ell=1$,\ldots, $m$ be the $T_{Y}$ characters corresponding to $u_1$,\ldots, $u_m$, 
$\beta'_{\ell}$, $\ell=1$,\ldots, $m$ the corresponding $T_{X}$ characters 
(i.e., composing with $T_{X}\longrightarrow T_{Y}$), and let $\mu'_j$ be the $T_X$-character of $D_j$.
Set

\begin{equation}
\label{eqn:S-j-set-def}
S_j:=\left\{ \mu'_j + \sum_{\ell=1}^{m} c_{\ell}\beta_{\ell} \st c_{\ell}\in \NN\right\}.
\end{equation}

The description of $I_{D_j,\sigma}$ above shows that 
image of the complex $\CComplex_{Q_j}(\Ish_{D_j})$ in $\CComplex_{Q_j}(L|_{\phi^{-1}(V)})$
under the injections \eqref{eqn:I_sigma-natural-map} is the sum, over all $\alpha\in S_j$ of the complexes
$\mbox{\eqref{eqn:alpha-Q-complex}}_{\alpha}$.

Thus, the corresponding map on cohomology, 

\begin{equation}
\label{eqn:Q_j-cohomology-map} 
H^{i}(C^{\smallbullet}_{Q_j}(\Ish_{D_j})(V_j))\longrightarrow H^{i}(\phi^{-1}(V_j),L|_{\phi^{-1}(V_j)})_{\mu} 
\end{equation}

is injective.   By construction, all of the weights $\mu_{t}'$ are in $S_j$, and thus the map
\eqref{eqn:Q_j-cohomology-map} above is also surjective, since its image contains all elements of 
$H^{i}(\phi^{-1}(V_j),L|_{\phi^{-1}(V_j)})_{\mu}$ of weights $\mu_{1}'$,\ldots, $\mu_{s}'$, and thus
contains all the generators $m_{1}$, \ldots, $m_{s}$ of $H^{i}(\phi^{-1}(V_j),L|_{\phi^{-1}(V_j)})_{\mu}$ 
chosen above. 

For a complex of coherent sheaves on an affine open set, computing cohomology sheaves commutes with the global
sections functor. 
Thus the isomorphism \eqref{eqn:Q_j-cohomology-map} is also an isomorphism 
$$\Hsh^{i}\left(C^{\smallbullet}_{Q_j}(\Ish_{D_j})\right)(V_j)
\stackrel{\sim}{\longrightarrow} (R^{i}\phi_{*}L|_{V_j})_{\mu}(V_j).$$
Finally, since both sides are coherent sheaves, and $V_j$ affine, we conclude that, on $V_j$, the natural map 

\begin{equation}
\label{eqn:sheaf-Q-j-isomorphism}
\Hsh^{i}\left(C^{\smallbullet}_{Q_j}(\Ish_{D_j})\right)\longrightarrow (R^{i}\phi_{*}L)_{\mu}
\end{equation}

from \S\ref{sec:basic-construction} is an isomorphism.

This completes step I.

{\em Construction in Step II}. We now assume we have been given $D_1$,\ldots, $D_k$, each $D_j$ a divisor on $X$,
such that the map \eqref{eqn:sheaf-Q-j-isomorphism} over $V_j$ is an isomorphism, as in step I.  
Let $F_1$,\ldots, $F_r$ be all of the $T_{Y}$-stable irreducible divisors on $Y$.
We apply the same idea as in step I, with the $F_1$,\ldots, $F_r$ replacing $\div(u_1)$,\ldots, $\div(u_m)$.  
We will use additive notation for the divisors. 

For each $j$, $j=1$,\ldots, $k$, $D_{j}-D_{1}$ is a $T_X$-stable Cartier divisor, whose associated
line bundle is $\Osh_{X}$, and which has $T_{K}$-character $\mu-\mu=0$.  Thus $D_{j}-D_{1}=\div(\phi^{*}f_j)$
for a $T_Y$-eigenfunction $f_j\in \CC(Y)$, and we may write 

\begin{equation}
\label{eqn:f-j-minimums}
\div(f_j) = \sum_{\ell=1}^{r} a_{\ell,j} F_{\ell}
\end{equation}

for unique $a_{\ell,j}\in \ZZ$.  When $j=1$ we have $(a_{1,1},\,\ldots,\, a_{r,1})=(0,\,\ldots,\,0)$.
As before, for each $\ell=1$,\ldots, $r$ we set $b_{\ell}=\min(a_{\ell,1},\,a_{\ell,2},\,\ldots,\, a_{\ell,k})$,
i.e., $b_{\ell}$ is the minimum of the coefficients of $F_{\ell}$ appearing above.  As before, each $b_{\ell}\leq 0$.
Setting 
$$
\begin{array}{rcl}
E & := &  -\sum_{\ell=1}^{r} b_{\ell} F_{\ell},\rule{0.25cm}{0cm}\mbox{and} \\[1mm]
D & : = & D_1-\phi^{*}(E) = D_1+\sum_{\ell=1}^{r} b_{\ell} \phi^{*}(F_{\ell}),
\end{array}
$$
then $E$ is an effective $T_Y$-stable divisor, $D$ a $T_{X}$-stable divisor, 
$\Osh_{X}(D)=L\otimes \phi^{*}(\Msh^{*})$, with $\Msh:=\Osh_{Y}(E)$, and 
$D$ has $T_{K}$-character $\mu$. 

By construction for each $D_j$

\begin{equation}
\label{eqn:D-j-to-D-1-relation}
D_j = D_1 + \sum_{\ell=1}^{r} a_{\ell,j}\phi^{*}F_{\ell} 
= D + \sum_{\ell=1}^{r} (a_{\ell,j}-b_{\ell}) \phi^{*}(F_{\ell}),
\end{equation}

with each $a_{\ell,j}-b_{\ell}\geq 0$.

Applying the variation of the basic construction from \S\ref{sec:variation-on-the-construction} we obtain a map 

\begin{equation}
\label{eqn:D-E-sheaf-map-proof}
\Hsh^{i}\left(C^{\smallbullet}_{P}(\Ish_{D})\right)\otimes\Msh \longrightarrow \left(R^{i}\phi_{*}L\right)_{\mu}.
\end{equation}

We now check that this map is an isomorphism of sheaves.  To do so it suffices to check that 
\eqref{eqn:D-E-sheaf-map-proof} is an isomorphism over each $V_j$ in the cover chosen above, 
and we may also check that the map is an isomorphism after tensoring with $\Msh^{*}$. 

So, fix such a $V_j$.  Since $V_j\cong \AA^{m}$, $\Msh|_{V_j}$ is trivial.  In particular there is a unique
$T_{Y}$-stable divisor $E_j$, supported on $Y\setminus V_j$, such that $\Osh_{Y}(E_j) = \Msh$.   Multiplying by the
corresponding section $s_{E_j}|_{V_j}$ thus gives an isomorphism 
$\Msh^{*}|_{V_j}{\longrightarrow} \Osh_{Y}|_{V_j}$.

Combining these steps, it suffices to check that the composite map 
$$
\Hsh^{i}\left(C^{\smallbullet}_{P}(\Ish_{D})\right)|_{V_j} \longrightarrow 
\left((R^{i}\phi_{*}L)\otimes \Msh^{*}\right)_{\mu}|_{V_j} \stackrel{\cdot s_{E_j}}{\longrightarrow} 
\left(R^{i}\phi_{*}L\right)_{\mu}|_{V_j}
$$
is an isomorphism.   Since the sheaves are coherent, and $V_j$ affine, it suffices to check that 
the corresponding map on global sections

\begin{equation}
\label{eqn:composite-over-V-global-sections}
H^{i}\left(C^{\smallbullet}_{P}(\Ish_{D})|_{V_j}\right)
\longrightarrow
H^{i}(\phi^{-1}(V_j),L|_{V_j})_{\mu}.
\end{equation}

is an isomorphism

As in step I, the reduced \v{C}ech complex given by $(P,\psi)$ (restricted to $\phi^{-1}(V_j)$) 
computing $H^{i}(\phi^{-1}(V_j),L)$ splits as a direct sum of weight spaces over all $\alpha\in M$~:

\begin{equation-w-subscript}{\alpha}
\label{eqn:alpha-P-complex}
\begin{array}{ccccccccc}
&& \rule{1.70cm}{0cm} \und{i-1} 
&& \rule{1.20cm}{0cm} \und{i} 
&& \rule{1.70cm}{0cm} \und{i+1} \\[2mm]
\cdots 
& \longrightarrow & \displaystyle \opp_{\sigma\in P, |\sigma|=i-1} \Gamma(U_{\sigma},L)_{\alpha}
& \stackrel{d^{i-1}}{\longrightarrow} & \displaystyle \opp_{\sigma\in P, |\sigma|=i} \Gamma(U_{\sigma},L)_{\alpha}
& \stackrel{d^{i}}{\longrightarrow} & \displaystyle \opp_{\sigma\in P, |\sigma|=i+1} \Gamma(U_{\sigma},L)_{\alpha}
& \longrightarrow & \cdots \\
\end{array}.
\end{equation-w-subscript}

For each $\sigma$, the map $\Ish_{\sigma,D}(V_j)\longrightarrow \Gamma(U_{\sigma},L)$ is the composition of the maps

\begin{equation}
\label{eqn:function-inclusion-maps}
f \mapsto \phi^{*}(f)\cdot s_{D}|_{\phi^{-1}(V_j)} \mapsto \phi^{*}(f) \cdot s_{D}\cdot s_{E_j}|_{\phi^{-1}(V_j)}
\end{equation}

(the multiplication by $s_{E_j}$ arises from the isomorphism $\Msh^{*}|_{V_j}\longrightarrow \Osh_{Y}|_{V_j}$). 

As before, this implies that the image of $C^{\smallbullet}_{P}(\Ish_{D})(V_j)$ under these maps is a sum
of the subcomplexes $\mbox{\eqref{eqn:alpha-P-complex}}_{\alpha}$.   Let $\mu'_j$ be the $T_{X}$-character of $D_j$
from step I.  We will show below that there are $e_1$,\ldots, $e_{m}\geq 0$ such that, setting

\begin{equation}
\label{eqn:S-set-def}
\overline{S}_j:=\left\{ \mu'_{j} + \sum_{\ell=1}^{m} (c_{\ell}-e_{\ell})\beta'_{\ell} \st c_{\ell}\in \NN\right\},
\end{equation}

the image of 
$C^{\smallbullet}_{P}(\Ish_{D})(V_j)$ under the maps \eqref{eqn:function-inclusion-maps}
is the sum of the $\mbox{\eqref{eqn:alpha-P-complex}}_{\alpha}$ for all $\alpha\in \overline{S}_j$. 
Here as in step I, $\beta'_1$,\ldots, $\beta'_m$ are the 
$T_X$-characters corresponding to the $T_Y$-eigenfunctions $u_1$,\ldots, $u_m$ 
generating $\Osh_{Y}(V_j)=\CC[u_1,\ldots, u_{m}]$.   

As in step I, we deduce from the direct sum decomposition that \eqref{eqn:composite-over-V-global-sections} 
is injective.  
The fact that $e_1$, \ldots, $e_{m}$ are all $\geq 0$ implies that $S_j\subseteq \overline{S}_j$, with $S_j$
as in \eqref{eqn:S-j-set-def}.
In step I we have chosen $S_j$ so that the weights in $S_j$ contain all the weights of the $\Osh_{Y}(V_j)$-generators 
of $H^{i}(\phi^{-1}(V_j),L)_{\mu}$. Thus, given that this description of $\overline{S}_j$ is correct, 
\eqref{eqn:composite-over-V-global-sections} is surjective. 

In summary, to prove that \eqref{eqn:composite-over-V-global-sections}, and therefore also that 
\eqref{eqn:D-E-sheaf-map-proof} is an isomorphism, it is sufficient to establish 
that the
the image of 
$C^{\smallbullet}_{P}(\Ish_{D})(V_j)$ under the maps \eqref{eqn:function-inclusion-maps}
is the sum of the $\mbox{\eqref{eqn:alpha-P-complex}}_{\alpha}$ for all $\alpha\in \overline{S}_j$ with
$\overline{S}_j$ given by \eqref{eqn:S-set-def}. 

Choose a $T_Y$-linearization of $\Msh$, and thus (via $T_{X}\longrightarrow T_{Y}$) a $T_{X}$-linearization
of $\phi^{*}\Msh$.  Since $\Osh_{X}(D) = L\otimes \phi^{*}(\Msh^{*})$, it then makes sense to ask for the 
$T_X$-character of $D$, which we denote by $\mu'_{\circ}$. 

By \eqref{eqn:D-j-to-D-1-relation} there is a $T_{Y}$-stable, effective divisor $E'_j$ on $Y$ (a nonnegative
sum of $F_1$,\ldots, $F_{r}$) such that $D_j = D + \phi^{*}(E'_j)$.  Since $\Osh_{X}(D_j)=L$ and 
$\Osh_{X}(D)=L\otimes\Msh^{*}$, and since $\phi$ is a fibration, we deduce that $\Osh_{Y}(E'_j)=\Msh$.

Since $E'_j$ is effective, its restriction to $V_j$ is a $T_Y$-stable, effective divisor.  Thus there are
$e_1$,\ldots, $e_m\geq 0$ such that $E'_j|_{V_j} = \div(u_1^{e_1}\cdots u_{m}^{e_m})|_{V_j}$.   Set
$E''_j=E'_j-\div(u_1^{e_1}\cdots u_m^{e_m})$.  Then $E''_j$ is a $T_{Y}$-stable divisor, supported on $Y\setminus V_j$,
and, since $E''_j$ is linearly equivalent to $E'_j$, $\Osh_{Y}(E''_j) = \Msh$.  As remarked 
after \eqref{eqn:D-E-sheaf-map-proof}, there is a unique
such divisor, and thus $E''_j = E_j$.  From $D_j=D + \phi^{*}(E'_j)$ we get that the $T_{X}$-character 
of $\phi^{*}(E'_j)$ is $\mu'_j-\mu'_{\circ}$.   
From $E'_j = \div(u_1^{e_1}\cdots u_{m}^{e_m})+E_j$ we get that the $T_{X}$-character of
$E_j$ is $\mu'_j-\mu'_{0}-\sum_{\ell=1}^{r} e_{\ell}\beta'_{\ell}$. 

Thus under the maps in \eqref{eqn:function-inclusion-maps}, starting with a $T_{Y}$-eigenfunction 
$f=u_1^{c_1}\cdots u_{m}^{c_m}\in \Osh_{Y}(V_j)$, 
$\phi^{*}(f)\cdot s_{D}\cdot s_{E_j}$ has $T_{X}$-character
$$\left(\sum_{\ell=1}^{r} c_{\ell}\beta'_{\ell}\right) + \mu'_{\circ} + 
\left( \mu'_j-\mu'_{0}-\sum_{\ell=1}^{r} e_{\ell}\beta'_{\ell}\right) 
= 
\mu'_{j} + \sum_{\ell=1}^{m} (c_{\ell}-e_{\ell})\beta'_{\ell}.
$$

Since $f\in \Osh_{Y}(V_j) = \CC[u_1,\ldots, u_m]$, all $c_{\ell}$ are $\geq 0$. This completes the argument
for the description of $\overline{S}_j$, and so also the proof of the existence of a pair $(D,E)$ with the 
desired properties.

The final statement, that if $(D,E)$ is a pair (for a given $\mu$) then for any effective $T_{Y}$-stable divisor
$E'$ on $Y$, that $(D-\phi^{*}(E'),E+E')$ is also such a pair, is clear, since over each $V_j$ the effect
is just to enlarge the set $\overline{S}_j$ of $T_{X}$-characters appearing in the image. \epf

\point 
\label{sec:Remarks-on-proof-of-D-E-theorem}
{\bf Remarks:} 
(\Euler{1}) By Corollary \ref{cor:reduced-cech-complex-on-a-subspace},
taking the reduced \v{C}ech complex on $X$ given by $P$, and restricting it to $\phi^{-1}(V_j)$ gives
a reduced \v{C}ech complex which also computes cohomology of coherent sheaves.  This was used implicitly in
the proof above, to understand the map \eqref{eqn:composite-over-V-global-sections}.

(\Euler{2}) In step I, we could have used this restriction of the reduced complex coming from $P$, instead
of using the complex given by $(Q_j,\psi_{j})$ on each $\phi^{-1}(V_j)$.

The advantage of using the $Q_j$ appears when one wants to make this construction explicit, so that it may be made 
into an algorithm.  The $Q_j$ are sub cell-complexes of $P$, and so in general simpler.
Using the associated smaller \v{C}ech complex makes finding the weights of $\Osh_{Y}(V_j)$-module 
computationally easier.  

(\Euler{3}) The divisor $E$ constructed in the proof of Theorem \ref{thm:D-E-theorem} is effective.

\appendix

\section{From theorem to algorithm} 
\label{sec:appendix-algorithm}

\bpoint{Overview}
The purpose of this appendix is to provide an effective method for computing higher direct images of line bundles for toric morphisms. We accomplish this by leveraging the success of Theorem~\ref{thm:D-E-theorem}, which gives a broadly constructive approach to computing higher direct images. However, there are two ingredients which are not explicitly constructed: the first is a finite set of characters $C(L,i)$ (see paragraph \ref{sec:algorithm-description}) and the second is a collection of divisors associated to these characters. We fill these gaps by providing an algorithmic construction both of these and further adapt the construction in Theorem~\ref{thm:D-E-theorem} into an algorithm for computing higher direct images. The algorithms presented in this section are already implemented in the computer algebra software \textit{Macaulay2}; see \cite{Zo} for further examples and usage.

\point
Fix a toric fibration $\varphi \colon X \rightarrow Y$, that is, a torus-equivariant morphism such that $\varphi_* \Osh_X = \Osh_Y$.  an integer $i$, and let $L$ be a line bundle on $X$. The discussion in paragraph~\ref{sec:algorithm-description} guarantees the existence of a finite set of $T_K$-characters $C(L,i)$ and Theorem~\ref{thm:D-E-theorem} guarantees that for each $\mu \in C(L,i)$, pairs of torus-stable Cartier divisors $(D_\mu,E_\mu)$ so that the $i$-th cohomology of the complex of ideal sheaves $\bigoplus_{\mu \in C(L,i)} C^{\smallbullet}_{P}(\Ish_{D_\mu}) \otimes \Osh_Y(E_\mu)$ is isomorphic to the $i$-th higher direct image of $L$ under $\varphi$ (Theorem~\ref{thm:D-E-theorem}). We provide two algorithms. Algorithm~\ref{alg:TK-characters} computes the character set $C(L,i)$ and pairs of divisors $(D_\mu, E_\mu)$ for each $\mu \in C(L,i)$; and Algorithm~\ref{alg:cechcomplex} constructs a complex of ideals in the Cox ring of $Y$ whose sheafification yields the desired complex.

\bpoint{Notation} 
Denote the Cox ring of $X$ and $Y$ as being $S_X$ and $S_Y$ respectively. We use $M_K := \chi(T_K)$ and $\cdot|_{T_K} \colon M_X \rightarrow M_K$ to denote the projection map of character lattices. We often use $\rho_j$ to denote the $j$-th ray of the fan of $X$. Given any maximal cone $\sigma$ in the fan of $Y$, there is an embedding $\varphi_\sigma \colon M_Y \rightarrow M_X$, described as follows. We have a natural inclusion $\iota_\sigma \colon M_Y \rightarrow \Pic^{T_Y}(Y)$; see \cite[Equation (5.3.1)]{CLS}. The toric map $\varphi$ induces a morphism $\varphi \colon \Pic^{T_Y}(Y) \rightarrow \Pic^{T_X}(X)$ by pulling back divisors and so we obtain the commutative square
\begin{equation*}
\begin{tikzcd}
M_Y \arrow[r, "\iota_\sigma"] \arrow[d,"\varphi_\sigma" left] & \Pic^{T_Y}(Y) \arrow[d, "\varphi" right] \\ M_X \arrow[r, "\div"] & \Pic^{T_X}(X),
\end{tikzcd}
\end{equation*}
where $\div$ is the map from \eqref{eqn:linearization-exact-sequence}.

\bpoint{Negative sets}
\label{sec:negative-sets}
We identify the group of $T_X$-stable divisors $\Pic^{T_X}(X)$ with $\ZZ^s$, where $s$ is the number of variables of $S_X$, which is equivalent to fixing an ordering of the rays of the fan of $X$. As shown in \cite[Theorem~1.1]{EMS}, the cohomology of a line bundle corresponding to a divisor on a toric variety depends only on the negative components of the divisor. To this end, for any $T_X$-stable divisor $D \in \ZZ^s$, we define $\neg(D) \coloneqq \bigl\{j \in \{1,2,\ldots,s\} \colon D_j < 0 \bigr\}$. For any toric variety $X$ and integer $i$, the collection of \textit{negative sets} is 
\begin{equation*}
    \Neg^i(X) \coloneqq \bigl\{ \neg(D) \subset \{1,2,\ldots,s\} \colon \text{$D \in \Pic^{T_X}$ and $H^i(X,\Osh_X(D)) \neq 0$} \bigr\}.
\end{equation*}
For a toric morphism $\varphi \colon X \rightarrow Y$ and maximal cone $\tau$ in the fan of $Y$, we define the \textit{restricted negative sets} by further imposing the condition $\rho_j \subset \varphi^{-1}(\tau)$ for all $j \in J$. Notationally, we set
\begin{align*}
    \neg_\tau(D) &\coloneqq \bigl\{j \in \{1,2,\ldots,s\} \colon D_j < 0 \text{ and } \rho_j \subset \varphi^{-1}(\tau) \bigr\}, \\
    \Neg_\tau^i(X) &\coloneqq \bigl\{ J \in \Neg^i(X) \colon \text{$\rho_j \subset \varphi^{-1}(\tau)$ for each $j \in J$} \bigr\}.
\end{align*}
We have $H^i \bigl( \varphi^{-1}(V_\tau), \Osh_X(D) \bigr) \neq 0$ if and only if $\neg_\tau(D) \in \Neg_\tau^i(X)$. See Example~\ref{ex:F1-to-P2} for a demonstration of (restricted) negative sets.

\bpoint{Computing $C(L,i)$ and $(D,E)$}
\label{sec:appendix-alg1}
Algorithm~\ref{alg:TK-characters} computes the finite set $C(L,i)$ from \S\ref{sec:TK-decomposition} and, for each character $\mu \in C(L,i)$, the pair $(D_\mu, E_\mu)$ as guaranteed from Theorem~\ref{thm:D-E-theorem}.

\renewcommand{\figurename}{Algorithm}
\renewcommand{\thefigure}{A.5}
\begin{figure}[!ht]
\label{alg:TK-characters}
\caption{Algorithm for computing the set $C(L,i)$ and $(D,E)$.}
\framebox{
\vbox{\begin{tabbing}
Output: \= \kill
Input: \> An integer $i$, a toric fibration $\varphi \colon X \rightarrow Y$ of smooth projective toric \\
\> varieties, and a Cartier divisor $D \in \ZZ^{s}$ on $X$ such that $\Osh_X(D) \cong L$. \\
Output: \= The finite set $C(L,i)$ and pairs $(D_\mu, E_\mu)$ for each character $\mu \in C(L,i)$. \\[5pt]
iiii \= iiii \= iiii \= \kill
For each maximal cone $\sigma$ in the fan of $Y$ do \\
\> For each index set $J$ in $\Neg_\sigma^i(X)$ do \\
\> \> Make the polyhedron given by all $\mu \in M_X \otimes \RR$ such that $\neg_\sigma(\div \mu + D) = J$; \\
\> \> Compute the lattice points in the bounded part of this polyhedron. \\
\> Let $\Gamma_\sigma$ to be the union of all the lattice points above and set $C_\sigma = \{ m|_{T_K} \colon m \in \Gamma_\sigma\}$; \\
\> For each $\mu \in C_\sigma$ do \\
\> \> Set $D_{\sigma,\mu} = \iota_\sigma(f_{\min})$, where $f_{\min}$ is the componentwise minimum vector of \\
\> \> \> the set $\{f \in M_Y \mathbin{|} \varphi_\sigma f = v - \mu \text{ for some } v\in \Gamma_\sigma\} \cup \{0\}$; \\
Set $C(L,i) = \bigcup C_\sigma$; \\
For each $\mu \in C(L,i)$ do \\
\> Fix a maximal cone $\tau$, and set $-E_\mu$ to be the componentwise minimum vector of \\
\> \> the set $\{E \in \Pic^{T_Y}(Y) \mathbin{|} E = D_{\sigma, \mu} - D_{\tau, \mu} \text{ for some maximal cone } \sigma \neq \tau\} \cup \{0\}$; \\
\> Set $D_\mu = \div \mu + D + \varphi^*(D_{\tau, \mu} - E_\mu)$;
\end{tabbing}}}
\end{figure}

\begin{proof}[Proof of Correctness of Algorithm~\ref{alg:TK-characters}]
Fix $\sigma$ being a maximal cone of $Y$ and consider the affine open subset $V_\sigma$ of $Y$. Locally, we have $(R^i \varphi_* L) (V_\sigma) = H^i \bigl( \varphi^{-1}(V_\sigma),\mathcal{O}_X(D) \bigr)$. It follows that this module inherits the fine grading of the Cox ring $S_X$, where degrees are indexed by torus-stable Cartier divisors $D \in \Pic^{T_X}(X)$. For a torus character $m \in M_X$, we naturally have the torus-stable divisor $\div m$, and so we may examine the $(\div m)$-graded piece of $(R^i \varphi_* L) (V_\sigma)$. In particular, consider the set of all characters $m \in M_X$ such that $(R^i \varphi_* L) (V_\sigma)_{\div m} \neq 0$. Let $C_\sigma$ be the projection of this set via the map $\cdot |_{T_K}$. We have that $C_\sigma$ is contained in $C(L,i)$, because $C(L,i)$ indexes where the higher direct images are nonzero. Moreover, we necessarily have that $\bigcup_{\text{$\sigma$ maximal}} C_\sigma = C(L,i)$ since the $V_\sigma$ cover $Y$.

Therefore, we need to see that $C_\sigma$ as described in the algorithm is the same as the projection of those $m \in M_X$ such that $(R^i \varphi_* L) (V_\sigma)_{\div m} \neq 0$. Indeed, by \ref{sec:negative-sets}, we have that the $(\div m)$-graded piece is nonzero precisely when $\neg_\sigma(\div m) \in \Neg^i_\sigma(X)$. For any fixed negative set $J \in \Neg^i_\sigma(X)$, the condition $\neg_\sigma(\div m) = J$ is a collection of hyperplane inequalities and hence defines a polyhedron. The union of these polyhedra over all negative sets $J$ therefore projects to give us $C_\sigma$. To see why it suffices to take only the bounded part of these polyhedra, we note that every polyhedron can be decomposed as the Minkowski sum of a polytope and recession cone; see \cite[Corollary 7.1b]{S}. However, the recession cone necessarily projects to a finite set, since otherwise $C(L,i)$ would be infinite.

Now we construct the pair $(D_\mu, E_\mu)$. Following the proof of Theorem~\ref{thm:D-E-theorem}, for a fixed character $\mu \in C(L,i)$, there is a finite set of $T_X$-eigenvectors which generates $(R^i \varphi_* L)_\mu$. In fact, we explicitly have that these are contained amongst the characters $m \in \Gamma_\sigma$ satisfying $m|_{T_K} = \mu$ by the previous paragraph (though it is possible that some of the characters in $\Gamma_\sigma$ may not be generators of $(R^i \varphi_* L)$ as an $\Osh_Y(V_\sigma)$-module). Since these characters all satisfy $m|_{T_K} = \mu$, they are indexed by eigenfunctions in $M_Y$, and are in correspondence with the set
\begin{equation*}
    \{f \in M_Y \mathbin{|} \varphi_\sigma f = v - \mu \text{ for some } v\in \Gamma_\sigma\} \cup \{0\}.
\end{equation*}
The elements of this set are precisely the eigenfunctions $f_t$ after passing through the isomorphism in Equation~\eqref{eqn:T_K-invt-functions}. Hence, the componentwise minimum $f_{\min}$ has its components given by the same $b_\ell$ as in the proof of Theorem~\ref{thm:D-E-theorem}. The proof goes on to define the $T_X$-stable divisor $D_j$, but for the purpose of computing $E_\mu$, it is more convenient to use the $T_Y$-eigenfunctions, which is why there is a type-discrepancy. Indeed, the pairwise differences $\iota_\sigma(f_{\min}) - \iota_\tau(f_{\min}')$ are themselves the divisors of $T_Y$-characters, which match the divisor in Equation~\eqref{eqn:f-j-minimums}. Hence, the definition of $E_\mu$ in Algorithm~\ref{alg:TK-characters} matches that of the proof of Theorem~\ref{thm:D-E-theorem}, as well as that of $D_\mu$ with the additional detail that because we were working with $T_Y$-characters, we need to translate by the $T_K$-character $\mu$.
\end{proof}

\bpoint{Remark}
Each of the polytopes $\Gamma_\sigma$ are necessarily convex, but their union need not be, and so in general it is unclear if $C(L,i)$ is a convex set. If $i = 0$, it is immediately true that it is convex, since there is only one possible negative set $J = \varnothing$. However, when $i \neq 0$, we can find examples where $C(L,i)$ is not convex, even when $i$ is equal to the fibre-dimension and we might expect some type of relative duality to hold. See Example~\ref{ex:convexity-of-weights}.

\bpoint{Complex of Ideals to Complex of Sheaves}
\label{sec:appendix-ideals-to-sheaves}
An obstacle to implementation of the algorithm is presenting sheaves in a finite way. In particular, our goal is to present the complex of ideal sheaves $\bigoplus_{\mu \in C(L,i)} C^{\smallbullet}_{P}(\Ish_{D_\mu}) \otimes \Osh_Y(E_\mu)$ in a way which can be parsed by a computer. Fortunately, \cite[Proposition~3.1]{C} gives us that the sheafification of multigraded modules over the Cox ring is an exact functor. Consequently we construct a complex of ideals in the Cox ring whose sheafification is the desired complex. Due to exactness, the homology of this complex of ideals will be a module whose sheafification yields the higher direct image sheaf.

\bpoint{Computing Higher Direct Images}
\label{sec:appendix-alg2}
Algorithm~\ref{alg:cechcomplex} constructs the complex of ideals which sheafifies to $\bigoplus_{\mu \in C(L,i)} C^{\smallbullet}_{P}(\Ish_{D_\mu}) \otimes \Osh_Y(E_\mu)$. By taking cohomology of this complex, we obtain an algorithm for computing higher direct images of line bundles for toric fibrations.

\renewcommand{\thefigure}{A.8}
\begin{figure}[!ht]
\label{alg:cechcomplex}
\caption{Algorithm for computing the $i$-th higher direct images of $L$.}
\framebox{
\vbox{\begin{tabbing}
Output: \= \kill
Input: \> An integer $i$, a toric fibration $\varphi \colon X \rightarrow Y$ of smooth projective toric \\
\> varieties, and a Cartier divisor $D \in \ZZ^{s}$ on $X$ such that $\Osh_X(D) \cong L$ \\
Output: \= A presentation for an $S_Y$-module which sheafifies to $R^i \varphi_* L$. \\[5pt]
iiii \= iiii \= iiii \= \kill
Construct the reduced \v{C}ech complex $(P,\psi)$ of $X$ by resolving the irrelevant ideal; \\
For each $\mu \in C(L,i)$ do \\
\> For each term in $P$ indexed by a face $\sigma$ \\
\> \> Construct the polyhedron given by all $f \in \Pic^{T_Y}(Y) \otimes \RR_{\geq 0}$ such that the $i$-th \\
\> \> \> component of the vector $\varphi^* (f) + D_\mu$ is nonnegative for each $i \in \sigma$; \\
\> \> Set $Q_{s,\sigma}$ to be the list of lattice points of the bounded part of this polyhedron; \\
\> \> Form the ideal $I_{\mu, \sigma} = \langle \textbf{y}^m \in S_Y \mathbin{|} m \in Q_{s,\sigma} \rangle$; \\
Construct the complex $C^\bullet_P(\Ish_D)$ by replacing each term in the complex $P$ by the direct \\
\> sum $\bigoplus_{\mu \in \mu_{D,i}} I_{\mu,\sigma} \otimes \Osh_Y(-E_\mu)$ and the maps by the induced inclusions of ideals; \\
Return the $i$-th cohomology of $I$;
\end{tabbing}}}
\end{figure}

\begin{proof}[Proof of Correctness of Algorithm~\ref{alg:cechcomplex}]
Since the maps of the complex are given by a composition of an incidence function and the natural inclusions, we just need to make a choice of incidence function and compute the terms. The most convenient way to obtain a choice of incidence function is by resolving the irrelevant ideal of $X$, which has the advantage of being very computationally efficient in \textit{Macaulay2}. To see why the resolution of the irrelevant ideal provides an incidence function on $P_X$, see \cite[Proposition 4.5 and Exercise 4.5]{MS}. One notes that these address the case of a projective toric variety, though the argument generalizes verbatim to the complete case.

To obtain the terms of $\Ish_{D_\mu}$, we first observe that it suffices to find ideals in $S_Y$ which sheafify to the corresponding ideal sheaf as described in \S\ref{sec:basic-construction}. Fix a maximal cone $\sigma$. We claim that the ideal $I_{\mu, \sigma}$ in the algorithm sheafifies to the ideal sheaf $\Ish_{D,\sigma}$. The primary condition to verify is in the second line of Equation~\ref{eqn:sigma-ideal-def}. In particular, we need that $I_{\mu, \sigma}$ contains all elements $f \in S_Y$ such that $D_{\mu} + \div f \geq 0$. This is precisely reflected in the algorithm, by constructing the polyhedron of all such elements which are regular in $S_Y$. Since any polyhedron decomposes into a Minkowski sum of a polytope and recession cone, the polytope of this polyhedron consists of generators for the ideal. By arranging these into a complex as described in \S\ref{sec:basic-construction}, we obtain the required complex.
\end{proof}

\section{Examples}
\label{sec:appendix-examples}

In this appendix we give examples of the algorithm, as well as some of the nuances in the combinatorics of the characters $C(L,i)$. For a demonstration of these examples in \textit{Macaulay2}, see the documentation for the \texttt{ToricHigherDirectImages} package.

\bpoint{Example}
\label{ex:F1-to-P2}
Let $Y = \FF_1$, the first Hirzebruch surface, and $Z = \PP^2$. We interpret $Y$ as the blow-up of $Z$ at a torus-fixed point and hence we have a toric blowdown map $\pi \colon Y \rightarrow Z$. The matrix of rays of $Y$ is $\div = \begin{bsmallmatrix} 1 & 0 & -1 & 0 \\ 0 & 1 & 1 & -1 \end{bsmallmatrix}^\text{\sffamily T}$ and the matrix of rays of $Z$ is $\begin{bsmallmatrix} 1 & 0 & 0 \\ 0 & 1 & -1 \end{bsmallmatrix}^\text{\sffamily T}$. The morphism $\pi$ corresponds to a map between lattices $\overline{\pi} \colon N_Y \rightarrow N_{Z}$ given by $\overline{\pi} = \begin{bsmallmatrix} 0 & -1 \\ 1 & 0 \end{bsmallmatrix}$. Let $D = 5B$ be five copies of the $(-1)$-curve, which corresponds to the vector $(0,5,0,0) \in \Pic^{T_Y}(Y)$. We compute the first higher direct images of $L = \Osh_Y(D)$ by following Algorithms~\ref{alg:TK-characters} and \ref{alg:cechcomplex}.

Denote the Cox ring of $Z$ by $S = \CC[z_0, z_1, z_2]$ with the standard grading. In this case, the only relevant maximal cone corresponds to the open subset containing the blown up point, so fix this maximal cone $\sigma$. We begin by computing the the first negative sets of $Y$ via the \texttt{NormalToricVarieties} package in \textit{Macaulay2}, which gives us $\Neg^1(Y) = \{\{0, 2\}, \{1, 3\}\}$. The rays of $\Sigma_Y$ which map into $\sigma$ are $\rho_0, \rho_1, \rho_2$, leaving us with the restricted negative sets $\Neg_\sigma^i(Y) = \{\{0, 2\}\}$. 

Given a vector $\mu = (a,b) \in M_Y \otimes \RR$, we see that $\neg_\sigma(\div \mu + D) = \{0, 2\}$ if and only if
\begin{align*}
    a &< 0, & b+5 &\geq 0, & b-a &< 0.
\end{align*}
These hyperplane inequalities define a bounded polyhedron, the lattice points of which are $\Gamma_\sigma = \left\{ \begin{bsmallmatrix} -i \\ -j \end{bsmallmatrix} \;\middle\vert\; i,j \in \ZZ, 1 \leq i \leq j \leq 5 \right\}$.
Since the dimension of $Y$ and $Z$ are the same and $\pi$ is surjective, the kernel torus and its character lattice $M_K$ are zero-dimensional, so we know $C_\sigma = \{0\}$ and hence $C(L,i) = \{0\}$.

The inclusion $\iota_\sigma$ is given by $\iota_\sigma = \begin{bsmallmatrix} 1 & 0 \\ -1 & -1 \\ 0 & 1 \end{bsmallmatrix}$. Solving for $\pi_\sigma$ yields $\pi_\sigma = \begin{bsmallmatrix} 0 & 1 \\ 1 & 1 \end{bsmallmatrix}$. To compute $f_{\min}$, we solve $\pi_\sigma m = v$ for each $v \in \Gamma_\sigma$ and obtain the set $\left\{ \begin{bsmallmatrix} -i \\ -j \end{bsmallmatrix} \;\middle\vert\;  i,j \in \ZZ, 1 \leq i, j \leq 4 \right\}$. We get $f_\text{min} = \begin{bsmallmatrix} -4 \\ -4 \end{bsmallmatrix}$ by taking the minimum componentwise values. It follows that $D_{\sigma,0} = \iota_\sigma f_\text{min} = (-4, 8, -4)$. For the other maximal cones $\tau$ of $Z$, we have that $C_\tau$ is empty and so there is no corresponding $D_{\tau, 0}$.

Hence, we compute $-E_\mu = 0$ as the componentwise negative minimum of the set $\{0\}$ since there are no $D_{\tau,0}$ for $\tau \neq \sigma$. We have $\pi^* \colon \Pic^{T_Z}(Z) \rightarrow \Pic^{T_Y}(Y)$ is given by the matrix $\begin{bsmallmatrix} 0 & 0 & 1 \\ 1 & 0 & 1 \\ 1 & 0 & 0 \\ 0 & 1 & 0 \end{bsmallmatrix}$, so that $D_\mu = D + \pi^*(E_\mu + D_{\sigma,0}) = (-4, -3, -4, 8)$. This concludes Algorithm~\ref{alg:TK-characters}.

Proceeding with Algorithm~\ref{alg:cechcomplex}, we construct the reduced \v{C}ech complex $(P,\psi)$ given by
\begin{equation*}
\begin{tikzcd}[column sep = 20pt]
0 \arrow[r] & \begin{matrix} H^0(U_{01},L) \oplus H^0(U_{03},L) \\ \oplus \\ H^0(U_{12},L) \oplus H^0(U_{23},L) \end{matrix} \arrow[r, "d^0"] & \begin{matrix} H^0(U_0,L) \oplus H^0(U_1,L) \\ \oplus \\ H^0(U_2,L) \oplus H^0(U_3,L) \end{matrix} \arrow[r, "d^1"] & H^0(T_Y,L) \arrow[r] & 0,
\end{tikzcd}
\end{equation*}
with differentials
\begin{align*}
    d^0 = \begin{bmatrix} -1 & 0 & 1 & 0 \\ -1 & 1 & 0 & 0 \\ 0 & -1 & 0 & 1 \\ 0 & 0 & -1 & 1 \end{bmatrix}, && d^1 = \begin{bmatrix} 1 & -1 & -1 & 1 \end{bmatrix}.
\end{align*}
For a vector $f = (a,b,c) \in \Pic^{T_Z}(Z)$ (indexing a monomial $z_0^a z_1^b z_2^c$), we consider 
\begin{equation*}
    \pi^* (f) + D_\lambda = (c - 4, a + c - 3, a - 4, b + 8).
\end{equation*}
Each of the terms in $\check{C}(P,L)$ is indexed by a maximal cone of $\Sigma_Y$, and we obtain polyhedral inequalities from these cones. For instance, for the $U_{01}$ term, we obtain the conditions $c - 4 \geq 0, a + c - 3 \geq 0$. We interpret these conditions as degree conditions for a monomial ideal, which yields the ideal $I_{0,\{0,1\}} = \langle z_2^4 \rangle$. Here are tables listing all of the terms we obtain.
\begin{align*}
\begin{tabular}{|c|c|c|}
\hline
Cone & Conditions & Ideal \\
\hline
01 & $c-4 \geq 0, a+c-3 \geq 0$ & $\langle z_2^4 \rangle$ \\
\hline
03 & $c-4 \geq 0, b+8 \geq 0$ & $\langle z_2^4 \rangle$ \\
\hline
12 & $a+c-3 \geq 0, a-4 \geq 0$ & $\langle z_0^4 \rangle$ \\
\hline
23 & $a-4 \geq 0, b+8 \geq 0$ & $\langle z_0^4 \rangle$ \\
\hline
\end{tabular}
&&
\begin{tabular}{|c|c|c|}
\hline
Cone & Conditions & Ideal \\
\hline
0 & $c-4 \geq 0$ & $\langle z_2^4 \rangle$ \\
\hline
1 & $a+c-3 \geq 0$ & $\langle z_0, z_2 \rangle^3$ \\
\hline
2 & $a-4 \geq 0$ & $\langle z_0^4 \rangle$ \\
\hline
3 & $b+8 \geq 0$ & $\langle 1 \rangle$ \\
\hline
$\varnothing$ & None & $\langle 1 \rangle$ \\
\hline
\end{tabular}
\end{align*}

There are two terms where something interesting happens. The cone $1$ corresponds to the $(-1)$-curve, and in this case we have the condition $a + c - 3 \geq 0$, which gives us the ideal $I_{0, \{1\}} = \langle z_0^3, z_0^2 z_2, z_0 z_2^2, z_2^3 \rangle = \langle z_0, z_2 \rangle^3$. The other interesting term is $U_3$, where the condition is $b + 8 \geq 0$, which is vacuous. In this case, we get $I_{0, \{3\}} = \langle 1 \rangle = S$.

Replacing the terms of the reduced \v{C}ech yields the complex of $S$-ideals $\check{C}^\bullet_P(\Ish_D)$
\begin{equation*}
\begin{tikzcd}
0 \arrow[r] & \begin{matrix} \langle z_2^4 \rangle \oplus \langle z_2^4 \rangle \\ \oplus \\ \langle z_0^4 \rangle \oplus \langle z_0^4 \rangle \end{matrix} \arrow[r, "d^0"] & \begingroup \setlength\arraycolsep{1pt} \begin{matrix} \langle z_2^4 \rangle & \oplus & \langle z_0, z_2 \rangle^3 \\ ~ & \oplus & ~ \\ \langle z_0^4 \rangle & \oplus & S \end{matrix} \endgroup \arrow[r, "d^1"] & S \arrow[r] & 0.
\end{tikzcd}
\end{equation*}
Taking first cohomology of this complex, we obtain that 
\begin{equation*}
\operatorname{H}^1 \bigl( \check{C}^\bullet_P(\Ish_D) \bigr) = \frac{\langle z_0, z_2 \rangle^3}{\langle z_0^4 \rangle + \langle z_2^4 \rangle}
\end{equation*}
as an $S$\nobreakdash-module. Since we computed earlier that $E_\mu = 0$, there is no further twist, and so this is a module sheafifying to the higher direct image $R^1 \pi_* L$. Notably, this is a finitely-generated presentation of the module.

The following example demonstrates that the twist $E_\mu$ need not be zero.

\bpoint{Example}
\label{ex:F1toP1}
Let $Y = \FF_1$, the first Hirzebruch surface, and $W = \PP^1$. We interpret $Y$ as the projective bundle $\PP\bigl( \Osh_W \oplus \Osh_W(1) \bigr)$, and let $\pi$ be the natural projection to the base. This is a toric morphism, with the corresponding map between lattices $\overline{\pi} \colon N_Y \rightarrow N_W$ given by $\overline{\pi} = \begin{bsmallmatrix} 1 & 0 \end{bsmallmatrix}$. Let $D = -2E$ be $(-2)$ copies of the $(-1)$-curve, which corresponds to the vector $(0,-2,0,0) \in \Pic^{T_Y}(Y)$. We compute the first higher direct images of $L = \Osh_Y(D)$.

Denote the Cox ring of $W$ by $A = \CC[w_0, w_1]$ with the standard grading and the maximal cones of $W$ as $\{1\}$ and $\{-1\}$. In this case the restricted negative sets for both maximal cones of $W$ are $\Neg_\bullet^i(X) = \{\{1,3\}\}$. The rays which map into $\{1\}$ are $0,1$ and $3$, and the rays which map into $\{-1\}$ are $1,2$, and $3$. The polyhedra we get are defined by the following hyperplane inequalities:
\begin{align*}
    \neg_{\{1\}}(\div \mu + D) = \{1,3\} && \Leftrightarrow && a &\geq 0, & b-2 &< 0, & -b &< 0, \\
    \neg_{\{-1\}}(\div \mu + D) = \{1,3\} && \Leftrightarrow && b-2 &< 0, & b-a &\geq 0, & -b &< 0.
\end{align*}
In this case, both of these polyhedra are unbounded as $a$ is a free parameter, but the lattice points obtained from the bounded parts are
$\begin{bsmallmatrix} 0 \\ 1 \end{bsmallmatrix}$ and $\begin{bsmallmatrix} 1 \\ 1 \end{bsmallmatrix}$ respectively. Since the restriction map $\cdot |_{T_K}$ is given by $\begin{bsmallmatrix} 0 & 1 \end{bsmallmatrix}$ there is only one character $\mu = 1 \in C(L,1)$.

The following analysis is the same for both maximal cones up to sign. For $\sigma = \{\pm 1\}$, the natural inclusion $\iota_\sigma$ is given by $\begin{bsmallmatrix} \pm 1 \\ \mp 1 \end{bsmallmatrix}$ and $\pi_\sigma$ by $\begin{bsmallmatrix} \pm 1 \\ 0 \end{bsmallmatrix}$. To compute $f_{\min}$, we solve for the preimages of $v - \mu$ for the lattice points above $v = \begin{bsmallmatrix} 0 \\ 1 \end{bsmallmatrix}, \begin{bsmallmatrix} 1 \\ 1 \end{bsmallmatrix}$. This set of preimages is $\{0,\pm 1\} \subset \ZZ^1$. Since $f_{\min}$ is the (component-wise) minimum of this set,  for $\sigma = \{1\}$ we get $f_{\min} = 0$ and for $\sigma = \{-1\}$ we get $f_{\min} = -1$. It follows that $D_{\{1\}, \mu} = \iota_{\{1\}}(0) = (0, 0)$ and $D_{\{-1\}, \mu} = \iota_{\{-1\}}(-1) = (-1, 1)$. 

Fixing $\tau = \{1\}$, we get that  $-E_\mu = (-1,0)$ is the componentwise minimum of the set $\{(0,0), (-1, 1)\}$ and so $E_\mu = (1, 0)$. The map $\pi^* \colon \Pic^{T_W}(W) \rightarrow \Pic^{T_Y}(Y)$ is given by the matrix $\begin{bsmallmatrix} 0 & 1 \\  0 & 0 \\ 1 & 0 \\ 0 & 0 \end{bsmallmatrix}$ and it follows that $D_\mu = \div \mu + D + \pi^*(D_{\{1\},\mu} - E_\mu) = (0, -1, 0, -1)$. 

For a vector $s = (a,b) \in CDiv(W)$, we have 
\begin{equation*}
    \pi^* (s) + D_\mu = (a, -1, b, -1).
\end{equation*}
Again we carry out the construction of $\check{C}(P,L)$. Here are tables listing the terms obtained.
\begin{align*}
\begin{tabular}{|c|c|c|}
\hline
Cone & Conditions & Ideal \\
\hline
01 & $a \geq 0, -1 \geq 0$ & $0$ \\
\hline
03 & $a \geq 0, -1 \geq 0$ & $0$ \\
\hline
12 & $-1 \geq 0, b \geq 0$ & $0$ \\
\hline
23 & $b \geq 0, -1 \geq 0$ & $0$ \\
\hline
\end{tabular}
&&
\begin{tabular}{|c|c|c|}
\hline
Cone & Conditions & Ideal \\
\hline
0 & $a \geq 0$ & $\langle 1 \rangle$ \\
\hline
1 & $-1 \geq 0$ & $0$ \\
\hline
2 & $b \geq 0$ & $\langle 1 \rangle$ \\
\hline
3 & $-1 \geq 0$ & $0$ \\
\hline
$\varnothing$ & None & $\langle 1 \rangle$ \\
\hline
\end{tabular}
\end{align*}
Replacing the terms of the reduced \v{C}ech complex yields the complex of $S$-ideals
\begin{equation*}
\begin{tikzcd}
0 \arrow[r] & 0 \arrow[r, "d^0"] & S \oplus S \arrow[r, "d^1"] & S \arrow[r] & 0
\end{tikzcd}
\end{equation*}
with $d^1 = \begin{bsmallmatrix} 1 & -1 \end{bsmallmatrix}$. Taking first cohomology of this complex, we obtain that 
\begin{equation*}
\operatorname{H}^1(\Ish_P(D_\lambda)) = S
\end{equation*}
as an $S$\nobreakdash-module, which sheafifies to $\Osh_W$. However, we also need to twist by the bundle $\Osh_W(-E_\mu) \cong \Osh_W(1)$. As a result, we conclude that $R^1 \pi_* \Osh_Y(D) \cong \Osh_W(1)$.

The next example is tied to Example~\ref{sec:X-to-F1-example}, which demonstrates that $C(L,i)$ can have more than one element.

\bpoint{Example}
\label{ex:X-to-F1-characters}
Let $Y = \FF_1$ be the first Hirzebruch surface and set $X$ to be the blowup of $Y \times \PP^1$ at the torus-fixed curve $C = B \times [1:0]$, where $B$ is the unique irreducible curve in $Y$ of self-intersection $(-1)$. The matrix of rays of $X$ is
\begin{equation*}
    \div^\text{\sffamily T} = \begin{bmatrix}
        1 & 0 & -1 & 0 & 0 & 0 & 0 \\
        0 & 1 & 1 & -1 & 0 & 0 & 1 \\
        0 & 0 & 0 & 0 & 1 & -1 & 1
    \end{bmatrix}^\text{\sffamily T}.
\end{equation*}
There is a surjective map $\pi \colon X \rightarrow Y$ given by the composition of the blowdown with the projection to the first factor. As a map between lattices, $\pi$ is given by the matrix $\begin{bsmallmatrix} 1 & 0 & 0 \\ 0 & 1 & 0 \end{bsmallmatrix}$. Consider the divisor $D = -2E - 2H$, where $E$ is the exceptional divisor of the blow-up, and $H$ is the pullback of a hyperplane class in $\PP^1$. With the choice of rays for $X$ above, this divisor is indexed by the vector $(0,0,0,0,0,-2,-2)$ in $\Pic^{T_X}(X)$. We proceed with Algorithm~\ref{alg:TK-characters} to compute the set $C(L,i)$ for $L = \Osh_X(D)$ and $i = 1$.

Denote the Cox ring of $Y$ by $S_Y = \CC[w_0,w_1,w_2,w_3]$. Consider the maximal cone $\sigma$ spanned by the rays $\rho_1$ and $\rho_2$; these are the middle two columns in the matrix of rays of $\FF_1$ shown in Example~\ref{ex:F1-to-P2}. The rays of $X$ which map into $\sigma$ are $\rho_1, \rho_2, \rho_4, \rho_5$, and $\rho_6$. We obtain (via \textit{Macaulay2}) that the restricted negative sets are
\begin{equation*}
    \bigl\{ \{1, 4\}, \{4, 5\}, \{1, 4, 5\}, \{5, 6\}, \{4, 5, 6\} \bigr\}.
\end{equation*}
For each such negative set, we construct a polyhedron. Let $\mu = (a,b,c) \in M_X$ denote a character. For two of the negative sets, $\{1,4\}$ and $\{1,4,5\}$, this polyhedron is empty; the conditions $(\div \mu + D)_1 = b < 0$ and $(\div \mu + D)_4 = -b < 0$ contradict each other. For the remaining negative sets, we get a finite collection of lattice points. For instance, the negative set $\{5,6\}$ gives the conditions
\begin{align*}
    b \geq 0, && b-a \geq 0, && c \geq 0, && -c-2 < 0, && b+c-2 < 0.
\end{align*}
This is an unbounded polyhedron, whose bounded part has three lattice points $\begin{bsmallmatrix} 0 \\ 0 \\ 0 \end{bsmallmatrix}, \begin{bsmallmatrix} 0 \\ 0 \\ 1 \end{bsmallmatrix}$, and $\begin{bsmallmatrix} 1 \\ 1 \\ 0 \end{bsmallmatrix}$. Taking the lattice points from the other two negative sets as well, we obtain the collection
\begin{align*}
    \Gamma_\sigma = \left\{\begin{bsmallmatrix} 3 \\ 3 \\ -1 \end{bsmallmatrix}, \begin{bsmallmatrix} 0 \\ 0 \\ 0 \end{bsmallmatrix}, \begin{bsmallmatrix} 0 \\ 0 \\ 1 \end{bsmallmatrix},\begin{bsmallmatrix} 1 \\ 1 \\ 0 \end{bsmallmatrix}, \begin{bsmallmatrix} 0 \\ 0 \\ -1 \end{bsmallmatrix}, \begin{bsmallmatrix} 1 \\ 1 \\ -1 \end{bsmallmatrix}, \begin{bsmallmatrix} 2 \\ 2 \\ -1 \end{bsmallmatrix} \right\}.
\end{align*}
Since the restriction map $\cdot|_{T_K}$ is given by $\begin{bsmallmatrix} 0 \\ 0 \\ 1 \end{bsmallmatrix}$, we get three characters $\{-1,0,1\} \subset \ZZ^1$. In fact, if one performs the same computation for each maximal cone $\sigma$ in the fan of $\FF_1$, one finds that $C(L,i) = \{-1,0,1\}$.

\bpoint{Example} 
\label{ex:convexity-of-weights}
For $i=0$, $T_{K}$-characters appearing in $\phi_{*}L$ (i.e., in $C(L,0)$) are a convex set. On the other hand, for $i\neq 0$, $T_{K}$-characters appearing in $R^{i}\phi_{*}L$ do not have to be a convex set, even for the top dimension. For example, consider the toric surface $W$ which is the blow-up of $\PP^1 \times \PP^1$ at the four torus fixed points. This surface comes with a map $\pi_1 \colon W \rightarrow \PP^1$ which is the composition of the blowdown maps with the projection to the first factor. Set $D = -2E_1 - 2E_3 + \pi_1^{-1}[1:0]$, i.e. two of the exceptional divisors which don't intersect and the pullback of a point; see Figure~\ref{fig:nonconvex}.

\begin{center}
\begin{tabular}{c}
\begin{tikzpicture}[scale=0.75]
\draw plot [smooth] coordinates {(1.2,3) (0,2.8) (-1.2,3)};
\draw plot [smooth] coordinates {(0.9,3.1) (1.8,1.8) (3.1,0.9)};
\draw plot [smooth] coordinates {(3,1.2) (2.8,0) (3,-1.2)};
\draw plot [smooth] coordinates {(3.1,-0.9) (1.8,-1.8) (0.9,-3.1)};
\draw plot [smooth] coordinates {(1.2,-3) (0,-2.8) (-1.2,-3)};
\draw plot [smooth] coordinates {(-0.9,-3.1) (-1.8,-1.8) (-3.1,-0.9)};
\draw plot [smooth] coordinates {(-3,1.2) (-2.8,0) (-3,-1.2)};
\draw plot [smooth] coordinates {(-3.1,0.9) (-1.8,1.8) (-0.9,3.1)};

\node [below] at (0,2.8) {$\pi_2^{-1}[1:0]$};
\node [above right] at (1.8,1.8) {$E_1$};
\node [left] at (2.8,0) {$\pi_1^{-1}[1:0]$};
\node [below right] at (1.8,-1.8) {$E_2$};
\node [above] at (0,-2.8) {$\pi_2^{-1}[0:1]$};
\node [below left] at (-1.8,-1.8) {$E_3$};
\node [right] at (-2.8,0) {$\pi_1^{-1}[0:1]$};
\node [above left] at (-1.8,1.8) {$E_4$};
\end{tikzpicture}
\\
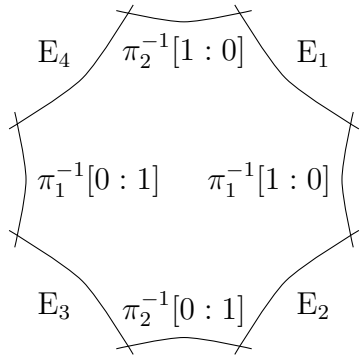
\Fig\label{fig:nonconvex} The boundary of $T_W$, which consists of a graph of $\PP^1$s.
\end{tabular}
\end{center}

Following the computation of Algorithm~\ref{alg:TK-characters}, one finds that for the maximal cone $\{1\}$, there is a single $T_K$ character $1$, and for the maximal cone $\{-1\}$ there is also a single $T_K$ character $-1$. Consequently $C(L,i) = \{-1, 1\}$, which is not a convex set. If we didn't include the pullback of the point, both characters would have been $-1$, which suggests that the behavior of the character sets can be independent amongst maximal affine opens.

\end{document}